\documentclass[12pt]{article}
\usepackage[a4paper,margin=2cm]{geometry}
\usepackage{amsfonts,graphicx,amsmath,amssymb,amsthm,color,nicefrac,enumerate, layout, bbm, listings,  hyperref}
\usepackage{cleveref} 
\usepackage{autonum} 
\usepackage{yfonts}
\usepackage[latin1]{inputenc}
\usepackage[numbers]{natbib}
\usepackage{tikz,caption}
\usetikzlibrary{patterns,shapes.symbols}
\usepackage[title]{appendix}
\usetikzlibrary{arrows.meta}
\numberwithin{equation}{section}
\newtheorem{proposition}{Proposition}
\newtheorem{lemma}{Lemma}
\newtheorem{theorem}{Theorem}

\theoremstyle{definition}
\newtheorem{definition}{Definition}
\theoremstyle{remark}
\newtheorem{example}{Example}
\newtheorem{hypothesis}{Hypothesis}
\newtheorem{remark}{Remark}

\newcommand{\convas}[1][]{\xrightarrow[#1]{\mathrm{a.s.}}}
\newcommand{\convp}[1][]{\xrightarrow[#1]{\PP}}
\newcommand{\convsl}[1][]{\xrightarrow[#1]{\mathcal{L}-\mathrm{s}}}
\newcommand{\convucp}[1][]{\xrightarrow[#1]{\mathrm{u.c.p.}}}
\newcommand{\NN}{\mathbb{N}}
\newcommand{\RR}{\mathbb{R}}
\newcommand{\PP}{\mathbb{P}}
\newcommand{\EE}{\mathbb{E}}
\newcommand{\cM}{\mathcal{M}}
\newcommand{\cQ}{\mathcal{Q}}
\newcommand{\cF}{\mathcal{F}}
\newcommand{\cG}{\mathcal{G}}
\newcommand{\ind}[1]{\mathbf{1}_{#1}}

\newcommand{\vd}{\,\mathrm{d}}
\newcommand{\Lone}{\mathcal{L}^1}
\newcommand{\Loneg}[1]{\mathcal{L}^1(\lambda^{(#1)})}
\newcommand{\Lgb}[1]{\mathcal{L}^{1}_{b}(\lambda^{(#1)})}
\newcommand{\Lg}[1]{\mathcal{I}_{#1}}
\newcommand{\bfF}[1]{\mathbf{E}^X_{#1}}
\newcommand{\bfH}[1]{\mathbf{E}^Y_{#1}}
\newcommand{\bfG}[1]{\mathbf{G}_{#1}}
\newcommand{\Lm}[3]{\mathbb{L}^{(#2)}(#1,#3)}
\DeclareMathOperator{\sgn}{sgn}
\def\eqlaw{\stackrel{\mathrm{law}}{=}}

\definecolor{rltgreen}{rgb}{0,0.5,0}
\definecolor{transpblue}{rgb}{0.6,0.8,0.95}
\definecolor{transpred}{rgb}{1,0.7,0.6}
\definecolor{trared}{rgb}{1,0.6,0.5}

\begin{document}

\title{Rates of convergence  
to the local time \\
of Oscillating and Skew Brownian Motions
}
\author{
Sara Mazzonetto\footnote{Universit\'e de Lorraine, CNRS, Inria, IECL, F-54000 Nancy, France; \texttt{sara.mazzonetto@univ-lorraine.fr}}
}

\date{09/01/2024}

\maketitle

\begin{abstract}
In this paper, a class of statistics based on high frequency observations of oscillating and skew Brownian motion is considered. 
Their convergence rate towards the local time of the underlying process is obtained in form of a functional limit theorem.
Oscillating 
and skew Brownian motion are solutions to stochastic differential equations with {\em singular\/} coefficients: piecewise constant diffusion coefficient or additive local time finite variation term. 
The result is applied to provide estimators of the skewness parameter and study their asymptotic behavior, 
and diffusion coefficient estimation is discussed as well.
Moreover, in the case of the 
classical 
statistics given by the
normalized number of crossings, 
the result is proved to hold for a larger class of It\^o processes with singular coefficients.\\
Up to our knowledge, 
this is the first result proving the convergence rates for estimators of the skewness parameter of skew Brownian motion.
\end{abstract}

\bigskip


\bigskip

\noindent{\textbf{AMS 2020: }} Primary 60F17, 
60J55;  
Secondary 
62F12, 
60F05. 


\section{Introduction}

It is well known that the normalized number of crossings of the level $r\in \RR$ of the time discretization (high frequency) of a Brownian motion (BM) provides an estimator for its local time at $r$. 
Roughly speaking the local time at the point $r$ 
measures the time the process spends around $r$ (see~\eqref{def:localtime} for a precise definition), so a rescaled number of crossings for high frequency data is a natural approximation of the local time also for more general processes. 
The normalized number of crossings has been extensively studied as an approximation of the local time of Brownian diffusions solutions of stochastic differential equations (SDEs) whose 
drift coefficient and diffusion coefficient $\sigma$ are sufficiently regular (in particular $\sigma$ is continuous). 
The convergence was proven, for instance, in~\cite{azais,azais2}. 
%
%

In this document we 
allow the presence of some singularities. We focus on one-dimensional It\^o processes solutions of stochastic differential equations (SDEs) with singular coefficients: 
discontinuous coefficients or drift in form of a weighted local time of the process at a given level. 
Solutions to such SDEs are often called {\em threshold\/} or {\em skew\/} diffusions. 
Two key cases belonging to this class of processes are oscillating Brownian motion (OBM) and skew Brownian motion (SBM).
They are generalizations of BM and of reflected BM as well, with distributions which are possibly singular with respect to BM. 
They change behavior when they reach a point, called barrier or {\em threshold\/}, which then becomes a {\em discontinuity point\/} for the local time $x\mapsto L^x_t(X)$ (see~\eqref{def:localtime}).
More precisely, OBM behaves like a BM with a different volatility above and below the threshold (causing a {\em regime-switch\/}) while SBM behaves like a BM everywhere except when it reaches the threshold, which plays the role of semi-permeable and semi-reflecting barrier. 
Note that OBM and SBM are {\em null-recurrent\/} processes, and so is the diffusion presenting both behaviors that we call oscillating-skew BM (OSBM).

More general functionals of discrete observations, which include the normalized number of crossings of a given level, are considered in this document. 
Given a stochastic process 
$(X_t)_{t\in [0,\infty)}$, let us consider the following statistics for high frequency observations:
\begin{equation} \label{eq:stat}
	\varepsilon_{n,t}^{(r,f,X)}:=\frac{1}{\sqrt{n}} \sum_{k=1}^{\lfloor n t \rfloor} f(\sqrt{n}(X_{(k-1)/n}-r),\sqrt{n}(X_{k/n}-r)) 
\end{equation}
where $f\colon \RR^2 \to \RR$ is a measurable function satisfying suitable integrability conditions.
The normalized number of crossings of the level $r$ corresponds to considering the function $f(x,y):= \ind{(-\infty,0)}(xy)$ (the statistic is explicitly provided in~\eqref{def:LT:discr:1}).

In the case $X$ is a BM, for the kind of estimators in \eqref{eq:stat}, convergence towards the local time and Central Limit Theorem (CLT) were obtained in~\cite{borodin86,borodin89}.
In the context of Brownian diffusions with regular coefficients mentioned above, convergence results are proven for specific functions $f$ in \cite{florens-zmirou} 
and for more general statistics of multivariate diffusions in \cite{j1}. In the latter article the associated CLT is proved.  
It is shown via semigroup estimates and martingale limit theorems,
that $\varepsilon_{n,t}^{(r,f,X)}$ behaves asymptotically (in $n$) as a mixed Gaussian distribution (see \eqref{intro:th} below).
The techniques developed so far can be adapted to study more general processes, for instance fractional BM in \cite{podolskij2017comment} and \cite[Corollary~14]{altmeyer2017estimating}.
In this document we adapt them to study another class of processes.

Let us have a look at the following simplified statement of our key result, that is Theorem~\ref{th:gen:jacod} in Section~\ref{sec:conv:rate}: 
let $X$ be an OSBM and $r$ be the threshold, let $(L^r_t(X))_{t\geq 0}$ be its symmetric local time at the threshold. 
Then for appropriate constants $c,K\in \RR$ it holds for all $t\in [0,\infty)$ that
\begin{equation}  \label{intro:th}
	n^{1/4} \left( \frac{ \varepsilon_{n,t}^{(r,f,X)} - c \, L^r_t(X)}{K^2 \sqrt{L^r_t(X)}}\right) \overset{n\to\infty}{\sim} \mathcal{N}(0,1).
\end{equation}
Although, at the threshold $r$, OSBM (as OBM and SBM) behave differently with respect to BM and in particular the local time $x \mapsto L^x_t(X)$ is discontinuous at $r$,
the speed of convergence is the same as for BM: $n^{1/4}$. 
\begin{remark}[Heuristics on the rate]
The convergence is different from what one might expect (i.e.~$n^{1/2}$), because the local time at $r$, and so its estimators, depends on the behavior of the process around $r$.
Indeed, as $n\to\infty$, among $n$ observations of the process on a fixed interval, the number of those which are sufficiently close to $r$ to matter is of order $n^{1/2}$.
Think of the random walk approximation of BM.
\end{remark}

The result in Theorem~\ref{th:gen:jacod} is actually more general. First, it is a {\em functional limit theorem\/}: the processes are seen as random variables with values in the Skorokhod space of c\`adl\`ag functions. Second, it holds also for 
processes combining pure skew and oscillating behavior and their drifted version,
under some suitable assumptions on the drift.
Therefore, this result extends the existing results on BM to solutions to SDEs with singular coefficients such as OBM, SBM, OSBM, and reflected BM (with suitable drift, by Girsanov's result). 
Extensions to more general statistics of 
more general singular diffusions 
is object of further research.
Nevertheless, in Theorem~\ref{th:crossings}, we prove the analogous of Theorem~\ref{th:gen:jacod} for the well known estimator of the normalized number of crossings when the process is a more general threshold diffusion. 

Let us illustrate, in the next section, few applications of Theorem~\ref{th:gen:jacod} on skewness parameter estimation for SBM. 
Note that, some among them are contributions of this document. 
In addition to those applications, in Section~\ref{sec:applications:base}, we discuss diffusion coefficients estimation for OSBM and propose, in Proposition~\ref{lem:estim:sigma}, estimators for the diffusions coefficients.

\subsection{Motivation and applications}

Since the seventies OBM and SBM together with their local time 
have been studied in the context 
of threshold diffusions.
In probability and stochastic analysis~\cite{barlow88skew, legall86yor, ramirezthesis, hajri15, forien19, taranto2021bi}, \ldots, recently in SPDEs \cite{BZ14, athreya2020well}, in simulation 
 \cite{EMloc, roberts2022skew}, \ldots, 
and we refer the reader to the introduction of~\cite{lejayPichot16} for some more applications of threshold diffusions in astrophysics, brain imaging, ecology, geophysics, fluid/gas dynamics, meteorology, molecular dynamics, oceanography. 

Some models in financial mathematics and econometrics are threshold diffusions, for instance continuous-time versions of SETAR ({\em self-exciting threshold auto-regressive}) models, see e.g.~\cite{decamps06, liptonsepp}. SBM and OBM and their local time have been recently investigated in the context of option pricing, as for instance in~\cite{GS17} and \cite{Cui21}. 
In \cite{lplaverage} it is shown that a time series of threshold diffusion type captures leverage and mean-reverting effects. We refer to the introduction of the latter article for further 
references. 

Statistical studies for threshold diffusions are partially motivated by calibration of such econometric models (see e.g.~\cite{SuChan15, lp}).
Indeed, study of (quasi) maximum likelihood estimators (MLE) of drift coefficients from high frequency observations 
depends on the approximations of occupation times 
and local times of the process. 
This is quite naturally explained by the fact that the behavior of the process changes at the threshold. 
Less heuristically, and more quantitatively, since 
a threshold diffusion 
behaves differently on two semi-axes, say $(-\infty,0)$ and $(0,+\infty)$, it is natural to look at the dynamics of the process in these semi-axes and this means considering $\max\{X_t,0\}$ and $\min\{X_t,0\}$. 
And if $X_t$ satisfies a SDE then by It\^o-Tanaka formula $\max\{X_t,0\}$ and $\min\{X_t,0\}$ satisfy a SDE involving the local time of the process $X_t$. (See~e.g.~\eqref{eq_Ito_Tanaka_plain} for It\^o-Tanaka formula.)
Theorem~\ref{th:gen:jacod} has been applied, in~\cite{mazzonetto2020drift}, to exhibit the asymptotic behavior in high frequency of (quasi) MLE of the drift parameters of a threshold diffusion which is a continuous-time SETAR model: a threshold Ornstein-Uhlenbeck process which follows two different Ornstein-Uhlenbeck dynamics above and below a fixed threshold. Similar applications are possible for other econometric models.

The latter application was the original motivation of this document. Nevertheless, Theorem~\ref{th:gen:jacod} contributes on one side to 
parameter estimation of SBM~answering to a conjecture in~\cite{lmt1} 
and on another side to estimation of the volatility jump of some levy-processes (see~\cite{robert_2023}, where a special case of Theorem~\ref{th:gen:jacod} is applied. 
Let us be more precise. 

Statistical analysis for SBM and OBM is quite recent: parameters estimation for skew reflected BM from continuous time observations was provided in~\cite{bardou2010statistical} and convergence properties (in long time) were studied by exploiting ergodic properties which SBM and OBM do not satisfy. 
Estimators based on high frequency observations of the skewness parameter of SBM are provided in~\cite{lmt,lmt1} and of the diffusion coefficients of OBM in~\cite{lp}. The estimators proposed are based on 
approximations of the local time and occupation times from high frequency observations of the process itself. 

Theorem~\ref{th:gen:jacod}
allows to establish 
the speed of convergence and limit distribution of several estimators of the skewness parameter from high frequency observation of OSBM.
In Theorem~\ref{th:skew:estim} we focus on SBM and introduce new estimators of the skewness parameter of SBM for which we show an asymptotic mixed normality property. 
This is the first result in the literature to prove the convergence speed for estimators of the skewness parameter of SBM.
The speed of convergence and the limit distribution for the MLE was conjectured in~\cite{lmt1}, based on the results for BM and illustrated using skew random walk in~\cite{lejay2018estimation}. 
In~\cite{lm23}, Theorem~\ref{th:gen:jacod} is applied to refine and prove the conjectured result.

In Section~\ref{sec:applications}, we apply Theorem~\ref{th:gen:jacod}, to two classical approximations of the local time of BM such as the normalized number of crossings.
Since standard BM (and reflected BM as well) is a special case of OSBM, as a by-product, we recover the classical results on the convergence (rates) for BM (see~\cite{borodin89,j1}).
In the case of the number of crossings,  
we can 
consider
threshold diffusions 
involving the local time and 
whose diffusion coefficient is piecewise differentiable and admits a finite jump 
(see Theorem~\ref{th:crossings}).
In a series of works Gikhman, Portenko, and Goshko study the convergence in law (without rate) 
towards the local time of the normalized number of crossings in a setting allowing for singularities (see~\cite{port1} and references therein) which seems to include OBM, SBM, and also sticky BM. 
Theorem~\ref{th:crossings} implies a stronger convergence (which was first proven in \cite{lmt1} for SBM) and exhibits also the rate of convergence. Up to our knowledge, Theorem~\ref{th:crossings} is the first result of its kind for such general one-dimensional threshold diffusions.

\subsection{Outline of the paper}
The paper is organized as follows.
First we introduce the processes OSBM as unique strong solutions to some SDE in Section~\ref{sec:main:processes}, 
then 
we state the main result: 
Theorem~\ref{th:gen:jacod} 
in Section~\ref{sec:conv:rate}. 
%
Section~\ref{sec:applications:base}
is devoted to applications: 
study of the number of crossings for a more general class of processes in Theorem~\ref{th:crossings}, estimation of the skewness parameter of SBM in Theorem~\ref{th:skew:estim}, and estimation of the parameters of OSBM in Proposition~\ref{lem:estim:sigma}.

Proofs are provided in Section~\ref{sec:proofs}. 
There it is shown that the proof of Theorem~\ref{th:gen:jacod} is based on an auxiliary proposition (Proposition~\ref{th:j2:2}), 
whose proof is so technical that it is provided in Appendix~\ref{sec:proof:j1}.
Appendix~\ref{sec:OBM:basic} deals with useful properties of OBM relevant in this article and in Appendix~\ref{sec:proof:j1}.
Appendix~\ref{app:prop:lmt:2} is an introduction to Appendix~\ref{sec:proof:j1}: some of the main ideas are already given in this section through a proof of the convergence (without rates) towards the local time of the statistics.

\subsection{Notation and notions of convergence} \label{sec:notion}

\subsubsection{Notation}
Throughout this document for every measurable functions $g\colon \RR \to \RR$ and measure $\mu$ on the Borel space $(\RR, \mathcal{B}(\RR))$ we denote by $\langle \mu , g \rangle$ the integral of $g$ with respect to the measure $\mu$:
\begin{equation}\langle \mu , g \rangle := \int_{-\infty}^{\infty} g(x) \mu(\!\vd x).\end{equation}

For every $\gamma \in [0,\infty)$ let 
$\lambda^{(\gamma)}$ be the measure on $(\RR, \mathcal{B}(\RR))$ absolutely continuous with respect to the Lebesgue measure satisfying 
${\lambda^{(\gamma)}(\vd x)} = |x|^\gamma \vd x$
and
let 
\begin{equation}( \Loneg{\gamma}, \|\cdot \|_{1,\gamma} )\end{equation} 
be the set of Borel measurable $\lambda^{(\gamma)}$-integrable functions and its norm.
If $\gamma=0$, we simply denote by $( \Lone, \|\cdot \|_1):=( \Loneg{0}, \|\cdot \|_{1,0})$ the normed space of Lebesgue integrable functions.

\begin{definition} \label{def:integral:gamma}
Let $\gamma \in [0,\infty)$. 
We denote by $\Lgb{\gamma}$ the following subspace of $\Lone$:
\begin{equation}
	\Lgb{\gamma} 
	= \{ f \colon \RR \to \RR , \text{ measurable and bounded} 
	\text{ s.t. } f\in \Loneg{\gamma}\}.
\end{equation}
We denote by $\Lg{\gamma}$ the following space of bi-variate measurable functions
\begin{equation}
	\Lg{\gamma} = \{ h \colon \RR^2 \to \RR , \exists \, \bar h \in \Lgb{\gamma}, \exists \, a\in [0,\infty) \text{ s.t. } \forall x,y\in \RR : |h(x,y)|\leq \bar h (x) e^{a |y-x|} \}.
\end{equation}
\end{definition}

\begin{remark} \label{rem:spaces:relation}
If $0 \leq \gamma_1 \leq \gamma_2 < +\infty$, then $\Lgb{\gamma_2} \subseteq  \Lgb{\gamma_1}$ and $\Lg{\gamma_2} \subseteq  \Lg{\gamma_1}$.
\end{remark}

\begin{remark} \label{rem:Lgb:remarks}
Let $\gamma \in [0,\infty)$ and 
let a function $\nu := \nu_+ \ind{[0,\infty)} + \nu_- \ind{(-\infty,0)}>0$.
Then for every functions $f\colon \RR \to \RR$ and $g\colon \RR^2 \to \RR$ it holds
\begin{itemize}
	\item $f\in \Lgb{\gamma}$ if and only if $f_\nu \in \Lgb{\gamma}$, where $f_\nu(x):= f(\nu(x) x)$,
	\item $h\in \Lg{\gamma}$ if and only if $h_\nu \in \Lg{\gamma}$, where $h_\nu(x,y):= h(\nu(x) x, \nu(y)y)$. 
\end{itemize}
\end{remark}
%
%
%
\begin{proof}[Proof of Remark~\ref{rem:Lgb:remarks}]
It suffices to prove one implication. Let us prove sufficiency for both statements.

Since $f$ and $\nu$ are measurable, then $f_\nu$ is measurable.
Boundedness for $f_\nu$ follows from boundedness of $f$.
A change of variable shows that there exists a constant $C\in (0,\infty)$ depending on $\nu_-,\nu_+$ such that
$
	\langle \lambda^{(\gamma)}, |f_\nu| \rangle 
	\leq C \langle \lambda^{(\gamma)}, |f| \rangle
$
which is finite by assumption.
Thus $f_\nu \in \Lgb{\gamma}$.

Since $h$ and $\nu$ are measurable, then $h_\nu$ is measurable.
Let us recall that there exist a non-negative function $\bar h \in \Lgb{\gamma}$ and a constant $a\in [0,+\infty)$ such that $|h(x,y)|\leq \bar h(x) e^{a|x-y|}$ for all $x,y\in \mathbb R$.
Hence, for all $x,y\in \mathbb R$ it holds that 
$
	|h_\nu(x,y)|\leq \bar h_\nu(x) e^{a|\nu(x) x-\nu(y)y|}
$ 
and it can be shown that
$
	|h_\nu(x,y)|\leq \bar h_\nu(x) e^{a \max\{\nu_-,\nu_+\}|x-y|}.
$ 
Since $\bar h\in \Lgb{\gamma}$, we have proven above that $\bar h_\nu \in \Lgb{\gamma}$, and $a \max\{\nu_-,\nu_+\} \in [0,\infty)$. 
Thus $h_\nu \in \Lg{\gamma}$.
\end{proof}

\subsubsection{The symmetric local time}
Let us give a more rigorous definition of the local time process.
Let $t\in [0,\infty)$ and let $(X_s)_{s\in [0,\infty)}$ be a one-dimensional semi-martingale. The {\em symmetric local time at the point $x$\/} accumulated on the time interval $[0,t]$ by the semi-martingale $X$ satisfies a.s.
\begin{equation} \label{def:localtime}
    L_t^x(X)=\lim_{\epsilon\to 0}\frac{1}{2\epsilon}\int_0^t
    \ind{\{-\epsilon\leq X_s-x\leq \epsilon\}}d\langle X\rangle_s
\end{equation}
and if $x=0$ we denote $L_t^0(X)$ by $L_t(X)$.

\subsubsection{Notions of convergence}
As already mentioned, the main aim of this article is studying, as $n\to\infty$, the convergence towards the local time together with its rate of the statistics $\varepsilon_{n,\cdot}^{(r,f,X)}$, with $X$ being an OSBM (see Section~\ref{sec:genera:proc}) and $f$ suitable function. 
Let us recall the notions of convergence used for the results of this paper.
The statement of the CLT involves the notion of {\em stable convergence\/} which was 
introduced and studied first in \cite{renyi} and \cite{AE78}.
We now specify it in the case used in this document.
\begin{definition}
Let $(D,d)$ be a metric space, $(\Omega^\prime ,\mathcal{F}^\prime,\mathbb{P}^\prime)$ be an extension of the probability space $(\Omega ,\mathcal{F},\mathbb{P})$, 
let $X_n \colon \Omega \to D$, $n\in \NN$, be a sequence of random variables, 
and let $X \colon \Omega^\prime \to D$ be a random variable. 
Then we say that $X_n$ converges stably in law to $X$
if for all $f\colon D \to \RR$ continuous and bounded and all bounded random variable $Y\colon \Omega \to \RR$ 
it holds that
\begin{equation}\lim_{n\to\infty} \EE\!\left[ f(X_n) Y\right] = \EE^\prime\!\left[ f(X) Y\right].\end{equation}
\end{definition}
Let $t \in [0,\infty)$, let $\mathbb{D}_t$, resp. $\mathbb{D}_\infty$, be the Skorokhod space of c\`adl\`ag functions from $[0,t]$, resp. $[0,\infty)$, to $\RR$ endowed with the Skorokhod topology.
When $D=\mathbb{D}_t$, $t\in [0,\infty]$ the functional stable convergence in law is usually denoted by
\begin{equation}X_n \convsl[n\to \infty] X.\end{equation}
Finally we recall the
notion used in the convergence results, i.e.~the {\em convergence in probability locally uniformly in time\/} or convergence {\em uniform on compacts in probability\/} (u.c.p.): 
let $X,X_n$, $n\in \NN$, be random variables with values in $\mathbb{D}_\infty$,
then
\begin{equation}X_n \convucp X,\end{equation} 
if for all $t\in [0,\infty)$ it holds that
$
	\sup_{s\in [0,t]} \left|  X_n(s) - X(s) \right| 
	\convp[n\to\infty] 0.
$

\section{Rates of convergence to the local time} \label{sec:main}

In the entire document
let $(\Omega,\mathcal{F}, (\mathcal{F}_t)_{t\in [0,\infty)},\PP)$ be a stochastic basis (for convenience, one may assume the usual conditions: completeness and right-continuity) 
and $W$ be an $(\mathcal{F}_t)_{t\in [0,\infty)}$-adapted standard BM.

In this section, we first introduce 
SBM as solution to the SDE~\eqref{eq:SBM} involving the local time of the process 
Then, OBM is specified as solution to the SDE~\eqref{eq:OBM} with discontinuous diffusion coefficient which is a special case of a wider class of processes extending both OBM and SBM, introduced in Section~\ref{sec:genera:proc}.
Finally we provide what can be considered the main result of this article, Theorem~\ref{th:gen:jacod}.
%
%

\subsection{The framework} \label{sec:main:processes}

\subsubsection{Skew Brownian motion}

Roughly speaking a SBM can be described trajectorially as 
a standard BM transformed by flipping its excursions from the origin with a certain probability.
In this document we refer to the characterization as unique strong solution to a SDE involving the local time, which was first considered by \cite{harrison-shepp_1981}. We refer the reader to a somehow recent survey paper on SBM~\cite{le}, where for instance, recalling Portenko's approach, it can be seen that the local time formally arises from a distributional drift: the Dirac $\delta$ at the threshold.\\

The SBM with skewness parameter 
$\beta\in [-1,1]$ at the threshold $r\in \RR$ is the diffusion which is strong solution to the following SDE
\begin{equation} \label{eq:SBM}
X_t=X_0 + W_t+ \beta L_t^{r}(X)
\end{equation}
where $L_t^{r}(X)$ is the symmetric local time of the process at $r$, $X_0\in \RR$, and $\beta X_0 \geq 0$ if $|\beta|=1$.
Some properties of the local time of SBM are object of the recent paper~\cite{BS19}.
\\
We call {\em standard SBM\/} a SBM with threshold $r=0$ starting at $X_0=0$. 
In this paper a SBM with skewness parameter $\beta\in (-1,1)$ is also denoted by $\beta$-SBM.
Note that a $0$-SBM is a BM. 
Moreover the $\pm 1$-SBM is a positively/negatively reflected BM.

The following quantities are important for the next sections.
Let us denote by $\mu_\beta$ the speed measure associated to the $\beta$-SBM with threshold $r=0$, that is
\begin{equation}
	\mu_\beta(\!\vd x) := \left( (1+\beta) \ind{(0,\infty)}(x) + (1-\beta) \ind{(-\infty,0)}(x) \right)\! \vd x =( 1+\sgn(x) \beta )  \vd x,
\end{equation}
and $p_\beta(t,x,y)$ denotes its transition density (first computed in \cite{walsh}): 
\begin{equation} \label{eq:SBM:density}
    p_\beta(t,x,y)
    =\frac{1}{\sqrt{2\pi t}}\exp\!\left(-\frac{(x-y)^2}{2t}\right)
    +\beta \sgn(y)\frac{1}{\sqrt{2\pi t}}\exp\!\left(-\frac{(|{x}|+|{y}|)^2}{2t}\right).
\end{equation}

\subsubsection{Oscillating-skew Brownian motion} \label{sec:genera:proc}

We now consider a generalisation of SBM, the oscillating-skew BM (OSBM) with threshold $r\in \RR$ which is solution to 
\begin{equation} \label{eq:OSBM}
	X_t = X_0 + \int_0^t  \sigma(X_s) \vd W_s + \beta L^r_t(X)  \quad t\geq 0 \quad \mathbb{P}\text{-a.s.}
\end{equation}
with $\beta \in [-1,1]$, 
the deterministic initial condition $X_0\beta\geq 0$ if $|\beta|=1$ and $X_0\in \RR$ otherwise, 
and $\sigma$ is the positive two-valued function discontinuous at the threshold  
\begin{equation}
    \label{sigmaOBM}
    \sigma := \sigma_- \ind{(-\infty,r)} + \sigma_+ \ind{[r,+\infty)}
\end{equation}
with $\sigma_-,\sigma_+\in (0,+\infty)$.
We refer to~\cite{legall} or \cite{BassChen05} for strong existence and pathwise uniqueness of solutions to~\eqref{eq:OSBM}.
In this document we also denote this process by $(\beta,\sigma_\pm)$-OSBM.
The speed measure of the OSBM with threshold $r=0$ (standard OSBM) is 
\begin{equation} \label{eq:speedm}
	\mu_X(\!\vd x) 
	= \frac{1+\sgn(x)\beta}{(\sigma(x))^2} \vd x
	= \left( \frac{1-\beta}{\sigma_-^2} \ind{(-\infty,0)}(x) + \frac{1+\beta}{\sigma_+^2} \ind{[0,\infty)}(x) \right)\! \vd x
\end{equation}
and (by It\^o-Tanaka formula in Lemma~\ref{lem:Ito-Tanaka}), the transition density function satisfies
\begin{equation} \label{eq:densities:X}
	p_X(t,x,y)=\frac{1}{\sigma(y)} \, p_{\beta_\sigma}\Big(t, \frac{x}{\sigma(x)}, \frac{y}{\sigma(y)}\Big) 
\quad \text{ with } \quad 
{\beta_\sigma}:= \frac{(1+\beta)\sigma_- -(1-\beta)\sigma_+}{(1+\beta)\sigma_-+(1-\beta)\sigma_+}
\end{equation}
where $p_{\beta_\sigma}$ is the density of the SBM recalled in~\eqref{eq:SBM:density}.

Let us specify the case of OBM which will be useful in the paper, 
although we could define it as a $(0,\sigma_\pm)$-OSBM.
\begin{remark}
When $\beta(\sigma_--\sigma_+)=0$ the $(\beta,\sigma_\pm)$-OSBM is either a SBM or an OBM. 
More precisely, a $(\beta,1)$-OSBM is a $\beta$-SBM and a $(0,\sigma_\pm)$-OSBM is a $\sigma_\pm$-OBM. 
\end{remark}

\subsubsection{Oscillating Brownian motion}

The OBM with threshold $r\in \RR$ solves the following special case of~\eqref{eq:OSBM}: 
\begin{equation}
    \label{eq:OBM}
    X_t=X_0+\int_0^t \sigma(X_s)\vd W_s, \quad t\geq 0,
\end{equation}
where $X_0\in \RR$ and the diffusion coefficient $\sigma$ is the function given in~\eqref{sigmaOBM}.
This process has been first defined and studied in~\cite{kw}. 
%
%
%

In this document, we also denote this process by $\sigma_\pm$-OBM and we call {\em standard OBM\/} an OBM $Y$ with threshold $r=0$ and starting point $Y_0=0$.
%
As for SBM and OSBM we can allow reflection. This would correspond to allow either $\sigma_-$ or $\sigma_+$ to be infinity: If $\sigma_+=1$, $\sigma_-=+\infty$, $Y_0 \geq r$ (resp. $\sigma_-=1$, $\sigma_+=+\infty$, $Y_0 \leq r$) then if $r=0$ it is a positively (resp. negatively) reflected BM.

We specify also the notation for~\eqref{eq:speedm} and~\eqref{eq:densities:X} for standard OBM (which corresponds to take $\beta=0$).
The speed measure and the transition density (first given in~\cite{kw}) are here denoted by
\begin{equation}	
\lambda_\sigma(\!\vd x) := \frac1{(\sigma(x))^2} \vd x 
\quad \text{ and } \quad
q_\sigma(t,x,y) =\frac{1}{\sigma(y)} \, p_{\beta_\sigma}\Big(t, \frac{x}{\sigma(x)}, \frac{y}{\sigma(y)}\Big)
\end{equation}
with ${\beta_\sigma}:= \frac{\sigma_- -\sigma_+}{\sigma_-+\sigma_+}$.

\subsection{Convergence towards the local time} \label{sec:main:SBM}


Let $X$ be the OSBM solution to~\eqref{eq:OSBM}, let $\mu_X$ denote the speed measure~\eqref{eq:speedm}, and let $p_X$ denote the transition density~\eqref{eq:densities:X}.
Note that 
$X$ is a one dimensional null-recurrent diffusion, therefore the speed measure $\mu_X$ is a stationary measure.
In particular one can take $X$ to be the $\sigma_\pm$-OBM solution to~\eqref{eq:OBM} or the $\beta$-SBM solution to~\eqref{eq:SBM}.

Given two measurable functions $f\colon \RR^2 \to \RR$ and $g\colon \RR \to \RR$, let
\begin{equation} \label{eq:muX}
	\bfF{f,g}(x) = \int_{-\infty}^{\infty} f(x,y) g(y) p_X(1,x,y) \vd y \quad
\text{ and } \quad \bfF{f} := \bfF{f,1}
\ (\text{i.e.~}g\equiv 1).
\end{equation}
We do not provide general assumptions for which the latter integrals are well defined. Assume for the moment that $f$ and $g$ are positives.
Note that, if $X$ denotes the solution with $r=0$, then
\begin{equation}
	\bfF{f,g}(x)=\EE\!\left[f(x,X_1) g(X_1) | X_0=x\right]
\end{equation}
and note that $\bfF{1,g}(x)$ is the transition semigroup (in usual notations, $P_1 g(x)$).

\begin{hypothesis} \label{hp:OSBM:1}
The measurable bi-variate 
function $f \colon \RR^2 \to \RR$ satisfies that
$\bfF{f}$, 
$\bfF{f^2} \in \Lgb{2}$
where $\Lgb{2}$ is given in~Definition~\ref{def:integral:gamma} in Section~\ref{sec:notion}.
\end{hypothesis}

\begin{remark} \label{rem:f_one_in_I}
With abuse of notation, if $f\colon \mathbb R \to \mathbb R$, we allow the slight abuse of notation 
$f \in \Lg{\gamma}$, $\bfF{f,g}$, and $\bfF{f}$ as above: one should interpret that $f$ coincides with the function 
$\left( (x,y) \mapsto f(x)\right)$.
Observe that, in this case, $\bfF{f}=f$.
\end{remark}

\begin{proposition}[Convergence towards the local time] \label{prop:lmt:2:OSBM}
Let $f \colon \RR^2 \to \RR$ be a function satisfying Hypothesis~\ref{hp:OSBM:1} and let $X$ be the OSBM solution to~\eqref{eq:OSBM}. 
Then 
\begin{equation} \label{prop:lmt:2:OSBM:eq}
	\varepsilon_{n,\cdot}^{(r,f,X)} \convucp[n\to\infty] \langle \mu_X , \bfF f \rangle L^r_{\cdot}(X).
\end{equation}
\end{proposition}
\noindent Note that the constant $\langle \mu_X , \bfF{f} \rangle$ can be rewritten as
\begin{equation}
	\langle \mu_X , \bfF{f} \rangle 
	= \EE\!\left[ f(X_0,X_1)|X_0 \sim \mu_X\right].
\end{equation}

\noindent Actually, the convergence in Proposition~\ref{prop:lmt:2:OSBM} is uniform in the parameters $\theta:= (\beta, \sigma_-,\sigma_+) \in [-1,1] \times (0,\infty) \times (0,\infty) :=\Theta$. 
More precisely, let $X^{(\theta)}$ denote the solution associated to the parameter $\theta \in \Theta$, then for all $t\in (0,\infty)$ it holds for all $\varepsilon \in (0,\infty)$ that 
\begin{equation} 
	\lim_{n\to\infty} \sup_{\theta \in \Theta} \mathbb{P}\left(  \sup_{s\in [0,t]} \left| \varepsilon_{n,s}^{(r,f,X^{(\theta)})} - \langle \mu_X , \bfF{f} \rangle L^r_{s}(X^{(\theta)})\right| \geq \varepsilon \right)
	= 0.
\end{equation}
In the case of SBM (with $\theta:=\beta$) the latter equation and Proposition~\ref{prop:lmt:2:OSBM} follow from \cite[Proposition~2]{lmt1} (with $T=1$) and the scaling property~\eqref{Yscaling} in Appendix~\ref{sec:scaling}.

\subsection{Rate of convergence to the local time} \label{sec:conv:rate}
We refine the above convergence showing that the speed of convergence is of order $1/4$.

\begin{theorem} \label{th:gen:jacod}
Let $f\in \Lg{\gamma}$, $\gamma > 3$, 
let $X$ be the $(\beta,\sigma_\pm)$-OSBM solution to~\eqref{eq:OSBM}.
Then there exists (possibly on an extension of the probability space) a BM $B$ independent of $X$ such that
\begin{equation} \label{th:gen:jacod:eq}
\begin{split}
	& n^{1/4} 
	\left(  \varepsilon_{n,\cdot}^{(r,f,X)} -  \langle \mu_X , \bfF{f} \rangle L^r_\cdot(X)\right)
	\convsl[n\to\infty]
\sqrt{K_f^X} B_{L^r_\cdot(X)}, 
\end{split} 
\end{equation}
where $K_f^X$ is the non-negative finite constant given by 
\begin{equation} \label{eq:k:th}
\begin{split}
   K_f^X
	& = \langle \mu_X , \bfF{f^2} + 2 \bfF{f,\mathcal{P}_{f}^X} \rangle  
			+ \frac{ 2 \sigma_- \sigma_+}{(1+\beta)\sigma_-+(1-\beta)\sigma_+} \frac{8}{3 \sqrt{2\pi}} (\langle \mu_X , \bfF{f} \rangle)^2 
			\\
	& \quad - 2 \sqrt{\frac{2}{\pi}} \frac{ 2 \sigma_- \sigma_+}{(1+\beta)\sigma_-+(1-\beta)\sigma_+} \langle \mu_X , \bfF{f} \rangle
			 \\
	& \qquad \cdot \int_{-\infty}^\infty \left( e^{-\frac{y^2}{2(\sigma(y))^2}} -\sqrt{2\pi} \frac{|y|}{\sigma(y)} \Phi\left(-\frac{|y|}{\sigma(y)}\right)\right) \mathcal{P}_{f}^X(y) \mu_X(\!\vd y)
	\\
	& \quad -
	2 \langle \mu_X , \bfF{f} \rangle 
	\left(\frac{ 2 \sigma_- \sigma_+}{(1+\beta)\sigma_-+(1-\beta)\sigma_+}\right)^{\!2}
	\\
	& \qquad \cdot 
\int_{-\infty}^{\infty} \int_0^1 \int_{-\infty}^{\infty} \tfrac{|x| e^{-\frac{x^2}{2(\sigma(x))^2}} \Phi\left(-\frac{|y|}{\sigma(y)}\right)}{\sqrt{2\pi}\sigma(x)}  \sqrt{\frac{1}t-1} f(x \sqrt{t}, y \sqrt{1-t}) \mu_X(\!\vd y) \vd t \mu_X(\!\vd x),
\end{split}\end{equation}
$\Phi$ is the cumulative distribution function of a standard Gaussian random variable, 
\begin{equation} \label{eq:cQ:limit}
	\mathcal{P}_{f}^X(x) = \sum_{j=0}^{\infty } \int_{-\infty}^{\infty} p_X(j,x, y)\left(\bfF{f}(y)- \langle \mu_X , \bfF{f} \rangle \bfF{g_\beta}(y)\right) \vd y,
\end{equation}
and 
$
	g_\beta(x,y) = \frac{1}{1+\sgn(y)\beta}\left(|y| - \frac{1+\sgn(y)\beta}{1+\sgn(x)\beta}|x|\right)
$.
\end{theorem}

\begin{remark}[Reflected BM]
	The assumptions of Proposition~\ref{prop:lmt:2:OSBM} and Theorem~\ref{th:gen:jacod} include reflected BM, i.e.~$|\beta|=1$. 
\end{remark}

\begin{remark}[About results for drifted OSBM] \label{rem:drift}
	Let us consider a drifted OBM solution of an SDE with drift coefficient given, for instance, by a bounded function $b$ of the process itself. Then, Proposition~\ref{prop:lmt:2:OSBM} and Theorem~\ref{th:gen:jacod} hold - with the same quantities and constants - for the drifted OBM. This follows by the same arguments in~\cite[Section~2.4]{j1}; to summarize, it suffices to prove the result on bounded time intervals $[0,t]$ (see the beginning of the proof if Theorem~\ref{th:gen:jacod}) and apply Girsanov's transform to pass from the driftless case to the drifted one. Since the drift is bounded, the Radon-Nikodym derivative is integrable. Then, by a dominated convergence argument, stable convergence holds for the drifted process.
The result for drifted OBM implies the one for drifted OSBM by a simple space transform as a consequence of Section~\ref{sec:interplay}.
\end{remark}

\begin{remark} \label{rem:th:BM}
If $\beta=0$ and $\sigma\equiv 1$ we recover the known result for BM: e.g.~\cite{borodin86,borodin89} and a special case of the already cited \cite[Theorem~1.2]{j1}.
The expression for the constant $K_f^X$ we propose is slightly more explicit. 
\end{remark}

\begin{remark} \label{comment:convergence}
Theorem~\ref{th:gen:jacod} implies a weaker version of Proposition~\ref{prop:lmt:2:OSBM}.
Proposition~\ref{prop:lmt:2:OSBM} requires Hypothesis~\ref{hp:OSBM:1} which is satisfied if for instance $f \in \Lg{2}$ (see Lemma~\ref{lem:link:theorems}). 
Theorem~\ref{th:gen:jacod} instead assumes $f\in \Lg{\gamma}$, $\gamma > 3$, which is a stronger condition.

Let us comment briefly on how to derive the u.c.p.~convergence from Theorem~\ref{th:gen:jacod}, 
although in Appendix~\ref{app:prop:lmt:2} we provide a direct proof of the weaker version of Proposition~\ref{prop:lmt:2:OSBM}.
The notions of convergence in law/stably in law/in probability coincide when the limit is constant and so  
$
	\varepsilon_{n,\cdot}^{(r,f,X)} -  \langle \mu_X , \bfF{f} \rangle L^r_\cdot(X) \convp[n\to\infty] 0
$
in the Skorokhod topology.
Since $L^r_\cdot(X)$ is (a.s.) continuous and increasing
it can be proven (splitting into positive and negative part of $f$, and so of $\bfF{f}$)
that 
$
	\varepsilon_{n,\cdot}^{(r,f,X)} \convucp[n\to\infty]  \langle \mu_X , \bfF{f} \rangle L^r_\cdot(X)
$
(see e.g.~(2.2.16) in \cite{jp}).
\end{remark}

\section{Applications} \label{sec:applications:base}

As mentioned in the introduction, several applications of Theorem~\ref{th:gen:jacod} and of the results of this section can be found in the literature (see for instance~\cite{lmt1,lmt,lm23,robert_2023}).
In this section, we first consider some classical estimators of the local time of BM and show that they are still estimators, up to a multiplicative constant, of the local time $L_T^r(X)$ of OSBM and a more general class of processes. Next, we provide estimators of the skewness parameter of SBM and analyze their asymptotic behavior.
This is the first time the convergence rate and an asymptotic mixed normality property for a skewness parameter of SBM is established and that joint estimation of OSBM is explored.

Let $r\in \RR$ be a fixed threshold, and let $X$ be a stochastic process. 
Let $T\in (0,\infty)$, $N\in \NN$, we observe the process on the discrete time grid $i \frac{T}{N}$. We denote by $X_i = X_{ i \frac{T}{N}}$.

\subsection{Estimating the local time via the number of crossings of the threshold }
\label{sec:applications}

Let $\alpha\in [0,\infty)$ and note that the function $h_\alpha$ given by $h_\alpha (x,y) = |y|^\alpha \ind{(-\infty,0)}(x y)$ is in $\Lg{\gamma}$ for all $\gamma\in [0,\infty)$. 
In fact $  h_\alpha(x,y) \leq c_\alpha e^{-|x|} e^{|y-x|}$ for some constant $c_\alpha$ depending on $\alpha$.
We consider two estimators obtained considering the functions proportional to $h_0$ and $h_1$:
\begin{align} 
& \label{def:LT:discr:1}
\mathcal{L}^r_{T,N}(X)
=  \varepsilon_{\frac{N}{T},T}^{(r,h_0,X)}
= \sqrt{ \frac{T}{N}}\sum_{i=0}^{N-1} \ind{\{(X_i-r)(X_{i+1}-r)<0\}}
\text{ and }
\\
& \label{def:LT:discr}
L^r_{T,N}(X)
= \varepsilon_{\frac{N}{T},T}^{(r, 2 h_1,X)}
= 2 \sum_{i=0}^{N -1} \ind{\{(X_i-r)(X_{i+1}-r)<0\}}|X_{i+1}-r|.
\end{align}
The first is concerned with the number of crossings of the threshold and the second takes into account the distance from it. Moreover note that $L^r_{T,N}$ requires only the knowledge of the $N+1$ observations $X_i$, $i=0,\ldots,N$, and not of $T/N$, nor of $T$.

As mentioned in the introduction, in the case of BM, and more general Brownian diffusions, these are consistent estimators of the local time up to a constant. 
We now show the consequence of Theorem~\ref{th:gen:jacod} for these estimators and, in the case of the statistic~\eqref{def:LT:discr:1}, we extend the result to skew diffusions with more general drift and diffusion coefficient. 


\subsubsection{Estimator counting the number of crossings of the threshold}

In this section, we consider the semi-martingale satisfying
\begin{equation} \label{eq:SDE:sigma}
	X_t = X_0 + \int_0^t b_s \vd s + \int_0^t  \sigma(X_s) \vd W_s + \beta L^r_t(X)  \quad t\geq 0 \quad \mathbb{P}\text{-a.s.}
\end{equation}
where $L^r_t(X)$ is the symmetric local time of the process at the fixed level $r\in \RR$, $X_0\in \RR$,
and $\beta$, $b$, and $\sigma$ satisfy Hypothesis~\ref{hyp:SDE:sigma}.
\begin{hypothesis} \label{hyp:SDE:sigma}
	Let $\beta \in [-1,1]$, $X_0 \beta \geq 0$ if $|\beta|=1$, 
and the diffusion coefficient and drift coefficient satisfy:
	$b_s=b(X_s)$, with $b$ bounded measurable function,
	$\sigma \in C^1(\RR\setminus\{r\})$ admits a finite jump at a fixed {\it threshold} $r\in \RR$, its derivative admits a finite jump at $r$ as well,
	$\sigma$ is bounded from below by a strictly positive constant and there exists a strictly increasing function $\Sigma$ such that $(\sigma(x)-\sigma(y))^2 \leq |\Sigma(x)-\Sigma(y)|$.
\end{hypothesis}

\begin{remark}
Hypothesis~\ref{hyp:SDE:sigma} is stronger than necessary. For instance, there is no need of a time-homogeneous drift, and the boundedness assumption 
can be relaxed as well.
\end{remark}

Under the above assumptions there exists a unique strong solution to the SDE~\eqref{eq:SDE:sigma} (see e.g.~\cite{BassChen05}).

\begin{theorem} \label{th:crossings}
Let $X$ be solution to~\eqref{eq:SDE:sigma} and let  $\sigma_\pm := \lim_{x\to r^\pm} \sigma(x)$. Let 
\begin{equation} 
	c_{\sigma_\pm,\beta} := \frac{2 (1-\beta^2)}{(1+\beta) \sigma_-+ (1-\beta) \sigma_+} \sqrt{\frac2\pi}  
\end{equation}
and $K_{\sigma_\pm,\beta}:=K^Y_{h_0}$ is the non-negative finite constant given by~\eqref{eq:k:th} with $Y$ $(\beta,\sigma_\pm)$-OSBM.
Then, 
$\mathcal{L}^r_{T,N}(X)$ in \eqref{def:LT:discr:1}, counting the number of times the process $X$ crosses its threshold $r$, satisfies
\begin{equation} 
	\mathcal{L}^r_{T,N}(X)  \convp[N\to\infty] c_{\sigma_\pm,\beta} L_T^r(X) 
\end{equation}
and, as $N$ goes to infinity,
\begin{equation} \label{eq:crossings_T}
	(N/T)^{1/4} \left( \mathcal{L}^r_{T,N}(X) - c_{\sigma_\pm,\beta} L^r_T(X)\right)  \ \text{ converges stably in law to } \ \sqrt{K_{\sigma_\pm,\beta}} \, \sqrt{L^r_{T}(X)} G 
\end{equation} 
where $G\sim \mathcal N(0,1)$ is independent from $L^r_T(X)$ (more generally from $\mathcal F_T^X$: the natural filtration). 
\end{theorem}

\begin{remark}
If we would not have assumed to have observations until the fixed time horizon $T$, the result could be provided as a functional limit theorem as Theorem~\ref{th:gen:jacod}:
\[\varepsilon^{(r,h_0,X)} \convucp[n\to\infty]  c_{\sigma_\pm,\beta} L^r_{\cdot}(X) \]
and then there exists (possibly on an extension of the probability space) a BM $B$ independent of $X$ such that
\[
\begin{split}
	& n^{1/4} 
	\left(  \varepsilon_{n,\cdot}^{(r,h_0,X)} -  c_{\sigma_\pm,\beta} L^r_\cdot(X)\right)
	\convsl[n\to\infty]
\sqrt{K_{\sigma_\pm,\beta}} B_{L^r_\cdot(X)}.
\end{split} 
\]
\end{remark}

\begin{remark}
We provide Theorem~\ref{th:crossings} for a more general class than OSBM, but an extension of the more general result is object of further research. This is related to the fact that the statistics of Theorem~\ref{th:crossings} basically depends only on the sign of the process.
\end{remark}

\subsubsection{Another estimator}

Let us consider the estimator $L^r_{T,N}$ in \eqref{def:LT:discr}. 
We specify its asymptotic behavior for OSBM.
For an OBM solution to \eqref{eq:OBM}, say $Y$, a proof that 
$L^r_{T,N}(Y)  \convp[N\to\infty] L_T^r(Y)$ 
can be found in \cite[Lemma~1]{lp2}.
Applying Proposition~\ref{prop:lmt:2:OSBM} to $\varepsilon^{(r, 2 h_1,Y)}$  (recall that $\beta=0$) we obtain a more general result. 
And applying Theorem~\ref{th:gen:jacod} 
we obtain the convergence rate.
We specify the asymptotic properties of $L^r_{T,N}(Y)$ in Proposition~\ref{prop:localtime:conv:OSBM} below.

The following proposition specifies, if the process is an OSBM, the constants in Proposition~\ref{prop:lmt:2:OSBM} and Theorem~\ref{th:gen:jacod} in case of the estimator \eqref{def:LT:discr} and shows its asymptotic properties.

\begin{proposition} \label{prop:localtime:conv:OSBM} 
Let $X$ be let $X$ be the $(\beta,\sigma_\pm)$-OSBM solution to~\eqref{eq:OSBM} and let ${L}^r_{T,N}(X)$ be the 
estimator of the local time $L^r_T(X)$ in \eqref{def:LT:discr}.
Then there exists (possibly on an extension of the probability space) a Gaussian random variable $G\sim \mathcal N(0,1)$ independent of $X$ (independent from $\mathcal{F}_T^X$) such that, as $N$ tends to infinity,
\begin{equation}
 	 (N/T)^{1/4} \left( L^r_{T,N}(X) -  \frac{(1-\beta^2) (\sigma_-+\sigma_+)}{(1+\beta)\sigma_- + (1-\beta )\sigma_+} L^r_T(X)\right)
\end{equation}
converges stably in law to $\sqrt {K_{\beta, \sigma_\pm}} \sqrt{L^r_T(X)} G$ where $K_{\beta, \sigma_\pm}:=K^X_{2 h_1}$ in~\eqref{eq:k:th}.

Moreover if $\beta=0$ ($X$ is an OBM) then $\frac{(1-\beta^2) (\sigma_-+\sigma_+)}{(1+\beta)\sigma_- + (1-\beta )\sigma_+}=1$ and
$K_{0,\sigma_\pm}= \frac{16}{3 \sqrt{2 \pi}} 
	 \frac{\sigma_-^2 +\sigma_+^2}{\sigma_-+\sigma_+}$.
And if $\sigma_-=\sigma_+ =1$ ($X$ is a SBM) then $\frac{(1-\beta^2) (\sigma_-+\sigma_+)}{(1+\beta)\sigma_- + (1-\beta )\sigma_+}=1-\beta^2$.
\end{proposition}

\subsection{Estimating the parameters of skew and oscillating Brownian motion} \label{sec:estimator_skew}

\subsubsection{Estimating the skewness parameter for SBM}
Theorem~\ref{th:gen:jacod} provides several estimations of the skewness parameter of SBM. They are also estimators for drifted SBM (see Remark~\ref{rem:drift}).
We consider here two new estimators and show an asymptotic mixed normality property as a straightforward application of the main result of this paper.
We choose two examples, one that involves statistics which check if there is a sign change between two consecutive time observations and one that does not. 
In both cases the computation of the constants involved in the CLT is not trivial at all, so we leave them implicit.
Comparing estimators goes beyond the purpose of this document.
And the study of the asymptotic behavior of the MLE estimator 
has been considered in~\cite{lm23} (after being conjectured in \cite{lmt1}), where a fine study of the constant involved in the CLT is provided as well. Let us mention that the latter results exploit our~Theorem~\ref{th:gen:jacod}.

Let $X$ be a $\beta$-SBM with threshold $r$, solution to~\eqref{eq:SBM}.
Recall that $L^r_{T,N}$ 
is given in~\eqref{def:LT:discr} and its asymptotic properties are provided in 
Proposition~\ref{prop:localtime:conv:OSBM}.
Let us introduce the following function $f_1(x,y):=\ind{[0,1)}(x)$. Then both $\bfF{f_1}$ in~\eqref{eq:muX} and $f_1^2$ coincide with $f_1$ and satisfy Hypothesis~\ref{hp:OSBM:1}. The same holds for the function $|f_1|(x):=f_1(|x|)$.

We consider the following estimators for the parameter $\beta$: 
\begin{equation} \label{eq:beta:estim:1}
\hat{\beta}_N 
	:= 
	\frac{ \sum_{k=0}^{N -1} \sgn(X_{k}-r) \ind{(-1,1)} (\sqrt{N/T}(X_{k}-r)) }{ \sum_{k=0}^{N -1} \ind{(-1,1)} (\sqrt{N/T}(X_{k}-r)) }
\end{equation}
or
\begin{equation} \label{eq:beta:estim:2}
\hat{\beta}_N 
	:= 
   1- \frac{L^r_{T,N}(X)}{ \sqrt{\frac{T}{N}} \sum_{k=0}^{N -1} \ind{[0,1)} (\sqrt{N/T}(X_{k}-r)) }.
\end{equation}
Note that both estimators above require the knowledge of the discretization step $T/N$. 

The following result is a direct application of Theorem~\ref{th:gen:jacod} and of the fact that joint stable convergence for two sequences holds if one of the sequence converges stably and the other converges in probability. 

\begin{theorem} \label{th:skew:estim}
Let $X$ be the $\beta$-SBM solution to~\eqref{eq:SBM} starting at $X_0=r$.
Then $\hat{\beta}_N$ in \eqref{eq:beta:estim:1} (resp.~\eqref{eq:beta:estim:2}) is a consistent estimator of $\beta$ and as $N\to\infty$
\begin{equation}
	N^{1/4} \left( \hat{\beta}_N  -\beta \right)  \ \text{ converges stably in law to } \  K \, T^{1/4} \, \frac{G}{\sqrt{L^r_T(X)}} 
\end{equation}
with $G\sim \mathcal N(0,1)$ is a standard Gaussian random variable independent of $X$ (of $\mathcal F_T^X$) and $K={\sqrt{K_{(1-\beta) f_1 -2 h_1}}}/{(1+\beta)}$ where $K_{(1-\beta)f_1-2 h_1}$ is provided by~\eqref{eq:k:th} (resp.~$K=\frac12 \sqrt{K_{(\sgn(\cdot) -\beta)|f_1|}}$).
\end{theorem}
\noindent The above result tells that $\hat{\beta}_N$ is a consistent estimator of the skewness parameter $\beta$, and $N^{1/4}(\hat{\beta}_N -\beta)$ for $N$ large behaves like a mixed Gaussian law.
Let us observe one again that, in the same way, we can provide other estimators satisfying Theorem~\ref{th:skew:estim} (with a suitable constant $K$), as the MLE studied in~\cite{lmt,lm23} respectively for $\beta=0$ and $\beta \in (-1,1)$. 
Moreover, as in~\cite{lmt}, this paves the way to hypothesis testing to check whether some observations come from (drifted) BM or (drifted) SBM.

\begin{remark}
The assumption that the $\beta$-SBM starts from $r$ is to ensure on one side that the estimator~\eqref{eq:beta:estim:1} (resp.~\eqref{eq:beta:estim:2}) is well defined and on the other side that the local time of the process at $r$ does not vanish at time $T$.
\end{remark}

\begin{remark}[Diffusion parameter estimation]
Similarly, we could provide an estimator for the parameters of OBM.
Let $Y$ be the solution to~\eqref{eq:OBM} starting at $X_0=r$, $f_+ \colon x \to \ind{[0,1)}(x)$ and $f_- \colon x \to \ind{(-1,0]}(x)$, and let $\pm\in \{-,+\}$. 
Then
\begin{equation}
\hat{\sigma}_\pm^2 
	:= 
	\frac{L^r_{T,N}(Y)}{ \sqrt{\frac{T}{N}} \sum_{k=0}^{N -1} f_\pm (\sqrt{N/T}(X_{k}-r))  }
\end{equation}
is a consistent estimator of $\sigma_\pm^2$ and 
$(N/T)^{1/4} \left( \hat{\sigma}_\pm^2   - \sigma_\pm^2 \right)$  
converges stably in law to 
a
(mixed) normal law.
%
Since the estimator proposed relies mostly on the behavior of the process around the threshold, it is slower than the ones based on quadratic variations and occupation times of the positive and negative part of the process proposed in~\cite{lp} for OBM which exploit the entire trajectory.
Therefore, we propose in Proposition~\ref{lem:estim:sigma}, in the next section, an estimator based on quadratic variation for the diffusion coefficient of OSBM.
\end{remark}

\subsubsection{Joint parameter estimation for oscillating-skew Brownian motion}

In this section, for simplicity, let the threshold $r=0$.
The estimators proposed in~\cite{lp} work also for estimating the diffusion coefficient of OSBM.
\begin{proposition} \label{lem:estim:sigma}
Let $X$ be the standard $(\beta,\sigma_\pm)$-OSBM solution to~\eqref{eq:OSBM}, and for any process $Z$ observed at times $j T/N$ let
$[Z,X]_N$ denote $\sum_{j=0}^{N-1} (Z_{j+1}-Z_j)(X_{j+1}-X_j)$.
Then, the estimator $((\widehat \sigma_+^V)^2, \ (\widehat \sigma_-^V)^2)$  of $(\sigma_+^2, \ \sigma_-^2)$ given by
\[ (\widehat \sigma_+^V) ^2 := \frac{ [ \max\{0,X\}, X ]_N}{\frac{T}{N} \sum_{j=0}^{N-1} \ind{\{X_{j} \geq 0\}} } , 
\quad  (\widehat \sigma_-^V) ^2 := \frac{ [ \max\{0, -X\}, X ]_N}{\frac{T}{N} \sum_{j=0}^{N-1} \ind{\{- X_{j} \geq 0\}} } \]
satisfies the following convergence property:
$\sqrt N \left( (\widehat \sigma_+^V)^2 -   \sigma_+ ^2, \quad (\widehat \sigma_-^V)^2 -   \sigma_- ^2\right)$ converges in law to
$
	\sqrt{2 T} \left( \frac{\sigma_+^2}{(1+\beta)^2} \frac{\int_0^T \ind{\{ X_s \geq 0 \}} \vd B_s}{\int_0^T \ind{\{ X_s \geq 0 \}} \vd s}, 
	\quad  \frac{\sigma_-^2}{(1-\beta)^2} \frac{\int_0^T \ind{\{ X_s < 0 \}} \vd B_s}{\int_0^T \ind{\{ X_s < 0 \}} \vd s} \right) 
$ 
where $B$ is a BM independent of $X$ (possibly, on an extension of the filtered probability space).
\end{proposition}
When $T=1$ and the process is OBM,  Proposition~\ref{lem:estim:sigma} is~\cite[Theorem~3.5]{lp}. 
By the scaling property~\eqref{Yscaling} one can take $T\neq 1$, and the result for OSBM follows by a spatial transformation given in the proofs section, more precisely Section~\ref{sec:interplay}.

Combining the main result of this article and Proposition~\ref{lem:estim:sigma} allows to estimate jointly $\beta, \sigma_-$ and $\sigma_+$ for an OSBM.
An example is the estimator 
\begin{equation}
	\widehat \theta_N:=( \widehat \beta, (\widehat \sigma_-^V)^2,(\widehat \sigma_+^V)^2) \quad \text{of} \quad \theta:=(\beta, \sigma_-^2,\sigma_+^2)
\end{equation} 
where
\begin{equation}
\widehat \beta := \frac{ (\widehat \sigma_+^V)^2 \varepsilon_{\frac{N}{T},T}^{(0,f_{1},X)} - (\widehat \sigma_-^V)^2  \varepsilon_{\frac{N}{T},T}^{(0,f_{-1},X)}}{(\widehat \sigma_+^V)^2  \varepsilon_{\frac{N}{T},T}^{(0,f_{1},X)}+ (\widehat \sigma_-^V)^2  \varepsilon_{\frac{N}{T},T}^{(0,f_{-1},X)}}
\end{equation}
and $f_{1}=\ind{(0,1)}$ and $f_{-1}=\ind{(-1,0)}$.
Then $N^{1/4} ( \widehat \theta_N- \theta)$ converges in law to the law of a vector $(Z,0,0)$ with $Z$ following a mixed Gaussian law.
Considering the convergence rate for $\widehat \sigma_\pm^V$ and further studying the joint behavior goes beyond the aims of this article.

\section{Proofs of the main results} \label{sec:proofs}

In this section we comment the results and their proofs. 
We first deal with the convergence in probability to the local time in Proposition~\ref{prop:lmt:2:OSBM}, which was already known for SBM. Another proof of Proposition~\ref{prop:lmt:2:OSBM} is also given in Appendix~\ref{app:prop:lmt:2}. 
Then we deal with the rate of (stable) convergence in the case of OSBM in Theorem~\ref{th:gen:jacod} whose proof is provided in Section~\ref{sec:proof:summary} relying on a well known CLT. Finally we prove Theorem~\ref{th:crossings}.

The fact that, via a spatial transformation, it is possible to change a {\em skewed\/} behavior into an {\em oscillating\/} behavior (i.e.~discontinuous diffusion coefficient) and viceversa is crucial in the proofs of the above results. 
In particular, we prove in Section~\ref{sec:interplay} the relationship between OSBM (in particular SBM) and OBM.

From now on, we take $r=0$ for simplicity. Indeed if $X$ is a $(\beta,\sigma_\pm)$-OSBM 
with threshold $r=0$ then $X-r$ is a $(\beta,\sigma_\pm)$-OSBM  
with threshold $r$. The same holds for solutions to~\eqref{eq:SDE:sigma} 
(after translation of the coefficients).

\subsection{It\^o-Tanaka formula and interplay between SBM and OBM} \label{sec:interplay}

SBM and OBM are strongly related in the following sense: 
Let $\sigma$ be the function in~\eqref{sigmaOBM}, let $Y$ the solution to \eqref{eq:OBM} and $X$ be a SBM solution to \eqref{eq:SBM} with skewness parameter ${{\beta_\sigma}:=\frac{\sigma_--\sigma_+}{\sigma_-+\sigma_+}}$ and suitable initial condition: 
Solution to the SDE
\begin{equation} \label{eq:SBM:2}
X_t
=\frac{Y_0}{\sigma(Y_0)}+ W_t+ \beta_\sigma L_t(X) .
\end{equation}
It holds that 
$ {  X_t =\frac{Y_t}{\sigma(Y_t)} }$,
and the local times satisfy  
\begin{equation}
    \label{eq:local_time:eq}
    L(X)=\frac{\sigma_++\sigma_-}{2\sigma_+\sigma_-}L(Y).
\end{equation}
Furthermore, since $ { X_t = \frac{Y_t}{\sigma(Y_t)}  }$ and $\sigma$ is positive, then $Y_t X_t \geq 0$. And, since $\sigma$ depends only on the sign of its argument,
then $\sigma(Y_t)=\sigma(X_s)$ and so it holds also ${ Y_t=\sigma(X_t) X_t }$.\\

These equalities can be proven via the {\it symmetric} (i.e.~we consider the symmetric local time) It\^o-Tanaka formula, which we specify in Lemma~\ref{lem:Ito-Tanaka} for the reader's convenience. 
In the following result $f'(y+)$ stands for right derivative of the function $f$ at $y\in \RR$ and $f'(y-)$ for left derivative and $f'(y)=\frac{f'(y+)+ f'(y-)}{2}$.
\begin{lemma} \label{lem:Ito-Tanaka}
Let $f$ be the difference of two convex functions such that $f\in C^{2}(\mathbb{R}\setminus\{0\})$ 
and $f(0)=0$,
and let $X$ be a continuous semi-martingale satisfying~\eqref{eq:SDE:sigma} (with drift $b_\cdot$, diffusion coefficient $\sigma$, skewness parameter $\beta$).
\\
Then $Y_t:=f(X_t)$ satisfies: 
\begin{multline}
	Y_t=f(X_t) = f(X_0) + \int_0^t \frac{f'( X_s+) + f'(  X_s -)}{2} b_s \vd s +  \int_0^t \frac{f'( X_s+) + f'(  X_s -)}{2} \sigma(X_s) \vd W_s
	\\+ \frac12 \int_0^t \frac{f''(X_s+) + f''(X_s-)}{2} \sigma(X_s)^2 \vd s 
	+ \frac{f'(0+) (1+\beta)  - f'(0-) (1-\beta)}{ 2 } L_t(X).
\end{multline}
If in addition $f$ is invertible 
and $f'(0+) (1+\beta) + f'(0-) (1-\beta) \neq 0$ then
\begin{multline}
	Y_t = Y_0 + \int_0^t \left( f'( f^{-1}(Y_s) ) b_s + \frac12 f''(f^{-1}(Y_s)) \sigma(f^{-1}(Y_s))^2 \right) \vd s 
	\\
	+ \int_0^t  f'(f^{-1}(Y_s)) \sigma(f^{-1}(Y_s))  \vd W_s 
	+ \frac{f'(0+) (1+\beta) - f'(0-) (1-\beta) }{ f'(0+) (1+\beta) + f'(0-) (1-\beta)} L_t(Y) 
\end{multline}
and $L_t(Y) = \frac{f'(0+) (1+\beta) + f'(0-) (1-\beta)}{2} L_t(X)$.
\end{lemma}
\begin{proof}
We first recall the well known symmetric It\^o-Tanaka formula. 
Let $ Z $ be a continuous semi-martingale and let $h$ be the difference of two convex functions.
Then, for any $t\ge 0 $,
		\begin{equation}\label{eq_Ito_Tanaka_plain}
			h(Z_t) = h(Z_0) + \int_{0}^{t} \frac{h'(Z_s+)+h'(Z_s-)}2 \vd Z_s + \frac{1}{2} 
			\int_{\mathbb R} L^y_t(Z) h''(\!\vd y)
		\end{equation}
where 
$h''(\!\vd y)$ is a measure satisfying 
\[
		\int_\RR g(y) h''(\!\vd y) = - \int_\RR g'(y) \frac{h'(y+) + h'(y-)}2 \vd y
\]
for all $g\in C^1(\mathbb{R})$ with compact support.
The occupation time formula helps dealing with the integral of the part of $h''(\!\vd y)$ which is absolutely continuous with respect to the Lebesgue measure. 
Note that $f''(\! \vd y) = f''(y)\vd y + \frac{f'(0+) - f'(0-)}{2} \delta_{0}(\vd y)$.
Hence, \eqref{eq_Ito_Tanaka_plain} applied to $f$ yields 
the first equation in the statement.
On one hand, this and \eqref{eq_Ito_Tanaka_plain} applied to $|f|$ yield
\begin{multline}
	|Y_t|=|f(X_t)| = |f(X_0)| +  \int_0^t  \sgn(f(X_s))  \vd Y_s  + \frac{f'(0+) (1+\beta) + f'(0-) (1-\beta) }{ 2 }   L_t(X).
\end{multline}
On the other hand, Tanaka formula for continuous semi-martingales states
\[
	|Y_t|= |Y_0| + \int_0^t \sgn(Y_s) \vd Y_s + L_t(Y).
\]
Hence the conclusion on the local times. 
The fact that $f$ is invertible is then used to complete the proof.
\end{proof}

\begin{example}[SBM associated to OSBM/$\sigma_\pm$-OBM]
\label{exa:OBM-SBM} 
Let $Y$ a $(\beta,\sigma_\pm)$-OSBM, in particular it satisfies~\eqref{eq:OSBM}.
Let $f(x):=\frac{x}{\sigma(x)}$.
Note that $f$ is continuous and for all $x\neq 0$ we have $f'(x)=\frac1{\sigma(x)}$ and $f''(x)=0$. 
Since $\sigma(x)$ depends only on the sign of $x$, then $f'\circ f^{-1}$ coincides with $f'$.
By Lemma~\ref{lem:Ito-Tanaka}, $X_t:=f(Y_t)$ satisfies
~\eqref{eq:SBM:2} with 
$
	\beta_\sigma= \frac{(1+\beta)\sigma_--(1-\beta)\sigma_+}{(1+\beta)\sigma_-+(1-\beta)\sigma_+}
$
and $L_{\cdot}(X) = \frac{(1+\beta)\sigma_-+(1-\beta)\sigma_+}{2\sigma_-\sigma_+}L_{\cdot}(Y)$ holds as well.
We call $X$, the {\em $\beta_\sigma$-SBM associated to the $(\beta,\sigma_\pm)$-OSBM $Y$\/}.
If $Y$ is a $\sigma_\pm$-OBM, we recognize $\beta_\sigma=\frac{\sigma_--\sigma_+}{\sigma_-+\sigma_+}$ and~\eqref{eq:local_time:eq}.
\end{example}

\begin{example}[OBM associated to a OSBM/$\beta$-SBM]
\label{exa:SBM-OBM}
Let $X$ be the solution to~\eqref{eq:SDE:sigma}. 
Take 
$\sigma_+^{(\beta)}:=\frac{1}{1 + \beta}$, 
$\sigma_-^{(\beta)}:=\frac{1}{1 - \beta}\in (0,\infty]$, 
then $\beta=\frac{\sigma_-^{(\beta)}-\sigma_+^{(\beta)}}{\sigma_-^{(\beta)}+\sigma_+^{(\beta)}}$  
and construct the diffusion coefficient 
${\sigma_{\beta}}:=\sigma_-^{(\beta)} \ind{(-\infty,0)}+\sigma_+^{(\beta)} \ind{[0,+\infty)}$ in~\eqref{sigmaOBM}. 
By Lemma~\ref{lem:Ito-Tanaka}, 
$\eta_t:=\sigma^{(\beta)}(X_t)X_t$ is solution to the SDE
\begin{equation}
	\eta_t = \sigma^{(\beta)}(X_0) (X_0) 
	+ \int_0^t \sigma^{(\beta)}(\eta_s) \sigma(\eta_s/\sigma^{(\beta)}(\eta_s)) \vd W_s 
	+ \int_0^t \sigma^{(\beta)}(\eta_s) b(\eta_s/\sigma^{(\beta)}(\eta_s)) \vd s .
\end{equation}
and
$
	L_t(\eta) = L_t(X).
$
In particular, if $X$ is the $(\beta,\sigma_\pm)$-OSBM solution to~\eqref{eq:OSBM},
then $Y_t:=\eta_t=\sigma^{(\beta)}(X_t) X_t$ is a 
$\sigma_\pm^{(\beta)}\sigma_\pm$-OBM with initial condition $Y_0={\sigma_{\beta}}(X_0) X_0$ and we refer to it as the {\em OBM associated to the $(\beta,\sigma_\pm)$-OSBM $X$\/}.
And if $X$ is the $\beta$-SBM solution to~\eqref{eq:SBM},
then $Y_t$ is the {\em $\sigma_\pm^{(\beta)}$-OBM associated to the $\beta$-SBM $X$\/}.
\end{example}

\subsection{Proof of Proposition~\ref{prop:lmt:2:OSBM}}

As already mentioned, in the case of SBM, Proposition~\ref{prop:lmt:2:OSBM} follows from \cite[Proposition~2]{lmt1} (with $T=1$) and the scaling property~\eqref{Yscaling} in Appendix~\ref{sec:scaling}.
By the transformation provided in Section~\ref{sec:interplay}, we are ready to prove 
Proposition~\ref{prop:lmt:2:OSBM} for the OSBM $Y$ solution to~\eqref{eq:OSBM}, and in particular for OBM solution to~\eqref{eq:OBM}.

However, in Appendix~\ref{app:prop:lmt:2}, we provide a proof of Proposition~\ref{prop:lmt:2:OSBM} directly for a standard OBM, with the purpose of illustrating the main ideas of the proof of Theorem~\ref{th:gen:jacod}.


Let us recall, given a measurable function $f\colon \RR^2 \to \RR$ and an OSBM $Y$ (with $r=0$), we denote by $p_Y$ its transition density~\eqref{eq:densities:X} and the function $\bfH{f}$ in \eqref{eq:muX} is computed using it. Moreover $\mu_Y$ denotes its speed measure~\eqref{eq:speedm}.

\begin{proof}[Proof of Proposition~\ref{prop:lmt:2:OSBM} for OSBM]
Let $X$ be the $\beta_\sigma$-SBM associated to the $(\beta,\sigma_\pm)$-OSBM $Y$ (see Example~\ref{exa:OBM-SBM} in Section~\ref{sec:interplay}) 
and let $f_\sigma$ be the function satisfying $f_\sigma(x,y):=f(\sigma(x)x, \sigma(y)y)$. 
For $\gamma=1,2$ it holds that
$
	\bfH {f^\gamma}(\sigma(x) x)= \bfF{(f_\sigma)^\gamma}(x)
$.
%
%
%
Hence, Remark~\ref{rem:Lgb:remarks} ensures that a function $f$ satisfies Hypothesis~\ref{hp:OSBM:1} for $Y$ (that is $\bfH{f}$, $\bfH{f^2} \in \Lgb{2}$) 
if and only if $f_\sigma$ satisfies Hypothesis~\ref{hp:OSBM:1} for $X$ (that is $\bfF{f_\sigma}$, $\bfF{f_\sigma^2} \in \Lgb{2}$). 
Moreover 
$\frac1{\sigma(x)}\mu_{\beta_\sigma}(\!\vd x) = \frac{2\sigma_-\sigma_+}{(1+\beta)\sigma_-+(1-\beta)\sigma_+} \mu_Y (\!\vd x)$.
So, 
$
		\langle \mu_{\beta_\sigma} , \bfF{(f_\sigma)^\gamma} \rangle  L_{\cdot}(X) = \langle \mu_Y , \bfH {f^\gamma} \rangle
 L_{\cdot}(Y)
$.
Applying Proposition~\ref{prop:lmt:2:OSBM} for the SBM $X$ and the function $f:=f_\sigma$ and taking into account the latter equalities complete the proof.
\end{proof}

We conclude the section with the lemma proving that the assumptions of Theorem~\ref{th:gen:jacod} are stronger than the ones of Proposition~\ref{prop:lmt:2:OSBM}, as stated in Remark~\ref{comment:convergence}.

\begin{lemma} \label{lem:link:theorems}
Let $\gamma \in [0,\infty)$ and $f\in \Lg{\gamma}$. Then $\bfH{f} \in \Lgb{\gamma}$ and $\bfH{f^2} \in \Lgb{\gamma}$. 
Let $g$ be a measurable function such that $x \mapsto |g(x)|/(1+|x|)$ is bounded,
then the function $x\mapsto \bfH{f,g}(x)/(1+|x|)$ belongs to $\Lgb{\gamma}$.
\end{lemma}
\begin{proof}
\emph{ We reduce to prove the statement for SBM. } 
Let $X$ denote the standard SBM associated to $Y$ in Example~\ref{exa:OBM-SBM}. 
Since $f\in \Lg{\gamma}$, Remark~\ref{rem:Lgb:remarks} ensures that $f_\sigma\in \Lg{\gamma}$.
Moreover 
$
	\bfH {f}(x)= \bfF{f_\sigma}(x/\sigma(x)).
$
If the first statement holds for the SBM $X$, then $\bfF{f_\sigma} \in \Lgb{\gamma}$
and Remark~\ref{rem:Lgb:remarks} ensures that $\bfH {f} \in \Lgb{\gamma}$.
We can conclude similarly for the second part of the statement, just notice that
$
	\bfH {f,g}(x)= \bfF{f_\sigma,g_\sigma}(x/\sigma(x))
$
and $|g_\sigma(x)|/(1+|x|)$ is bounded as well.

\emph{The SBM case: $Y$ is a Standard SBM.}
Let $f\in \Lg{\gamma}$. 
Let us recall that there exist a non-negative function $ \bar f \in \Lgb{\gamma}$ and constant $a \in [0,\infty)$ such that $|f(x,y)| \leq \bar f(x) e^{a|y-x|}$.
Since there exists a constant $C\in (0,\infty)$ such that for all $x,y\in \mathbb R$:
$
	p_\beta (1,x,y) \leq 
	C \frac1{\sqrt{2\pi}} e^{-\frac12(y-x)^2},
$ 
then for all $x\in \mathbb R$
\[	
	|\bfH{f}(x)| 
	\leq C \bar f(x) \int_{-\infty}^{+\infty} e^{a |y-x|}  \frac1{\sqrt{2\pi}} e^{-\frac12(y-x)^2} \vd y
	\leq 2 e^{a^2/2} C \bar f(x).
\]	
Thus $ \bfH{f} \in \Lgb{\gamma}$.
And $\bfH{f^2} \in \Lgb{\gamma}$ is proven in the same way, since $\bar f$ is bounded.
Similarly for the second part of the statement, there exists a constant $C_2\in (0,\infty)$ such that
for all $x\in \mathbb R$
\[	
\begin{split}
	|\bfH{f,g}(x)| 
	& \leq C_2 \bar f(x) \int_{-\infty}^{+\infty} (1+|y|) e^{a |y-x|}  \frac1{\sqrt{2\pi}} e^{-\frac12(y-x)^2} \vd y
	\\
	& \leq C_2 \bar f(x) \int_{-\infty}^{+\infty} (e^{|y-x|}+|x|) e^{a |y-x|}  \frac1{\sqrt{2\pi}} e^{-\frac12(y-x)^2} \vd y
	\\
	& \leq  2 e^{(a+1)^2/2} C_2 \bar f(x) +  2 e^{a^2/2} C_2 \bar f(x) |x| 
	\leq 2 e^{(a+1)^2/2} C_2 \bar f(x) (1+|x|)
\end{split}
\]	
where we recall that $ \bar f \in \Lgb{\gamma}$.
The proof is thus completed.
\end{proof}

\subsection{Proof of the Central Limit Theorem: Theorem~\ref{th:gen:jacod}}
 \label{sec:proof:summary}

Let us focus on proving Theorem~\ref{th:gen:jacod} for OBM. At the end of this section we prove Theorem~\ref{th:gen:jacod} for OSBM via its associated OBM from Section~\ref{sec:interplay}, say $Y$.

The idea of the proof is decomposing the left hand side of~\eqref{th:gen:jacod:eq}, namely 
\[n^{1/4}\left( \varepsilon_{n,t}^{(0,h,Y)} - \langle \lambda_\sigma , \bfH h \rangle L_t(Y)\right),\] 
into a discrete martingale and a vanishing term, see~\eqref{eq:decomposition}.
The decomposition is introduced in Proposition~\ref{th:j2:2}, where it is stated that the discrete martingale satisfies the assumptions of Proposition~\ref{th:j2} (a reformulation of a special case of Theorem~3.2 in \cite{j2}) which entails the limits in Theorem~\ref{th:gen:jacod}.
Proposition~\ref{th:j2} is stated just after Proposition~\ref{th:j2:2} and the comments to its proof.

\begin{proposition} \label{th:j2:2}
Let $(Y_t)_{t\in [0,1]}$ be the OBM with threshold $r=0$ strong solution to~\eqref{eq:OBM} on a stochastic basis $(\Omega, \cF, (\cF_t)_{t\in [0,1]},\PP)$ (where the filtration is the natural one associated to $Y$), let $\gamma >3$, and let $h\in \Lg{\gamma}$.
\\
Then constant $K_h^Y$ given by~\eqref{eq:k:th} (with $f=h$, $X=Y$, $\mu_X=\lambda_\sigma$, $\beta=0$) is a non-negative finite constant, 
and there exist a sequence of stochastic processes 
$(\mathcal{V}^n_t)_{t\in [0,1]}$, $n\in \NN$,
and a sequence $(\mathcal{M}^n_t)_{t\in [0,1]}$, $n\in \NN$, of $(\mathcal{F}_{{\lfloor n t \rfloor}/{n}})_{t\in [0,1]}$-martingales
such that for all $t\in [0,1]$, $n\in \NN$
\begin{equation} \label{eq:decomposition}
	n^{1/4}\left( \varepsilon_{n,t}^{(0,h,Y)} - \langle \lambda_\sigma , \bfH h \rangle L_t(Y)\right)
	=		 
	\mathcal{M}^n_t + n^{1/4} \mathcal{V}^n_t,
\end{equation}
$
	\sup_{s\in [0,1]} m^{1/4} |\mathcal{V}^m_s| \convp[m\to\infty] 0
$, 
and for all $t\in [0,1]$, $n\in \NN$ it holds that
$
	\mathcal{M}^n_t =  \sum_{k=1}^{\lfloor n t \rfloor} \chi^n_k
$
where $(\chi^n_k)_{k\in \{1, \ldots, \lfloor n t \rfloor\}}$, $n\in \NN$,
are random variables satisfying:
\begin{enumerate}[i)]
\item \label{th:j2:item1:2}
	for all $k\in \{1, \ldots, \lfloor n t \rfloor\}$
$\chi^n_k$ is square integrable $\mathcal{F}_{\frac{k}{n}}$-measurable
and
$
	\EE\!\left[ \chi_k^n | \cF_{(k-1)/n}\right] = 0$,
\item \label{th:j2:item2:2}
for all $t\in [0,1]$ it holds that
$
	\sum_{k=1}^{\lfloor n t \rfloor} \EE\!\left[ (\chi_k^n)^2 | \cF_{(k-1)/n}\right] \convp[n\to\infty] K_h^Y L_t(Y),
$
\item \label{th:j2:item3:2}
for all $t\in [0,1]$ it holds that
$
	\sum_{k=1}^{\lfloor n t \rfloor} \EE\!\left[ \chi_k^n (Y_{k/n}-Y_{{(k-1)}/n})| \cF_{(k-1)/n}\right] \convp[n\to\infty] 0
$,
\item \label{th:j2:item4:2}
for all $\varepsilon \in (0,\infty)$ it holds that
$
	 \sum_{k=1}^{n} \EE\!\left[ |\chi_k^n|^2 \ind{\{|\chi_k^n|\geq  \varepsilon \}}| \cF_{(k-1)/n}\right] \convp[n\to\infty] 0
$.
\end{enumerate}
\end{proposition}
\noindent The proof of this result is quite technical, therefore it is provided in Appendix~\ref{sec:proof:j1}. 
Basically, it consists in generalizing to the case of OBM the fundamental procedure used in~\cite{j1} for BM which exploits properties of the transition semigroup.  
To provide the ideas and step of the proof, and for the reader convenience, a proof of Proposition~\ref{prop:lmt:2:OSBM} for OBM is provided in Appendix~\ref{app:prop:lmt:2} as an introduction to the proof of Proposition~\ref{th:j2:2}.
Moreover, a proof of the finiteness of the constant $K_h^Y$ of Theorem~\ref{th:gen:jacod} is provided in Lemma~\ref{lem:finite:K}, so that the reader does not have to go through the whole proof of Proposition~\ref{th:j2:2} to check the finiteness.

Let us now assume that Proposition~\ref{th:j2:2} holds and let us prove Theorem~\ref{th:gen:jacod} first for OBM and then for OSBM. To do so, we need to state Proposition~\ref{th:j2} which provides the asymptotic behavior of the discrete martingale of Proposition~\ref{th:j2:2}.

\begin{proposition}[cf.~Theorem~3.2 in \cite{j2}] \label{th:j2}
Let $(Y_t)_{t\in [0,1]}$ be an $(\mathcal{F}_t)_{t\in [0,1]}$-local martingale on the stochastic basis $(\Omega,\mathcal{F}, (\mathcal{F}_t)_{t\in [0,1]},\PP)$.
For $n\in \NN$ and $k=1,\ldots, n$, let 
$\chi^n_k$ be square integrable $\mathcal{F}_{\frac{k}{n}}$-measurable random variables, and assume that there are $E$ and $F$
continuous processes on 
$(\Omega, \cF, (\mathcal{F}_t)_{t\in [0,1]},\PP)$ such that $E$ has bounded variation and
\begin{enumerate}[i)]
\item \label{th:j2:item1}
$
	\sup_{s\in [0,1]} \left|\sum_{k=1}^{\lfloor n s \rfloor} \EE\!\left[ \chi_k^n | \cF_{(k-1)/n}\right]- E_s\right| \convp[n\to\infty] 0
$,
\item \label{th:j2:item2}
for all $t\in [0,1]$ it holds that
$
		\sum_{k=1}^{\lfloor n t \rfloor} \EE\!\left[ \left(\chi_k^n- \EE\!\left[ \chi_k^n | \cF_{(k-1)/n}\right]\right)^{\!2} | \cF_{(k-1)/n}\right] \convp[n\to\infty] F_t,
$
\item \label{th:j2:item3}
for all $t\in [0,1]$ it holds that
$
	\sum_{k=1}^{\lfloor n t \rfloor} \EE\!\left[ \chi_k^n (Y_{k/n}-Y_{{(k-1)}/n})| \cF_{(k-1)/n}\right] \convp[n\to\infty] 0 
$,
\item \label{th:j2:item4}
for all $\varepsilon \in (0,\infty)$ it holds that
$
	 \sum_{k=1}^{n} \EE\!\left[ |\chi_k^n|^2 \ind{\{|\chi_k^n|\geq  \varepsilon \}}| \cF_{(k-1)/n}\right] \convp[n\to\infty] 0
$, and
\item \label{th:j2:item5}
for all $t\in [0,1]$ and every 
$M$ bounded $(\mathcal{F}_s)_{s\in [0,1]}$-martingale such that for all $s\in [0,1]$ the cross variation $\langle M, Y\rangle_s$ is $\mathbb P$-a.s.~equal to 0, 
it holds that
\begin{equation}
	\sum_{k=1}^{\lfloor n t \rfloor} \EE\!\left[ \chi_k^n  (M_{k/n}-M_{{(k-1)}/n}) | \cF_{(k-1)/n}\right] \convp[n\to\infty] 0.
\end{equation}
\end{enumerate}
Then there exists a BM $B$, possibly on an extension of the probability space $(\Omega,\cF,(\mathcal{F}_t)_{t\in [0,1]},\PP)$, such that $B$ and $Y$ are independent 
and
\begin{equation} \sum_{k=1}^{\lfloor n \cdot \rfloor} \chi^n_k
\convsl[n\to\infty] E_\cdot + B_{F_\cdot}.\end{equation}
\end{proposition}

\begin{proof}[Proof of Theorem~\ref{th:gen:jacod} for OBM]
Let $Y$ be the OBM strong solution to \eqref{eq:OBM}. 
Let us first assume that the filtration is the natural one associated to $Y$. 
Observe that
\[ 
	n^\frac14\left( \varepsilon^{(0,f,Y)}_{n,\cdot} - \langle \lambda_\sigma, \bfH{f} \rangle L_\cdot(Y)\right) \convsl[n\to\infty] \sqrt{K_f^Y} B_{L_\cdot(Y)}
\] 
in $\mathbb{D}_\infty$
if and only if, for all $t\in [0,\infty)$, 
$
	n^\frac14 (\varepsilon^{(0,f,Y)}_{n,\cdot} - \langle \lambda_\sigma, \bfH{f} \rangle L_\cdot(Y))|_{[0,t]} \convsl[n\to\infty] \sqrt{K_f^Y} B_{L_\cdot(Y)}|_{[0,t]}
$
in $\mathbb{D}_t$ (indeed since $B_{L_\cdot(Y)}$ has (a.s.) continuous trajectories 
this follows, e.g.~combining~\cite[Theorem~16.2]{billingsley} and~\cite[Proposition~2.2.4]{jp}).
Moreover, without loss of generality we can reduce ourselves to prove Theorem~\ref{th:gen:jacod} on the interval $[0,1]$ (on the Skorokhod space $\mathbb{D}_{[0,1]}$) for $n\in \NN$ tending to infinity. Indeed the scaling property for the OBM and its local time (see~\eqref{Yscaling}) yields the result for all non-negative times (on the Skorokhod space $\mathbb{D}_{[0,t]}$).
(The scaling property also ensures that in Theorem~\ref{th:gen:jacod} $n$ is not necessarily a natural number, but it can stay for a positive real number tending to infinity.)

Proposition~\ref{th:j2:2} 
implies that there exists a decomposition as in \eqref{eq:decomposition} and its desired stable limit as $n\in \NN$ goes to infinity coincides with the stable limit of the sequence 
$\mathcal{M}^n$ of c\`adl\`ag $(\mathcal{F}_{\frac{\lfloor n s\rfloor}n})_{s\in [0,1]}$-martingales, $n\in \NN$. 
Indeed the fact that 
$ 
	\sup_{s\in [0,1]} n^{1/4} |\mathcal{V}^n_s| \convp[n\to\infty] 0,
$ 
implies that
for every $h\colon \mathbb{D}_1 \to \RR$ continuous and bounded it holds that 
$ |h(\mathcal{M}^n + n^\frac14 \mathcal{V}^n)-h(\mathcal{M}^n)| \convp[n\to\infty] 0$
and so
for every bounded continuous function $h\colon \mathbb{D}_1 \to \RR$ and bounded measurable random variable $Y\colon \Omega \to \RR$ it holds that
\[
	\lim_{n\to\infty} \EE\!\left[|h(\mathcal{M}^n + n^\frac14 \mathcal{V}^n)-h(\mathcal{M}^n)|| Y|\right] = 0.
\]

Proposition~\ref{th:j2:2} 
also ensures that
$\mathcal{M}^n$, $n\in \NN$, satisfies all assumptions, except Item~\ref{th:j2:item5}, of Proposition~\ref{th:j2} 
(with 
$E\equiv 0$ 
and
$
	F = K_f^Y L(Y)
$
where $K_f^Y$ is the constant in~equation~\eqref{eq:k:th} with $X=Y$, $\beta=0$ and $\mu_X=\lambda_\sigma$).
Item~\ref{th:j2:item5} of Proposition~\ref{th:j2} is trivial because such bounded martingale $M$, measurable with respect to the natural filtration of the OBM and orthogonal to OBM, is nothing but a constant (follows from a martingale representation theorem, see Lemma~\ref{lem:trivial:OBM} in Appendix~\ref{sec:orthogonal}). 
Therefore, applying Proposition~\ref{th:j2} as described above, completes the proof of Theorem~\ref{th:gen:jacod} under the assumption that the filtration is the one associated to $Y$.
It remains to relax this assumption on the filtration, by the same argument in~\cite[Section 2.2]{j1}.
The proof of Theorem~\ref{th:gen:jacod} is thus completed.
\end{proof}

Since Theorem~\ref{th:gen:jacod} for OBM holds, we now prove it for OSBM and, in particular, for SBM.
The proof is based on the interplay between the processes in Section~\ref{sec:interplay}.

\begin{proof}[Proof of Theorem~\ref{th:gen:jacod} for OSBM]
Let 
$\gamma >3$, $X$ the $(\beta,\sigma_\pm)$-OSBM solution to \eqref{eq:OSBM}, and $f \in \Lg{\gamma}$. 
Let $Y$ be the OBM associated to $X$ (in the notation of Example~\ref{exa:SBM-OBM} in Section~\ref{sec:interplay}) 
and let $h \colon \RR^2 \to \RR$ be the function satisfying $h(x,y) :=f(x/{\sigma^{(\beta)}}(x), y/{\sigma^{(\beta)}}(y))$. 
Remark~\ref{rem:Lgb:remarks} ensures that $f\in \Lg{\gamma}$ if and only if $h \in \Lg{\gamma}$.
Theorem~\ref{th:gen:jacod} applied to the OBM $Y$ and the function $h$ yields 
\begin{equation}
n^{1/4} 
	\left(  \varepsilon_{n,\cdot}^{(0,h,Y)} 
-  \langle \lambda_{\sigma^{(\beta)}\sigma} , \bfH { h } \rangle L_\cdot(Y) \right)
	\convsl[n\to \infty] \sqrt{K_{h}^\beta} B_{L_\cdot(Y)}
\end{equation}
with $K_{h}^\beta$ given by~\eqref{eq:k:th} (taking $f=h$, $X=Y$, $\sigma_\pm=\sigma_\pm^{(\beta)}\sigma_\pm$, $\beta=0$, $\mu_X=\lambda_{\sigma^{(\beta)}\sigma}$).
As in the proof of Proposition~\ref{prop:lmt:2:OSBM}, the fact that $Y=\sigma^{(\beta)}(X)X$, $L_{\cdot}(X)=L_{\cdot}(Y)$, and the transition densities~\eqref{eq:densities:X},
implies that for $\gamma=1,2$ it holds
$
	\langle \lambda_{\sigma^{(\beta)}\sigma} , \bfH { h^\gamma } \rangle 
	= \langle \mu_X , \bfF { f^\gamma } \rangle 
$ 
and so
\begin{equation}
\begin{split}
&n^{1/4} 
	\left(  \varepsilon_{n,\cdot}^{(0,f,X)} 
-  \langle \mu_X , \bfF { f } \rangle L_\cdot(X) \right)
	\convsl[n\to \infty] \sqrt{K_{h}^\beta} B_{L_\cdot(X)}.
\end{split}\end{equation}
The same arguments are used to rewrite the constant $K_{h}^\beta$ obtaining $K_f^X$: more precisely note that $\bfH{h}(x \sigma^{(\beta)}(x))=\bfF { f } (x)$, $\bfH{g_0}(x \sigma^{(\beta)}(x)) = \bfF { g_\beta } (x)$, and $\mathcal{P}_{h}^Y(x \sigma^{(\beta)}(x)) = \mathcal{P}_{f}^X(x) $.
\end{proof}

\subsection{Proof of Theorem~\ref{th:crossings}}

The proof relies on It\^o-Tanaka formula and Girsanov's transform and reduces to prove the result for OBM, i.e.~applying Proposition~\ref{prop:lmt:2:OSBM} and Theorem~\ref{th:gen:jacod} for $f(x,y)=h_0(x,y)= \ind{(-\infty,0)}(x y)$.
The fact that we can extend so easily the result for more general diffusions solutions to~\eqref{eq:SDE:sigma} is due to the special structure of the estimator which counts the number of crossings of the threshold level.
Extensions of Theorem~\ref{th:gen:jacod} to more general diffusions is object of further research.

\begin{proof}[Proof of Theorem~\ref{th:crossings}]
Let us recall that Hypothesis~\ref{hyp:SDE:sigma} holds. 
In particular it holds that $\sigma \in C^1(\RR\setminus\{0\})$ is strictly positive, 
admits a finite jump at the threshold $0$, 
$\sigma_\pm:=\lim_{x\to 0^\pm} \sigma(x)\in (0,\infty)$, and its derivative admits a finite jump at $0$ as well.

{\it We first prove the statement for $\beta=0$ in several steps.}

{\it First step.} Note that if $X$ is a driftless $\sigma_\pm$-OBM, then the statement follows from 
Proposition~\ref{prop:lmt:2:OSBM} 
and
Theorem~\ref{th:gen:jacod} taking $f=h_0$. 
The result holds with constants $c_{\sigma_\pm,0}$ and $K_{\sigma_\pm,0}$.

{\it Second step.} Assume $\sigma$ and its derivative are bounded (the derivative is defined on ${\RR\setminus\{0\}}$, one may take $\sigma'(0)=(\sigma'(0^+)+\sigma'(0^-)$). 
Let us consider the process $Z$ satisfying that $\PP$-a.s.~for all $t\in [0,\infty)$
\[ 
	Z_t= Z_0 + \int_0^t \sigma(Z_s) \vd W_s + \frac12 \int_0^t \sigma(Z_s) \sigma'(Z_s) \vd s.
\]
Let 
\[
	S(x):=  \sigma_+ \int_0^x \frac{1}{\sigma(y)}\vd y \, \ind{[0,+\infty)}(x) -  \sigma_- \int_x^0 \frac{1}{\sigma(y)}\vd y  \, \ind{(-\infty,0)}(x).
\]
Note that $S$ is difference of two convex functions, $S\in C^1(\RR) \cap C^2(\RR\setminus\{0\})$, $S(0)=0$, and $S$ is strictly increasing.
By Lemma~\ref{lem:Ito-Tanaka}, $Y_t:=S(Z_t)$ is an OBM 
starting at $S(Z_0)$.
Since $x S^{-1}(x) \geq 0$ for all $x\in \RR$, then for all $t$ it holds that $Z_t$ and $S^{-1}(Z_t)=Y_t$ have the same sign and, since $\mathcal{L}^0_{T,N}$ depends only on the observations through their sign, 
$\mathcal{L}^0_{T,N}(Y) =\mathcal{L}^0_{T,N}(Z)$.
Therefore we proved that Theorem~\ref{th:crossings} holds for the process $Z$ with the same constants $c_{\sigma_\pm,0}$ and $K_{\sigma_\pm,0}$. 

{\it Third step.} Assume $\sigma$ and $\sigma'$ are bounded, as in the previous step. Combining stable converge and Girsanov's transform (see Remark~\ref{rem:drift}) ensures that the same result holds for $X$ solution to~\eqref{eq:SDE:sigma} if $\beta=0$ with the same constants $c_{\sigma_\pm,0}$ and $K_{\sigma_\pm,0}$. 
Let us be more precise. 
Let $Z$ as in the previous step and let $X$ solve~\eqref{eq:SDE:sigma} with $\beta=0$.
The Dol\'eans-Dade exponential, expressing the Radon-Nikodym derivative between the law of $X$ and the one of $Z$, is $\mathcal{E}(\xi)$ with 
$
	\xi_t := \int_0^t \frac{b(Z_s)-\frac12 \sigma(Z_s) \sigma'(Z_s)}{\sigma(Z_s)} \vd W_s.
$ 
It is an exponential martingale (see e.g.~Novikov's condition in~\cite[Corollary 3.5.14]{ks}) because $b$, $\sigma'$ are bounded and $\sigma$ is bounded from below by a strictly-positive constant.
By a localization argument, the boundedness conditions for $\sigma$ and $\sigma'$ can be removed.

{\it Now let us consider $\beta\neq 0$.}

Consider the process $\eta_s:=\sigma^{(\beta)}(X_s) X_s$ of Example~\ref{exa:SBM-OBM} in Section~\ref{sec:interplay}.
Note that 
$\mathcal{L}^0_{T,N}(\eta) =\mathcal{L}^0_{T,N}(X)$ because the statistics depends only on the sign of the observations of $\eta$ and $X$ and this sign is preserved by the transformation we apply.
We have proved above that the process $\eta$ satisfies Theorem~\ref{th:crossings} with constants $c_{\sigma^{(\beta)}_\pm \sigma_\pm,0}$ and $K_{\sigma^{(\beta)}_\pm \sigma_\pm,0}$
hence
$X$ satisfies the theorem with the same constants that we just rename in order to explicit the dependence from $\sigma^{(\beta)}$ as a dependence from $\beta$.
The proof is thus completed.
\end{proof}

\appendix

\section{Properties of OBM} \label{sec:OBM:basic}

We recall known properties of the law of OBM and its local time.
And we establish bounds for the transition semigroup which are useful in the proof of the key Proposition~\ref{th:j2:2} in Appendix~\ref{sec:proof:j1}.
Note that, by the interplay between OSBM and OBM (see Section~\ref{sec:interplay}), one can derive results for OSBM analogous to the ones obtained in this section.

Unless explicitly stated, in this section we consider $Y$ to be a driftless OBM with threshold $r=0$.

\subsection{Scaling property} \label{sec:scaling}
Assume $Y$ starts 
from a deterministic point $Y_0$, let $c\in (0,\infty)$,
and let $Z$ denote the OBM with threshold $r=0$ starting from $\sqrt{c} Y_0$.
Let us mention the following well known diffusive scaling properties for OBM:
$
\bigg(\frac{1}{\sqrt{c}} Z_{ct} \bigg)_{t\geq 0} \eqlaw \big( Y_{t}\big)_{t\geq 0} 
$
(i.e.~``the rescaled OBM is still a OBM with rescaled starting point")
%
%
and
\begin{equation}\label{Yscaling}
\bigg(\frac{1}{\sqrt{c}} Z_{ct} , \frac{1}{\sqrt{c}} L_{ct}(Z)\bigg)_{t\geq 0} \eqlaw \big( Y_{t}, L_t(Y) \big)_{t\geq 0}.
\end{equation} 
\begin{remark}
The transition density satisfies for every $c>0$ that $q_\sigma(t,x,y)= \sqrt{c} q_\sigma(c t, \sqrt{c} x, \sqrt{c} y ).$
This and the Markov property imply that for all $k =0,1,\ldots$
\[
\begin{split} 
	& \frac{1}{\sqrt{n}}\bfH{f,g}(\sqrt{n}Y_{k/n}) 
	= \frac1{\sqrt n}  \EE\!\left[ f(\sqrt{n} Y_{k/n}, \sqrt{n} Y_{(k+1)/n}) g(\sqrt{n} Y_{(k+1)/n}) | Y_{k/n} \right].
\end{split}
\]
And $\rho^\sigma_t(y,\ell)$, the joint density of a standard OBM and its local time introduced in the next section, satisfies for every $c>0$: $\rho^\sigma_t(y,\ell) = c \rho^\sigma_{c t}(\sqrt c y, \sqrt c \ell)$.
\end{remark}

\subsection{The joint density of a standard OBM and its local time}
The joint density of a standard OBM (i.e.~$Y_0=0$) and its local time at time $t$,
$\rho^\sigma_t(y,\ell)$ coincides with
\begin{equation} \label{eq:OBMl}
		\rho^\sigma_t(y,\ell) = \frac{1}{(\sigma(y))^2} \rho_t\left(\frac{y}{\sigma(y)},  \frac{\sigma_-+\sigma_+}{2\sigma_-\sigma_+} \ell\right)
\end{equation}
for $y\neq 0$, where $\rho$ is the joint density of the BM and its local time at time $t$:
\begin{equation} \label{eq:BMl}
	\rho_t(y,\ell) = \frac{|y| + \ell}{\sqrt{2\pi t^3}} \exp\!\left(- \frac{(|y| + \ell)^2}{2 t}\right) \ind{(0,\infty)}(\ell).
\end{equation}
In particular
$
		\rho^\sigma_t(y,\ell) \vd y \vd \ell  =  \rho_t\left(\frac{y}{\sigma(y)},  \frac{\sigma_-+\sigma_+}{2\sigma_-\sigma_+} \ell\right) \lambda_\sigma(\!\vd y) \vd \ell.
$

\subsection{Bounds for the semigroup}
Let ${\beta_\sigma}:=\frac{\sigma_--\sigma_+}{\sigma_-+\sigma_+}$ and recall that, for every function $f\colon \RR \to \RR$, $f_\sigma(x):=f(\sigma(x) x)$. 
For a measurable, bounded function $f\colon \RR \to \RR$ set 
\begin{equation}\label{eq:def:semig}
	Q^{\sigma}_t f(x) := \int_{-\infty}^{\infty} q_\sigma(t,x, y)f(y) \vd y
	\quad { and } \quad 
	P^{\beta_\sigma}_t f(x) := \int_{-\infty}^{\infty} p_{\beta_\sigma}(t,x, y)f(y) \vd y
\end{equation}
for all $t\in [0,\infty)$.
They are respectively the semigroup of the standard OBM and of the standard $\beta_\sigma$-SBM and they satisfy
$
	Q^{\sigma}_t f(x) = P^{{\beta_\sigma}}_t f_\sigma (x/\sigma(x)).
$
Note that $P_t:=P^0_t=Q_t^{1}$ is the semigroup of the BM and
\begin{equation} \label{eq:semig}
	Q^\sigma_t f(x) = P_t f_\sigma(x/\sigma(x))+ {\beta_\sigma} P_t (f_\sigma \ind{[0,\infty)})(-|x|/\sigma(x)) - {\beta_\sigma} P_t (f_\sigma \ind{(-\infty,0)}) (|x|/\sigma(x)).
\end{equation}
From this relationship between the semigroups of OBM and BM we derive the following properties.
\begin{lemma} \label{lem:semig:estim}
Let $f\in \Lgb{2}$, and let us denote by $p(t,\cdot)$ is the density of a Gaussian random variable with variance $t$.
Then there exists a positive constant $K\in (0,\infty)$ such that for all $x,y\in \RR$, $0 \leq s\leq t$ it holds that
\begin{enumerate}[i)]
	\item  \label{item1:semig:estim}
$
	|Q^\sigma_t f (x)| \leq \frac1{\min\{\sigma_-, \sigma_+\}} \left(1+\frac{|\sigma_- -\sigma_+|}{\sigma_-+\sigma_+}\right) \frac1{\sqrt{2 \pi}} \frac{\|f\|_{1}}{\sqrt t}, 
$
	\item \label{item2:semig:estim}
$ 
	\left| Q^\sigma_t f (x)  - \frac{2 \sigma_-\sigma_+}{\sigma_-+\sigma_+} \langle \lambda_\sigma , f \rangle p(t,{x}/{\sigma(x)}) \right|\leq \frac{K}{\sqrt{t^3}} \left( \| f \|_{1,2} +\| f \|_{1,1} |x|\right),
$
	\item \label{item2bis:semig:estim} for all $\zeta \geq 0 $ there exists a positive constant $K_\zeta$ such that\\
$ 
	\left| Q^\sigma_t f (x)  - \frac{2 \sigma_-\sigma_+}{\sigma_-+\sigma_+} \langle \lambda_\sigma , f \rangle p(t,{x}/{\sigma(x)}) \right|
\leq  
\frac{K_\zeta}{t} \left( \frac{\| f \|_{1,1}}{1+(|x|/(\sigma(x)\sqrt{t}))^\zeta} +\frac{\| f \|_{1,1+\zeta}}{1+(|x|/\sigma(x))^\zeta} \right),
$
	\item  \label{item3:semig:estim}
$ 
	\left| Q^\sigma_t f (x)  - Q^\sigma_t f (y) \right| 
\leq \frac1{\min\{\sigma_-, \sigma_+\}} \left(1+\frac{|\sigma_- -\sigma_+|}{\sigma_-+\sigma_+}\right) K \frac{|x-y|}{t} \|f\|_{1}
$, and
	\item  \label{item4:semig:estim}
$ 
	\left| Q^\sigma_t f (x)  - Q^\sigma_s f (x) \right| \leq  \frac1{\min\{\sigma_-, \sigma_+\}} \left(1+\frac{|\sigma_- -\sigma_+|}{\sigma_-+\sigma_+}\right) K \frac{t-s}{\sqrt{s^{3}}} \|f\|_{1}.
$
\end{enumerate}
\end{lemma}
\begin{proof}
Item~\ref{item1:semig:estim} is a straightforward consequence of \eqref{eq:densities:X} (with $\beta=0$) and of the fact that 
\begin{equation}
	p_{\beta_\sigma} (t,x,y) 
 \leq (1+|{\beta_\sigma}|) p(t,x-y) 
\leq \frac{1+|{\beta_\sigma}|}{\sqrt{2\pi t}}.
\end{equation}
To prove the other items we also use the fact for all $\alpha \geq 0$
$\| f(\sigma(\cdot)\cdot) \|_{1,\alpha} \leq \frac{\| f \|_{1,\alpha}}{(\min{\{\sigma_-,\sigma_+ \}})^{1+\alpha}}$ and 
$\| f \ind{[0,\infty)} \|_{1,\alpha} + \| f \ind{(-\infty,0)} \|_{1,\alpha}= \| f \|_{1,\alpha}$.
Item~\ref{item2:semig:estim} follows from \eqref{eq:semig} and \cite[Lemma~1]{lmt1} for SBM.
Item~\ref{item2bis:semig:estim} follows from \eqref{eq:semig} and the analogous result for BM: equation (3.2) in~\cite[Lemma~3.1]{j1}.
Item~\ref{item3:semig:estim} and Item~\ref{item4:semig:estim} follow from \eqref{eq:semig} and equations (3.4)-(3.5) in \cite[Lemma~3.1]{j1}.
\end{proof}

We now consider the aggregated action of the semigroup which is useful in the convergence results provided in Section~\ref{sec:aux:conv}.
Let $f\colon \RR \to \RR$ and
\[
	\Gamma_t(n,f, Y_0) 
	:= 
	\sum_{k=1}^{\lfloor n t \rfloor -1} Q^\sigma_k f (\sqrt{n} Y_0).
\] 
Note that if $Y$ is an OBM, then
$\Gamma_t(n,f,x) 
	= \sum_{k=1}^{\lfloor n t \rfloor -1} \EE\left[ f(\sqrt{n} Y_k) | Y_0=\sqrt{n} x\right].
$
And by the scaling property~\eqref{Yscaling} it holds that
$
	\Gamma_t(n,f,x) = \sum_{k=1}^{\lfloor n t \rfloor -1} \EE\left[ f(Y_{k/n}) | Y_0=x\right].
$

\begin{lemma} \label{lem:j1:3.3}
Let $f\colon \RR \to \RR$ and
$
	\Gamma_t(n,f, Y_0) 
$
as above.
There exists a positive constant $K$ (depending on $\sigma_\pm$) such that
$|\Gamma_t(n,f,Y_0)| \leq K \|f\|_{1} {\sqrt{n t}}$
and if $\langle \lambda_\sigma , f \rangle=0$ then it holds that
$ |\Gamma_t(n,f,Y_0)| \leq K \left( \| f \|_{1,2} + \| f \|_{1,1} |Y_0| {\sqrt{n}}\right)$ 
and
$
|\Gamma_t(n,f.Y_0)| \leq K \| f \|_{1,1} (1+ \log(n t)).
$
\end{lemma}
\begin{proof}
This proof is analogous to the one of \cite[Lemma~3.3]{j1}.
Item~\ref{item1:semig:estim} in Lemma~\ref{lem:semig:estim} and the fact that
$ \sum_{k=1}^{\lfloor n t \rfloor -1}\frac{1}{\sqrt{k}} \leq 2  \sqrt{n t}$ prove the first result.
The second part of the statement, under the assumption that $\langle \lambda_\sigma , f \rangle=0$, follows from Items~\ref{item2:semig:estim}-\ref{item2bis:semig:estim} in Lemma~\ref{lem:semig:estim} and the facts that $ \sum_{k=1}^{\infty}\frac{1}{\sqrt{k^3}} <\infty$ and
$ \sum_{k=1}^{\lfloor n t \rfloor -1} \frac1k \leq (1+\max\{0,\log(nt)\})$. 
\end{proof}

\subsection{Local time: on time continuity and moments}
In this section we explore some properties of the local time of OBM and its moments.

The following statement is well known for the local time of BM. The proof is standard (based on It\^o-Tanaka formula and Burkholder-Davis-Gundy inequality). 

\begin{lemma} \label{lem:holderL}
For all $q\in (2,\infty)$, $\alpha \in (0,\frac{q-2}{2 q})$ 
it holds that (the pathwise continuous version of) the local time $L_\cdot(Y)$ is locally $\alpha$-H\"older continuous.
In particular for all $\delta \in (-\infty,\frac12)$, $T\in [0,\infty)$ it holds that
\begin{equation} 
	\sup_{t\in [0,T]} n^{\delta} \left(L_{ t +\frac1n }(Y) - L_t(Y)\right) \convas[n\to\infty] 0.\end{equation}
\end{lemma}
\begin{proof}
In this proof let $T\in [0,\infty)$ and $s,t \in [0,T]$ with $s\leq t$ be fixed.
Let us first note that It\^o-Tanaka formula 
implies that 
\[L_t -L_s = |Y_t|-|Y_s| -\int_s^t \sgn(Y_u) \sigma(Y_u) \vd W_u.\]
The fact that for all $a,b\in \RR$, $p\in [1,\infty)$ it holds that $(a+b)^p \leq 2^{p-1}(|a|^p + |b|^p)$ and that $||a|-|b|| \leq |a-b|$,
and Burkholder-Davis-Gundy inequality ensure that for all $q\in [2,\infty)$ there exists a constant $K_q>0$ (dependent on $q,\sigma_-$ and $\sigma_+$) such that
\begin{equation}
\begin{split}
	\EE\!\left[|L_t -L_s|^q\right] 
	& \leq K_q \EE\!\left[ \left| \int_s^t \sigma^2(Y_u) \vd u \right|^{\frac{q}2} \right]
	\leq 
	\max{\{\sigma_-,\sigma_+\}}^q 
	K_q (t-s)^{q/2}.
\end{split}
\end{equation}
Finally Kolmogorov continuity theorem ensures that there exists a continuous version of the local time (that we took already) and it is locally H\"older continuous as required. 
\end{proof}

Let us recall that $Y$ is an OBM with threshold $r=0$.
For every $p\in [0,\infty)$, $x\in \RR$, assume that $Y_0=x$ and for every function $f\colon\RR\to\RR$ either non-negative or such that $\left( L_1(Y)\right)^p f(Y_1)\in \Lone(\PP)$,
let 
\begin{equation} \label{function:G}
	\Lm{f}{p}{x} := \EE\!\left[ \left( L_1(Y)\right)^p f(Y_1) | Y_0=x \right] \!.
\end{equation}
In this document we only consider functions $f\colon \RR \to \RR$ satisfying that there exist $K,\alpha\in [0,\infty)$ such that $|f(y)|\leq K e^{\alpha |y|}$ for all $y\in \RR$, so $\Lm{f}{\cdot}{\cdot}$ is well defined. 
The scaling property~\eqref{Yscaling} in Appendix~\ref{sec:scaling} implies that 
\begin{equation} \label{function:G:equal}
	\Lm{f}{p}{\sqrt{n} Y_{\frac{(k-1)}{n}}}= n^{\frac{p}2} \EE\!\left[ \left( L_{\frac{k}n}(Y) -  L_{\frac{k-1}n}(Y) \right)^p f(\sqrt{n} Y_{\frac{k}{n}}) | \cF_{\frac{(k-1)}{n}}\right]
\!.
\end{equation}
In particular note that $\Lm{1}{1}{\cdot} = \bfH{g}(\cdot)$ with $g(x,y)=|y|-|x|$.

\begin{lemma} \label{lem:prop:G}
Let $x\in \RR$,
let $W$ be a standard BM and $Y$ be a standard OBM defined on the same probability space,
let $f\colon \RR \to \RR$ be a function satisfying that there exist $K,\alpha\in [0,\infty)$ such that $|f(y)|\leq K e^{\alpha |y|}$ for all $y\in \RR$. 
Then for all $p\in \NN$ it hold that
\begin{equation}
\begin{split} 
	& \Lm{f}{p}{x}
	 = \int_0^1 \frac{|x|}{\sigma(x)} \frac{(1-t)^\frac{p}{2}}{\sqrt{2\pi} t^\frac32} e^{-\frac{x^2}{2 (\sigma(x))^2 t}} 
	\EE\!\left[ \left( L_1(Y)\right)^p f(Y_1 \sqrt{1-t})\right]\! \vd t
	\\
	& 
	= \left(\frac{2 \sigma_-\sigma_+}{\sigma_-+\sigma_+}\right)^{\! (p+1)} \frac{|x|}{\sigma(x)}
	\int_0^1  \frac{(1-t)^\frac{p}{2}}{\sqrt{2\pi} t^\frac32} e^{-\frac{x^2}{2 (\sigma(x))^2 t}} 
	\EE\!\left[ \left( L_1(W)\right)^p (\sigma(W_1))^{-1} f(\sigma( W_1) W_1 \sqrt{1-t})  \right] \! \vd t,
\end{split}
\end{equation}
if $x\neq 0$, and 
$ 
	\Lm{f}{p}{0}
	= \left(\frac{2 \sigma_-\sigma_+}{\sigma_-+\sigma_+}\right)^{\! (p+1)} 
	\EE\!\left[ \left( L_1(W)\right)^p (\sigma(W_1))^{-1} f(\sigma( W_1) W_1) \right] \!.
$\\
(If $\sigma_\pm=\infty$ then replace $f$ in the right hand side of last two equalities with $f \ind{\RR_\mp}$.)
\end{lemma}
\begin{proof}
We reduce to consider the case $x\neq 0$ because if $x=0$ then the statement follows from simple computations using the joint density of the OBM $Y$ and its local time~\eqref{eq:OBMl}.

Let $\beta_\sigma= \frac{\sigma_- - \sigma_+}{\sigma_- + \sigma_+}$ and let $X$ be the ${\beta_\sigma}$-SBM starting from $x/\sigma(x)$. 
Let $Z_t:=Y/\sigma(Y)$ a standard ${\beta_\sigma}$-SBM and let $B$ be a BM starting at $x/\sigma(x)$ independent of $Y$.
For a process $\xi$ let us denote by
$
	T_0(\xi):= \inf(\{\infty\}\cup \{t \geq 0 \colon \xi_t = 0\})
$
the first time it hits 0.

One well known property of SBM is that the process behaves as a BM until it reaches the barrier, which is 0. This means that $T_0(X) \eqlaw T_0(B)$.
Moreover, by the Markov property, it holds that
$X_{\cdot+T_0\left(X\right)}$ conditioned on $T_0\left(X\right)$ is distributed as $Z_\cdot$.

This and the relationship between the local times of OBM and SBM~\eqref{eq:local_time:eq} show that
\begin{equation}\begin{split}
	& \Lm{f}{p}{x} 
	= \left(\tfrac{2 \sigma_-\sigma_+}{\sigma_-+\sigma_+}\right)^p 
		\EE\!\left[ \ind{\{T_0\left(X\right)\leq 1\}} \EE\!\left[  \left( L_1(X)\right)^p f(\sigma(X_1)X_1) | T_0(X)\right] \right]
	\\
	& = \left(\tfrac{2 \sigma_-\sigma_+}{\sigma_-+\sigma_+}\right)^p 
		\EE\!\left[ \ind{\{T_0\left(B\right) \leq 1\}} \left( L_{1-T_0\left(B\right)}(Z)\right)^p f(\sigma(Z_{1-T_0\left(B\right)}) Z_{1-T_0\left(B\right)}) \right]\!.
\end{split}\end{equation}
Let us recall the well known fact that the random variable $T_0(B)$ has density w.r.t. the Lebesgue measure given by
$
	(0,\infty) \ni t\mapsto \frac{|x|}{\sigma(x)} \frac{1}{\sqrt{2\pi} t^\frac32} e^{-\frac{x^2}{2 (\sigma(x))^2 t}}.
$
Then the relationship between the local times of OBM and associated SBM~\eqref{eq:local_time:eq},
the scaling property~\eqref{Yscaling} and simple changes of variables 
imply that
\begin{equation}\begin{split}
	\Lm{f}{p}{x} 
& =	\int_0^1 \frac{|x|}{\sigma(x)} \frac{(1-t)^\frac{p}{2}}{\sqrt{2\pi} t^\frac32} e^{-\frac{x^2}{2 (\sigma(x))^2 t}} \EE\!\left[ \left( L_1(Y)\right)^p f(Y_1 \sqrt{1-t})\right] \! \vd t.
\end{split}\end{equation}
The relationship between the joint density of the standard OBM and its local time~\eqref{eq:OBMl} and the one for BM and its local time~\eqref{eq:BMl} 
yield the conclusion.
\end{proof}

\subsection{Orthogonal martingales} \label{sec:orthogonal}

In this section let Y be a OBM and $(\cF_t)_{t\in [0,1]}$ its natural filtration.
In this section we show another fact that OBM has in common with the one-dimensional BM: the only orthogonal square-integrable $(\cF_t)_{t\in [0,1]}$-martingales are the constants.

\begin{lemma}\label{lem:trivial:OBM}
Let $(M_t)_{t\in [0,1]}$ be a square-integrable $(\cF_t)_{t\in [0,1]}$-martingale such that for all $t\in [0,1]$ the cross variation satisfies $\PP(\langle M_t, Y_t\rangle_t = 0)=1$. Then $M$ is constant.
\end{lemma}
\begin{proof}
Without loss of generality we can assume $M_0=0$, otherwise consider $M_t-M_0$.
There exists an $(\mathcal{F}_t)_{t\in [0,1]}$-progressively measurable process $\nu$ such that $M_t = \int_0^t \nu_s \vd Y_s$ for all $t\in [0,1]$.
The orthogonality of $M$ to $Y$ rewrites as follows for all $t\in [0,1]$ it holds $\PP(\int_0^t \nu_s \sigma^2(Y_s) \vd s = 0)=1$. 
By continuity, we deduce that $\PP$-a.s.~for all $t\in[0,1]$ it holds that $\int_0^t \nu_s \sigma^2(Y_s) \vd s = 0$. 
Thus $\PP$-a.s. $\int_0^1 \ind{\{\nu_s \neq 0\}} \vd s =0$. Hence for all $t\in [0,1]$ it holds $\PP$-a.s.~that $M_t = \int_0^t \nu_s  \vd Y_s = 0$.
 Theorem~4.15 in \cite{ks} ensures the a.s.~uniqueness of the paths of the process $\nu$ in $L^2([0,1];\RR)$ and so $M_t \equiv 0$.
\end{proof}

\section{Towards proving Proposition~\ref{th:j2:2}: Sketch of a proof of Proposition~\ref{prop:lmt:2:OSBM}} \label{app:prop:lmt:2}

In this section we prove  
Proposition~\ref{prop:lmt:2:OSBM} for OBM with threshold $r=0$. 
More precisely we focus on proving that
\begin{equation} \label{eq:conv:decomp:fake}
	\sup_{s\in [0,t]} \left| \varepsilon_{n,s}^{(0,h,Y)} - \langle \lambda_\sigma , \bfH h \rangle L_s(Y)\right|
	\convp[n\to\infty] 0 
\end{equation}
under Hypothesis~\ref{hp:OSBM:1} (i.e.~$\bfH{h}$, $\bfH{h^2} \in \Lgb{2}$).
Since the proof of Proposition~\ref{th:j2:2} is quite technical, we exploit the proof of the latter convergence to provide a clearer picture to the reader and, at the same time, to introduce some results which are useful in the proof of Proposition~\ref{th:j2:2}.

First, we provide the following decomposition
\begin{equation} \label{eq:decomposition:bis}
	 \varepsilon_{n,t}^{(0,h,Y)} - \langle \lambda_\sigma , \bfH h \rangle L_t(Y)
	=		 
	\widehat{\mathcal{M}}^n_t + \mathcal{V}^n_t
\end{equation}
where $\widehat{\mathcal{M}}^n_\cdot$ is a $(\cF_{\lfloor n t \rfloor/n})_{t\geq 0}$-martingale.
The decomposition in Proposition~\ref{th:j2:2} is such that $n^{1/4}\widehat{\mathcal{M}}^n_t$ has a non-trivial limit involving the local time.
In this section we are not interested on the limit, therefore the decomposition we take here, is slightly different that the one we consider in Section~\ref{sec:decomposition}. 
Indeed here both terms of the decomposition should vanish as $n$ grows.
How to find $\widehat{\mathcal{M}}^n_t$ and $\mathcal{V}^n_t$? 
A first step is to look for an approximation of the local time.

\subsection{A statistics approximating the local time}

The following lemma is crucial in because it identifies a statistics which provide an approximation of the local time.

\begin{lemma}\label{lem:conv:hatg}
Let $g\colon \RR^2 \to \RR$ be the real function satisfying
$
	g(x,y)=|y|-|x|.
$
Then for all $t\in [0,\infty)$ it holds that $\langle \lambda_\sigma , \bfH{g} \rangle=1$,
\begin{equation}
	\varepsilon_{n,t}^{(0,\bfH{g},Y)} := \frac1{\sqrt{n}}\sum_{k=0}^{\lfloor n t \rfloor-1} \bfH{g} (\sqrt{n}Y_{k/n})
	=  \sum_{k=0}^{\lfloor n t \rfloor-1}  \EE\!\left[  L_{(k+1)/n}(Y) - L_{k/n}(Y) | \cF_{k/n} \right]
\end{equation}
and 
$
	\sup_{s\in [0,t]} \Big| 
\varepsilon_{n,s}^{(0,\bfH{g},Y)} 
- L_s(Y) \Big| \convp[n\to\infty] 0.
$
\end{lemma}
\begin{proof}
Some computations show that $\langle \lambda_\sigma , \bfH{g} \rangle=1$. 
 In fact it holds that
\begin{equation} \begin{split}
	& \langle \lambda_\sigma , \bfH{g} \rangle 
	 = \int_{-\infty}^{\infty} \int_{-\infty}^{\infty} \frac{(|y|-|x|)}{(\sigma(x))^2} q_\sigma(1,x,y) \vd y \vd x
	\\
	& =  \frac{1-{\beta_\sigma}}{\sigma_- \sigma_+^2} \int_{-\infty}^{0} \int_{0}^{\infty} \frac{- (y+x)}{\sqrt{2 \pi}} e^{- \frac12 \left(\frac{x}{\sigma_+}-\frac{y}{\sigma_-}\right)^2} \vd x \vd y
	 +  \frac{1+{\beta_\sigma}}{\sigma_-^2 \sigma_+} \int_{-\infty}^{0} \int_{0}^{\infty}  \frac{(y+x)}{\sqrt{2 \pi}} e^{- \frac12 \left(\frac{x}{\sigma_-}-\frac{y}{\sigma_+}\right)^2} \vd y \vd x
	\\
	& \quad +  \int_{0}^{\infty} \int_{0}^{\infty} \frac{(y-x)}{\sigma_+^3 \sqrt{2 \pi}} 
	\left( e^{- \frac12 \left(\frac{x}{\sigma_+}-\frac{y}{\sigma_+}\right)^2} 
		+ {\beta_\sigma} e^{- \frac12 \left(\frac{x}{\sigma_+}+\frac{y}{\sigma_+}\right)^2} 
	\right)\vd x \vd y 
	\\
	& \quad + \int_{-\infty}^0\int_{-\infty}^0 \frac{-(y-x)}{\sigma_-^3 \sqrt{2 \pi}} 
	\left( e^{- \frac12 \left(\frac{x}{\sigma_-}-\frac{y}{\sigma_-}\right)^2} 
		- {\beta_\sigma} e^{- \frac12 \left(\frac{x}{\sigma_-}+\frac{y}{\sigma_-}\right)^2} 
	\right)\vd x \vd y.
\end{split}\end{equation}
The first two terms of the right-hand-side cancel because $\frac{1-{\beta_\sigma}}{\sigma_- \sigma_+^2} = \frac{1+{\beta_\sigma}}{\sigma_-^2 \sigma_+}$, so simple change of variables show that
\begin{equation} \begin{split}
	& \langle \lambda_\sigma , \bfH{g} \rangle 
	 = 
	\int_{0}^{\infty} \int_{0}^{\infty} \frac{(y-x)}{\sqrt{2 \pi}} 
	\left( e^{- \frac{(x-y)^2}2} 
		+ {\beta_\sigma} e^{- \frac{(x+y)^2}2} 
		+ e^{- \frac{(x-y)^2}2} 
		- {\beta_\sigma} e^{-\frac{(x+y)^2}2}
	\right)\vd x \vd y
	\\
	& = 2 \int_0^\infty   \int_{0}^{\infty} -\frac{(x-y)}{\sqrt{2 \pi}} 
	 e^{- \frac{(x-y)^2}2} \vd x \vd y
	 = 2 \int_0^\infty  \frac1{\sqrt{2 \pi}} 
	 e^{- \frac{(0-y)^2}2}  \vd y =1.
\end{split}\end{equation}

The scaling property~\eqref{Yscaling}, the Markov property, and It\^o-Tanaka formula (see Lemma~\ref{lem:Ito-Tanaka}) show 
$\frac1{\sqrt{n}}{\bfH{g}(\sqrt{n}Y_{k/n})} = \EE\!\left[  L_{(k+1)/n}(Y) - L_{k/n}(Y) | \cF_{k/n} \right]\!.$
Indeed, let us observe that a simple change of variable (corresponding to the scaling property~\eqref{Yscaling}) and the Markov property yield
\begin{equation}
\begin{split} 
	& \frac{1}{\sqrt{n}}\bfH{g}(\sqrt{n}Y_{\frac{k}{n}}) 
	= \frac1{\sqrt{n}} \int_{-\infty}^{\infty} \left( |y|- \sqrt{n}|Y_{\frac{k}{n}}| \right) q_\sigma(1,\sqrt{n}\, Y_{\frac{k}{n}}, y) \vd y
	\\
	& = \int_{-\infty}^{\infty} \left( |y|- |Y_{\frac{k}{n}}| \right) q_\sigma(\tfrac1{n},Y_{\frac{k}{n}}, y) \vd y
	= \EE\!\left[ |Y_{\frac{(k+1)}{n}}| - |Y_{\frac{k}{n}}| | |Y_{\frac{k}{n}}| \right] = \EE\!\left[ |Y_{\frac{k+1}{n}}| - |Y_{\frac{k}{n}}| | \cF_{\frac{k}{n}} \right] \!.
\end{split}
\end{equation}

Lemma~2.14 in \cite{j3} and Lemma~\ref{lem:holderL} ensure the desired convergence in probability.
\end{proof}

Observe that $\bfH{g}(x)$ is the difference between the expected value of $|Y_{1/n}|$ and $|Y_0|=|x|$
and $\varepsilon_{n,t}^{(0,\bfH{g},Y)}$ is a sort of average of all absolute value increments between two consecutive observations.

\subsection{Identifying the martingale term} \label{app:decomp:fake}

Lemma~\ref{lem:conv:hatg} suggests that the decomposition~\eqref{eq:decomposition:bis} could be
\begin{align} 
	 \varepsilon_{n,t}^{(0,h,Y)} - \langle \lambda_\sigma , \bfH h \rangle L_t(Y)
	& =		
	 \varepsilon_{n,t}^{(0,h,Y)} - \langle \lambda_\sigma , \bfH h \rangle 
	\varepsilon_{n,t}^{(0,\bfH{g},Y)} 
	+  \langle \lambda_\sigma , \bfH h \rangle \Big(
	\varepsilon_{n,t}^{(0,\bfH{g},Y)} 
	-L_t(Y) \Big).
\end{align}
Let us look more closely to the first two terms of the right-hand-side.
It would be convenient to add and remove a term involving the function $\bfH h$. 
The idea behind that is to prove, by dealing with martingales or functions of one variable, 
that the statistics $ \varepsilon_{n,t}^{(0,h,Y)} $, $ \varepsilon_{n,t}^{(0,\bfH{h},Y)} := \frac1{\sqrt{n}}\sum_{k=0}^{\lfloor n t \rfloor-1} \bfH{h} (\sqrt{n}Y_{k/n})$, 
and  $\langle \lambda_\sigma, \bfH h \rangle \varepsilon_{n,t}^{(0,\bfH{g},Y)}$ 
have the same limit (i.e.~$\langle \lambda_\sigma , \bfH h \rangle L_t(Y)$ by Lemma~\ref{lem:conv:hatg}). 

A candidate for $\widehat{\mathcal{M}}^n_t$ is 
\begin{equation}
	\varepsilon_{n,t}^{(0,h,Y)} - \varepsilon_{n,t}^{(0,\bfH{h},Y)}
	= 
	\frac{1}{\sqrt{n}} \sum_{k=1}^{\lfloor n t \rfloor} 
	\left( h(\sqrt{n}(Y_{(k-1)/n}-r),\sqrt{n}(Y_{k/n}-r)) 
		- \bfH{h} (\sqrt{n}Y_{(k-1)/n}) \right)
\end{equation}
which can be easily seen to be a $(\cF_{\lfloor n t \rfloor/n})_{t\geq 0}$-martingale: 
$	
	\widehat{\mathcal{M}}^n_t = \frac{1}{\sqrt{n}} \sum_{k=1}^{\lfloor n t \rfloor} \epsilon^n_k
$
with
\begin{equation}
	\epsilon^n_k
	:=	
h(\sqrt{n}(Y_{(k-1)/n}-r),\sqrt{n}(Y_{k/n}-r)) 
			- \EE\!\left[h(\sqrt{n}Y_{(k-1)/n},\sqrt{n}Y_{k/n}) | \sqrt{n}Y_{(k-1)/n}\right].
\end{equation}
Hence, $\mathcal{V}^n_t$ is given by
\begin{equation}
	\mathcal{V}^n_t = \langle \lambda_\sigma , \bfH h \rangle \left( 
\varepsilon_{n,t}^{(0,\bfH{g},Y)} 
-L_t(Y) \right)
	+  \frac1{\sqrt{n}}\sum_{k=0}^{\lfloor n t \rfloor-1} \bfG{h} (\sqrt{n}Y_{k/n})
\end{equation}
where
\begin{equation} \label{eq:kappa}
	\bfG{h} := \bfH {h} - \langle \lambda_\sigma , \bfH h \rangle \bfH{g}
\end{equation}
with
$
	g(x,y) = |y| - |x|.
$

By Lemma~\ref{lem:conv:hatg}, to prove~\eqref{eq:conv:decomp:fake}, it remains to study the convergence of both $\widehat{\mathcal{M}}^n_t$
and $\cG^n_t:=\frac1{\sqrt{n}}\sum_{k=0}^{\lfloor n t \rfloor-1} \bfG{h} (\sqrt{n}Y_{k/n})$: 
$\sup_{s\in [0,t]} (|\widehat{\mathcal{M}}^n_s| + |\cG^n_s|) \convp[n\to\infty] 0$.
In Section~\ref{sec:aux:conv}, we provide some auxiliary convergence results which are useful to deal with these two terms.
In particular, Proposition~\ref{prop:th4.1b:j1} (see Remark~\ref{rem:prop:th4:1:j1}) is suitable for dealing with both convergences. 
We provide the details in Section~\ref{sec:proof:details}.

To conclude we give some rough ideas of what it means that $\cG^n_t$ vanishes and why it is so. Note that the function $\bfG{h}$ is such that $\langle \lambda_\sigma , \bfG{h} \rangle =0$ which suggests that {``averaging''} over $\mathbf{E}^Y_{h}$ is not as different from ``averaging'' over couples consecutive increments crossing the threshold and then consider the average (integral) of the function $\mathbf{E}^Y_{h}$ with respect to the stationary measure of the process.

\subsection{Auxiliary convergence results} \label{sec:aux:conv}

We now provide several convergence results towards 0 or towards the local time of statistics of the form~$\frac1{\sqrt{n}} \sum_{k=0}^{\lfloor n t \rfloor-1} g_n(\sqrt{n} Y_{\frac{k}n})$ under different assumptions on the sequence $(g_n)_{n}$.
Note that $\widehat{\mathcal{M}}^n_t$ and  $\cG^n_t$ are of the latter form.

The first auxiliary result is the generalization to the case of OBM of~\cite[Lemma~4.2]{j1}. 
\begin{lemma} \label{lem:j1:4.2}
Let $(g_n)_{n\in \NN}$ be a sequence of real functions satisfying that $\langle \lambda_\sigma , g_n \rangle=0$ and for all $x\in \RR$
\begin{equation} \label{eq:j1:4.2}
	\lim_{n\to\infty} \frac{g_n(x\sqrt{n})^2}{n}+ \frac{\|g_n^2\|_{1}}{\sqrt{n}} + \frac{\| g_n \|_{1,1}|g_n(x \sqrt{n})|\log(n)}{n} + \frac{\| g_n \|_{1,1} \|g_n\|_{1} \log(n)}{\sqrt{n}}
=0.
\end{equation}
Then
$\lim_{n\to\infty} \EE\!\left[\left(\frac{1}{\sqrt{n}}  \sum_{k=0}^{\lfloor n t \rfloor-1} g_n(\sqrt{n} Y_{\frac{k}n}) \right)^2\right]=0$
for all $t\in [0,1]$.
\end{lemma}
\begin{proof}
The proof is analogous to the one of Lemma~4.2 in \cite{j1}.
It consists in bounding from above 
$\EE\!\left[\left( \sum_{k=0}^{\lfloor n t \rfloor-1} g_n(\sqrt{n} Y_{\frac{k}n}) \right)^2\right]$
by 
\[g^2(\sqrt n Y) + \Gamma_t(n,g^2,Y_0) + 2 (|g(\sqrt n Y_0)| + \Gamma_t(n,|g|,Y_0)) \sup_{x\in \mathbb{R},s\in[0,t]} |\Gamma_s(n,g,x)|\] 
and then using the different bounds of Lemma~\ref{lem:j1:3.3} to conclude. 
\end{proof}

The following propositions correspond to Theorem~4.1 a) and b) in \cite{j1}. 
\begin{proposition} \label{prop:th4.1a:j1}
Let $g_n \colon \RR \to \RR$, $n\in \NN$, be a sequence of functions satisfying
$\lim_{n\to\infty} \|g_n\|_{1} =0$ 
and for all $x\in \RR$ it holds that
$\lim_{n\to\infty} \frac{1}{\sqrt{n}} g_n(\sqrt{n} x) =0$.
Then for all $t\in (0,1]$ it holds
\begin{equation}\lim_{n\to\infty} \EE\Big[ \sup_{s\in [0,t]} \big|n^{-1/2}  \sum_{k=0}^{\lfloor n s \rfloor-1} g_n(\sqrt{n} Y_{\frac{k}n}) \big|\Big]=0.
\end{equation}
\end{proposition}
\begin{proof}
Notice that the integrand is bounded by $n^{-1/2}|g_n|(\sqrt{n}Y_0) + n^{-1/2} \Gamma_t(n,g_n,Y_0)$. The conclusion follows from the first statement in Lemma~\ref{lem:j1:3.3}. 
\end{proof}

\begin{proposition} \label{prop:th4.1b:j1}
Let $g_n \colon \RR \to \RR$, $n\in \NN$, be a sequence of functions satisfying \eqref{eq:j1:4.2}
and there exists $\lambda\in \RR$ such that $\lim_{n\to\infty} \langle \lambda_\sigma , g_n \rangle=\lambda$.
Then for all $t\in [0,1]$ it holds that 
\begin{equation}
	\frac1{\sqrt{n}}  \sum_{k=0}^{\lfloor n t \rfloor-1} g_n(\sqrt{n} Y_{\frac{k}n})  \convp[n\to\infty] \lambda L_t(Y).
\end{equation}
If in addition $\sup_{n\in \NN} \|g_n\|_{1} <\infty$ 
then 
\begin{equation}
	\sup_{s\in [0,1]} \Big| n^{-1/2}  \sum_{k=0}^{\lfloor n s \rfloor-1} g_n(\sqrt{n} Y_{\frac{k}n}) - \lambda L_s(Y)\Big|\convp[n\to\infty] 0.
\end{equation}
\end{proposition}
\begin{proof}
Let us set the sequence $f_n:=g_n - \langle \lambda_\sigma , g_n \rangle \bfH{g}$ with 
$g (x,y) := |y| - |x|$.
Note that Lemma~\ref{lem:conv:hatg} ensures that 
$\langle \lambda_\sigma , \bfH{g} \rangle=1$ and that $\frac{1}{\sqrt n}  \sum_{k=0}^{\lfloor n t \rfloor-1} \bfH{g}(\sqrt{n} Y_{\frac{k}n})  \convp[n\to\infty] L_t(Y)$.
Hence $\langle \lambda_\sigma , f_n \rangle=0$ and one can easily show that $f_n$ satisfies \eqref{eq:j1:4.2}.
Lemma~\ref{lem:j1:4.2} yields the result.
The additional statement is the same as \cite[Theorem~4.1]{j1}. It consists in observing that, up to taking a subsequence 
$\langle \lambda_\sigma , |g_{n}| \rangle$ converges and then applying the first statement to the sequences $g_{n}^\pm(x)=\max\{\pm g_{n}(x),0\}$. The fact that $\frac1{\sqrt{n}}  \sum_{k=0}^{\lfloor n t \rfloor-1} g_n^\pm(\sqrt{n} Y_{\frac{k}n}) $ are increasing and converge to a continuous (in time) limit implies the convergence in probability locally uniformly in time.
\end{proof}

\begin{remark}[Proposition~\ref{prop:th4.1a:j1} and \ref{prop:th4.1b:j1} for a constant sequence of functions] \label{rem:prop:th4:1:j1}
Let $f \in \Lone$ such that for all $x\in \RR$ it holds that $\lim_{n\to\infty}\frac{f(\sqrt{n}x)}{\sqrt{n}}=0$ (e.g.~$f\in \Lgb{0}$).
Then Proposition~\ref{prop:th4.1a:j1} states that if $\|f\|_1=0$ then for all $t\in (0,1]$ :
\begin{equation}
	\lim_{n\to\infty} \EE\Big[\sup_{s\in [0,t]}|n^{-1/2}  \sum_{k=0}^{\lfloor n s \rfloor-1} f(\sqrt{n} Y_{\frac{k}n}) |\Big]=0.
\end{equation}
And Proposition~\ref{prop:th4.1b:j1} states that if 
$f^2\in \Lone$ and $f\in \Loneg{1}$ (e.g.~$f\in\Lgb{1}$)
then 
\begin{equation}
	\sup_{s\in [0,1]} \Big| n^{-1/2}  \sum_{k=0}^{\lfloor n s \rfloor-1} f(\sqrt{n} Y_{\frac{k}n}) - \langle \lambda_\sigma , f \rangle L_s(Y)\Big|\convp[n\to\infty] 0.
\end{equation}
\end{remark}

\subsection{Proof of the consistency~\eqref{eq:conv:decomp:fake}} \label{sec:proof:details}

We wish to show $\sup_{s\in [0,t]} (|\widehat{\mathcal{M}}^n_s| + |\cG^n_s|) \convp[n\to\infty] 0$, 
where $\widehat{\mathcal{M}}^n_s$ and $\cG^n_s$ have been introduced in Section~\ref{app:decomp:fake}.
The next result, Lemma~\ref{lem:prop:kappa}, ensures that $\bfG{h}$ satisfies the assumptions in Remark~\ref{rem:prop:th4:1:j1}, thus $\sup_{s\in [0,t]} |\cG^n_s|\convp[n\to\infty] 0$.
For $\widehat{\mathcal{M}}^n_\cdot$, the proof is more elaborate. It is sketched after Lemma~\ref{lem:prop:kappa}.

\begin{lemma} \label{lem:prop:kappa}
Let $\gamma\in [0,\infty)$, $h \colon \mathbb R^2 \to \mathbb R$ be a measurable function such that $\bfH{h}$ is well defined and $\bfH{h}\in \Lgb{\gamma}$ (e.g.~$h\in \Lg{\gamma}$ by Lemma~\ref{lem:link:theorems}), and let
$\bfG{h}$ be the function defined in~\eqref{eq:kappa}.
Then $\langle \lambda_\sigma , \bfG{h} \rangle=0$ and $\bfG{h} \in \Lgb{\gamma}$.
\end{lemma}
\begin{proof}
Throughout this proof let $K_\sigma =\frac1{\min\{\sigma_-^2,\sigma_+^2\}}\frac{2 \sigma_- \sigma_+}{\sigma_-+\sigma_+} \in (0,\infty) $.
First note that the fact that
$
	q_\sigma(1,x,y) \leq K_\sigma \frac{1}{\sqrt{2\pi}} e^{-\frac{(x-y)^2}{2}}
$
implies that we reduce to the case of Brownian motion and so, some computations show:
\begin{equation}
	|\bfH{g}(x)| 
	\leq K_\sigma \frac{1}{\sqrt{2\pi}}  \int_{-\infty}^{\infty} ||y|-|x|| e^{-\frac{(x-y)^2}{2}} \vd y
	\in \Lgb{\alpha}
\end{equation}
for all $\alpha \geq 0$.
Moreover $\bfH h \in \Lgb{\gamma}$, hence in particular it holds $\langle \lambda_\sigma , \bfH h \rangle \leq \|\bfH h\|_{1}<\infty$. 
Therefore $|\bfG{h}| \leq |\bfH h| + \|\bfH h\|_{1} |\bfH g| \in \Lgb{\gamma}$ and so $\bfG{h}\in \Lgb{\gamma}$.

It remains to prove that $\langle \lambda_\sigma , \bfG{h} \rangle=0$.
This follows from the fact that $\langle \lambda_\sigma , \bfG{h} \rangle= \langle \lambda_\sigma , \bfH h \rangle (1-\langle \lambda_\sigma , \bfH{g} \rangle)$ and
$\langle \lambda_\sigma , \bfH{g} \rangle=1$ by~Lemma~\ref{lem:conv:hatg}.
\end{proof}

To conclude, we provide the arguments to deal with 
$	
	\widehat{\mathcal{M}}^n_t = \frac{1}{\sqrt{n}} \sum_{k=1}^{\lfloor n t \rfloor} \epsilon^n_k
$. 
By~\cite[Lemma~9]{GCJ93}, to prove that $\widehat{\mathcal{M}}^n_t$
converges to 0 in probability, it suffices to show that 
\begin{equation}
	\frac1n \sum_{k=1}^{\lfloor n t \rfloor} \EE\!\left[ (\epsilon^n_k)^2 | \cF_{(k-1)/n}\right] \convp[n\to\infty] 0.
\end{equation}
(Proposition~\ref{th:j2} provides the asymptotic behavior as well so it requires more restrictive assumptions.)
Since
\[
	\EE\!\left[ (\epsilon^n_k)^2 | \cF_{(k-1)/n}\right] 
	\leq 
	\EE\!\left[h(\sqrt{n}Y_{(k-1)/n},\sqrt{n}Y_{k/n})^2 | \sqrt{n}Y_{(k-1)/n}\right] 
	= 
	\bfH{h^2} (\sqrt{n}Y_{(k-1)/n}),
\]
it suffices to prove that 
\begin{equation}
	 \frac1{\sqrt n} \left( \frac1{\sqrt n} \sum_{k=1}^{\lfloor n t \rfloor} \bfH{h^2} (\sqrt{n}Y_{(k-1)/n}) \right) 
\convp[n\to\infty] 0.
\end{equation}
And this follows by applying the results of the following section, in particular by applying Proposition~\ref{prop:th4.1b:j1} (see Remark~\ref{rem:prop:th4:1:j1}) to $\bfH{h^2} \in \Lgb{1}$. 
The proof of~\eqref{eq:conv:decomp:fake} is thus completed.

\section{Proof of the key Proposition~\ref{th:j2:2}} \label{sec:proof:j1}

\begin{figure}[t]
\begin{center}
\begin{tikzpicture}

\fill[transpblue, rounded corners] (-3,2.5) rectangle (3,5);
\fill[transpred, rounded corners] (-6,0.5) rectangle (7,-3.5);
\fill[trared, rounded corners] (-2,-0.5) rectangle (6.2,-1.8);

\draw[blue, opacity=0.3] (-3,4.2) -- (3,4.2);
\path(-1.5,3.7) node{Lemma~\ref{lem:M:mart:2}};
\path(1.5,3.7) node{Lemma~\ref{lem:holderL}};
\draw[blue, opacity=0.3] (0,5) -- (0,2.5);
\path(-1.5,3) node{Lemma~\ref{lem:M:mart}};
\path(1.5,3) node{Lemma~\ref{lem:j1:5.2}};

\path[blue](-5.5,4.5) node{{\textbf{Decomposition~\eqref{eq:decomposition}} into }};
\path[red](-1.5,4.5) node{{martingale}};
\path[blue](1.5,4.5) node{{vanishing}};
\path[blue](3.7,4.5) node{{term}};

\path(-1,1.5) node{Lemma~\ref{lem:j1:4.2}};
\draw[-stealth, line width=0.3mm](-0.7,1.8) -- (1,2.7);

\path(3.2,1.5) node{Lemma~\ref{lem:j1:3.3}};
\draw[-stealth, line width=0.3mm] (3,1.2) -- (3,0.3);
\draw[-stealth, line width=0.3mm] (2.3,1.5) -- (-0.1,1.5);
\draw[-stealth, line width=0.3mm] (0,-0.6) -- (2.7,1.2);

\draw[red, opacity=0.3] (-6,-2.5) -- (7,-2.5);

\path[red](-4.8,-3) node{Item~\ref{th:j2:item1:2}};
\path(-4.8,0) node{Lemma~\ref{lem:conv:hatg}};
\draw[-stealth, line width=0.3mm] (-4.8,0.2) -- (-1.8,2.7);
\draw[red, opacity=0.3] (-3.5,0.5) -- (-3.5,-3.5);

\path[red](-1,-3) node{Item~\ref{th:j2:item2:2}};
\path(-1.8,0) node{Proposition~\ref{prop:th4.1b:j1}};
\draw[-stealth, line width=0.3mm] (-1,1.2) -- (-1.5,0.3);

\path(-0.3,-0.8) node{Lemma~\ref{lem:semig:estim}};
\draw[-stealth, line width=0.3mm](0,-0.6) -- (1.3,2.7);
\path(2,-1.1) node{Lemma~\ref{lem:prop:kappa}};
\path(4.3,-1.3) node{Lemma~\ref{lem:prop:G}};
\draw[-stealth, line width=0.3mm] (1.5,-0.9) -- (1.5,2.7);

\draw[red, opacity=0.3] (1,0.5) -- (1,-0.5);
\draw[red, opacity=0.3] (1,-1.8) -- (1,-3.5);

\path(-1.8,-2.2) node{Lemma~\ref{lem:cQ:conv}};
\draw[-stealth, line width=0.3mm] (-0.3,-1) -- (-1.6,-2);
\draw[-stealth, line width=0.3mm] (1,-1.2)-- (-0.8,-2);


\path[red](2.8,-3) node{Item~\ref{th:j2:item3:2}};
\path(2.8,0) node{Proposition~\ref{prop:th4.1a:j1}};
\draw[-stealth, line width=0.3mm](2.7,0.2) -- (1.7,2.7);
\draw[red, opacity=0.3] (4.5,0.5) -- (4.5,-0.5);
\draw[red, opacity=0.3] (4.5,-1.8) -- (4.5,-3.5);

\path[red](5.8,-3) node{Item~\ref{th:j2:item4:2}};

\end{tikzpicture}
\end{center}
\captionsetup{singlelinecheck=off}
\caption[Proof Prop map]{Map of the proof of Proposition~\ref{th:j2:2}.}\label{ProofMapProp}
In the appendix we introduce many auxiliary results. This map show how they intervene in the proof of Proposition~\ref{th:j2:2} and of other results of the appendix. 
\end{figure}
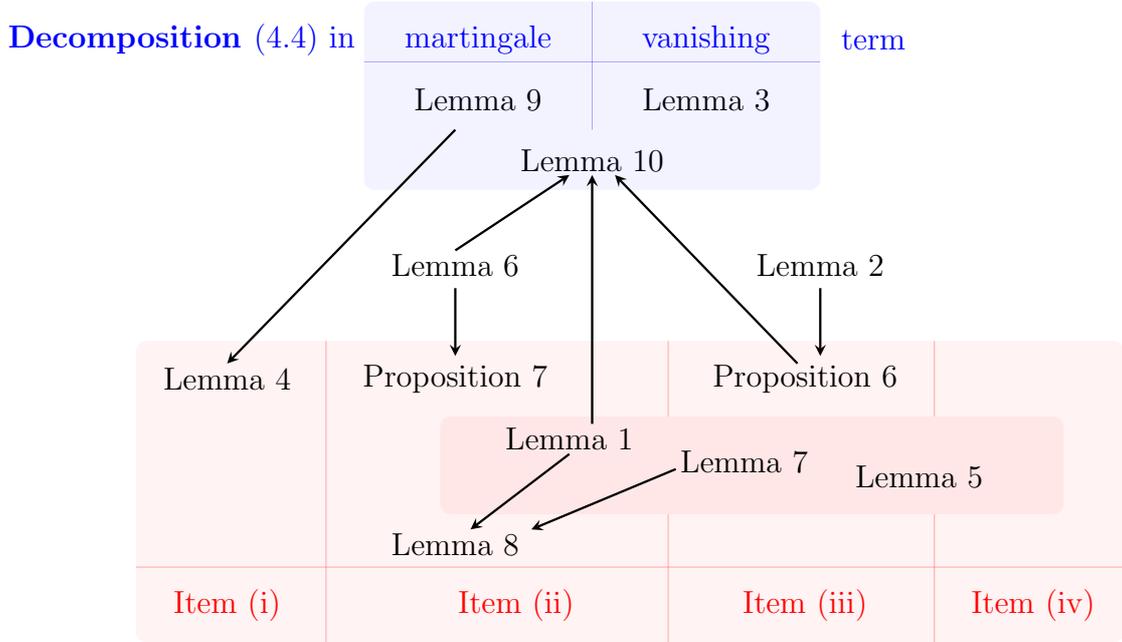

In this section we prove Proposition~\ref{th:j2:2} which was stated in Section~\ref{sec:proof:summary}.
The section is organized as follows: We split the proof of Proposition~\ref{th:j2:2} into three parts. The first part, in Section~\ref{sec:decomposition}, consists in finding the decomposition~\eqref{eq:decomposition} into a sum of martingale and a vanishing term. 
We deal with the vanishing term in Section~\ref{sec:th:j2:vanish}. Finally, in Section~\ref{sec:th:j2:hp}, we demonstrate that the martingale part satisfies Items~\ref{th:j2:item1:2}-\ref{th:j2:item4:2} of Proposition~\ref{th:j2:2}.
In Section~\ref{sec:finite:K}, we prove that the constant $K_h^Y$ defined by~\eqref{eq:k:th} is finite.
Figure~\ref{ProofMapProp} show how the results of the paper intervene in the proof of Proposition~\ref{th:j2:2}.

In this section, $(Y_t)_{t\in [0,1]}$ is a standard OBM.

Let us recall that by assumption $\gamma >3$ and $h\in \Lg{\gamma}$. 
By Lemma~\ref{lem:link:theorems}, $\bfH{f}$ and $\bfH{f^2}$ are well defined and belong to $\Lgb{\gamma}$.
For most of the following definitions and proof steps $\gamma\geq 0$ suffices. For the sake of precision, we specify when it does not.
%

\subsection{The decomposition into a sum of martingale and vanishing terms} \label{sec:decomposition}

In this section we determine the terms of the decomposition in equation~\eqref{eq:decomposition} into a sum of a vanishing term $n^{1/4}\mathcal{V}^n$ and a martingale $\mathcal{M}^n$.
Note that  for all $n\in \NN$, $t\in [0,1]$ :
\begin{equation} 
\begin{split}
	& \varepsilon_{n,t}^{(0,h,Y)}  - \langle \lambda_\sigma , \bfH h \rangle L_t(Y)
	 =		 
	M^{n,1}_t + \cG^n_t + \langle \lambda_\sigma , \bfH h \rangle \left(  L_{\lfloor n t \rfloor/n}(Y) - L_t(Y) \right)
\end{split}
\end{equation}
where $M^{n,1}$ and $\cG^n$ are the processes satisfying for all $t\in [0,1]$ that
\begin{equation}  
\begin{split}
	M^{n,1}_t 
	& = \varepsilon_{n,t}^{(0,h,Y)} - \varepsilon_{n,t}^{(0,\bfH h,Y)} + \langle \lambda_\sigma , \bfH h \rangle \left(  \varepsilon_{n,t}^{(0,\bfH g,Y)}	- L_{\frac{\lfloor n t \rfloor}n}(Y)\right)
\end{split}
\end{equation}
and 
$
	\cG^n_t =  \frac1{\sqrt{n}} \sum_{k=0}^{\lfloor n t \rfloor-1} \bfG{h}( \sqrt{n} Y_{\frac{k}n})	
$ 
($\bfG{h}$ is defined in \eqref{eq:kappa} and $g(x,y) := |y| - |x| $).
We used here the abuse of notation $\varepsilon_{n,t}^{(0,\bfH h,Y)}$ and $\varepsilon_{n,t}^{(0,\bfH g,Y)}$ as in Appendix~\ref{app:prop:lmt:2}.

This decomposition is not the same as the one proposed in Section~\ref{app:decomp:fake}.
Indeed, Lemma~\ref{lem:holderL} suggests that $L_{\lfloor n t \rfloor/n}(Y) - L_t(Y) $ will contribute to $\mathcal{V}^n_t$, so we split into two terms the quantity 
$
	 \langle \lambda_\sigma , \bfH h \rangle  \Big( \varepsilon_{n,t}^{(0,\bfH g,Y)} - L_t(Y)\Big).
$ 
Moreover, it has been shown in Section~\ref{sec:proof:details} that $\sup_{s\in [0,t]} |\cG^n_s|\convp[n\to\infty] 0$,
but the convergence speed have not been studied. In Section~\ref{sec:decomposition} we decompose $\cG^n_t $ into two terms, one contributing to $\mathcal{M}^n_t $ and the other to $\mathcal{V}^n_t $.
More precisely, $\cG^n= (\cG^n - M^{n,2} ) +M^{n,2}$
where
\begin{equation} \label{eq:W}
\begin{split}
& M^{n,2}_t  
	= \frac{1}{\sqrt{n}} \sum_{k=1}^{\lfloor n t \rfloor} \sum_{j=0}^{\lfloor n^\frac14\rfloor} \left(\EE\!\left[ \bfG{h}(\sqrt{n} Y_{(j+k)/n})| \cF_{k/n}\right]- \EE\!\left[ \bfG{h}(\sqrt{n} Y_{(j+k)/n})| \cF_{(k-1)/n}\right]\right).
\end{split}
\end{equation}
And so it holds for all $n\in \NN$, $t\in [0,1]$ that
\begin{equation} 
\begin{split}
	& \varepsilon_{n,t}^{(0,h,Y)}  - \langle \lambda_\sigma , \bfH h \rangle L_t(Y)
	 =		 
	M^{n,1}_t + M^{n,2}_t  + \cG^n_t - M^{n,2}_t + \langle \lambda_\sigma , \bfH h \rangle \left(  L_{\lfloor n t \rfloor/n}(Y) - L_t(Y) \right).
\end{split}
\end{equation}

We first show, in Lemma~\ref{lem:M:mart}, that $M^{n,1}_t$ is a good candidate for contributing to the martingale part of the decomposition.

\begin{lemma}  \label{lem:M:mart}
	$M^{n,1}$ is a martingale with respect to the filtration $(\cF_{\lfloor n t \rfloor/n})_{t\in [0,1]}$. 
\end{lemma}
\begin{proof}
Throughout this proof let $A^n$ and $B^n$ be the processes given by
\begin{equation}
	A^n_t 
	= \frac1{\sqrt{n}}  \sum_{k=0}^{\lfloor n t \rfloor-1} \bfH{g}(\sqrt{n} Y_{\frac{k}n}) - L_{\frac{\lfloor n t \rfloor}n}(Y) 
\text{ and }
	B^n_t = M^{n,1}_t - \langle \lambda_\sigma , \bfH h \rangle A^n_t.
\end{equation}
So
$ B^n_t = \frac1{\sqrt{n}} \sum_{k=0}^{\lfloor n t \rfloor-1} \left( h(\sqrt{n} Y_{\frac{k}{n}}, \sqrt{n} Y_{\frac{k+1}{n}} ) - \bfH h( \sqrt{n} Y_{\frac{k}n}) \right).
$
It suffices to show the martingale property for $A^n$ and $B^n$.
Let $t\in [0,1]$ be fixed. 
The martingale property for $A^n$ is an immediate consequence of Lemma~\ref{lem:conv:hatg}.
Let us explicit the case of the process $B^n$.
For all $j\in \{0,\ldots, \lfloor n t \rfloor-1\}$ 
it can be easily shown that
\begin{equation} 
\begin{split}
	\EE\!\left[  n^{\frac12} B^{n}_t | \cF_{\frac j n} \right]
 	& =
	\sum_{k=j}^{\lfloor n t \rfloor-1} \EE\!\left[  h(\sqrt{n} Y_{\frac{k}{n}}, \sqrt{n} Y_{\frac{k+1}{n}} ) - \bfH h( \sqrt{n} Y_{\frac{k}n}) | \cF_{\frac j n} \right] + B^{n}_{\frac j n}
= B^{n}_{\frac j n}.
\end{split}
\end{equation}
because for all $k \in \{j,\ldots, \lfloor n t \rfloor-1\}$ it holds that
\begin{equation} 
\begin{split}
	& \EE\!\left[  h(\sqrt{n} Y_{\frac{k}{n}}, \sqrt{n} Y_{\frac{k+1}{n}} ) - \bfH h( \sqrt{n} Y_{\frac{k}n}) | \cF_{\frac j n} \right]
	= \EE\!\left[  \EE\!\left[  h(\sqrt{n} Y_{\frac{k}{n}}, \sqrt{n} Y_{\frac{k+1}{n}} ) - \bfH h( \sqrt{n} Y_{\frac{k}n}) | \cF_{\frac k n} \right] | \cF_{\frac j n} \right]
	\\
	&
	= \EE\!\left[ \int_{-\infty}^{\infty} h(\sqrt{n} Y_{\frac{k}{n}}, \sqrt{n} y ) q_\sigma(1/n, \sqrt{n} Y_{\frac{k}{n}}, y)  \vd y - \bfH h( \sqrt{n} Y_{\frac{k}n})  | \cF_{\frac j n} \right] 
	= \EE\!\left[ 0  | \cF_{\frac j n} \right] =0.
\end{split}
\end{equation}
The proof is thus completed. 
\end{proof}

Definition~\eqref{eq:W} implies that $M^{n,2}$ is a martingale as well.
\begin{lemma}  \label{lem:M:mart:2}
	$M^{n,2}$ is a martingale with respect to the filtration $(\cF_{\lfloor n t \rfloor/n})_{t\in [0,1]}$. 
\end{lemma}

\begin{proof}
Using \eqref{eq:cQ} we rewrite $M^{n,2}$ in \eqref{eq:dec:m3} as
\begin{equation} \label{eq:W}
\begin{split}
& M^{n,2}_t  
	= \frac{1}{\sqrt{n}} \sum_{k=1}^{\lfloor n t \rfloor} \sum_{j=0}^{\lfloor n^\frac14\rfloor} \left(\EE\!\left[ \bfG{h}(\sqrt{n} Y_{(j+k)/n})| \cF_{k/n}\right]- \EE\!\left[ \bfG{h}(\sqrt{n} Y_{(j+k)/n})| \cF_{(k-1)/n}\right]\right)
\end{split}
\end{equation}
which shows that $M^{n,2}$ is a martingale with respect to the filtration $(\cF_{\lfloor n t \rfloor/n})_{t\in [0,1]}$.
\end{proof}

Lemma~\ref{lem:M:mart} and Lemma~\ref{lem:M:mart:2} identify the $(\cF_{\lfloor n t \rfloor/n})_{t\geq 0}$-martingale $\mathcal{M}^n$ 
and also the candidate to be the vanishing term with rate of order at least $1/4$ (denoted by $\mathcal{V}^n$):
\begin{equation} \label{eq:mv}
	\begin{split}
	&\mathcal{M}^n_t  := n^{1/4} \left( M^{n,1}_t + M^{n,2}_t \right) 
	\ \text{and\/}  \
	\\
	&
	\mathcal{V}^n_t  := \cG^n_t - M^{n,2}_t + \langle \lambda_\sigma , \bfH h \rangle \left(  L_{\frac{\lfloor n t \rfloor}n}(Y) - L_t(Y) \right).
\end{split}
\end{equation}
To prove Proposition~\ref{th:j2:2}, it remains to show that $\sup_{t\in [0,1]} n^{1/4} |\mathcal{V}^n_t| \convp[n\to\infty] 0$ and that $\mathcal{M}^n_\cdot$ satisfies Items~\ref{th:j2:item1:2}-\ref{th:j2:item4:2}. 
This is done in Section~\ref{sec:th:j2:vanish} and Section~\ref{sec:th:j2:hp}.

\subsection{Finiteness of the variance} \label{sec:finite:K}

In this section we provide the proof of the finiteness of the constant $K_h^Y$ defined by~\eqref{eq:k:th}. 
The proof is actually contained in Section~\ref{sec:th:j2:hp}, but we choose to show it separately.
Instead, we do not prove separately non-negativity of the constant.

The expression of $K_h^Y$ depends on a function, $\mathcal{P}_h^Y$ which is here expressed with the use of some auxiliary functions.
We first introduce and prove some properties of these functions that are also used to provide an alternative expression for $M^{n,2}$.
For all $i,j\in \{0,1\}$, $\eta\in [0,\infty)$ let 
\begin{equation} \label{eq:cQ}
	\cQ^{(\eta)}_{n, i , j} := \sum_{k=i}^{\lfloor n^{\eta}\rfloor + j } Q^\sigma_k \bfG{h} \quad \text{ and } \quad \cQ_{n, i , j}:=\cQ^{(\frac14)}_{n, i , j}
\end{equation}
where $Q^\sigma$ is the semigroup of the OBM given in~\eqref{eq:def:semig}. 
These quantities may be seen as the aggregated action of the semigroup on the function $\bfG{h}$ up to different times of order $n^\eta$. We reduce to consider only some precise values of $\eta$, in particular $\eta=1/4$ and $\eta=3/4$. 
(The results in the next sections would restrict the parameter choices, and these precise choices would respect the restrictions.) 

%
Note that the process $M^{n,2}$ satisfy for all $t\in [0,1]$ that
\begin{equation} \label{eq:dec:m3}
	M^{n,2}_t   =  \frac1{\sqrt{n}} \sum_{k=1}^{\lfloor n t \rfloor}  \left(\cQ_{n,0,0}(\sqrt{n} Y_{\frac k n}) - \cQ_{n,1,1}(\sqrt{n} Y_{\frac{k-1}{n}})\right).
\end{equation}
%

%
%
The following lemma relates the sequence $(\cQ_{m,0,0})_{m\in \NN}$ with $\mathcal{P}_{h}^{Y} $ in~\eqref{eq:cQ:limit} (with $\beta=0$ because $Y$ is an OBM).
\begin{lemma} \label{lem:cQ:conv}
Pointwise	
$\lim_{m\to\infty} \cQ_{m,0,0} = \mathcal{P}_{h}^{Y} $. 
Moreover, $x\mapsto \mathcal{P}_{h}^{Y}(x)/(1+|x|) $ is bounded.
\end{lemma}
\begin{proof}
 Let $x\in \mathbb{R}$.
	Note that $\mathcal{P}_{h}^{Y}$ in~\eqref{eq:cQ:limit} for OBM rewrites as
	\[
		\mathcal{P}_{h}^{Y}(x)=\sum_{k=0}^{\infty } \int_{-\infty}^{\infty} p_X(k, x, y) \bfG{h}(y) \vd y
		= \sum_{k=0}^{\infty } Q^\sigma_k \bfG{h}(x) .
	\]
Hence there exists a positive constant $K\in (0,\infty)$ such that
\[ 
	\left| \cQ_{m,0,0}(x) - \mathcal{P}_{h}^{Y}(x) \right| 
	\leq 
	\sum_{k=\lfloor m^{1/4}\rfloor+1}^{\infty} \left|  Q^\sigma_k \bfG{h}(x) \right| \leq K (1+|x|) \sum_{k=\lfloor m^{1/4}\rfloor+1}^{\infty} k^{-3/2}
\]
where the last inequality follows from 
Item~\ref{item2:semig:estim} in Lemma~\ref{lem:semig:estim} 
and Lemma~\ref{lem:prop:kappa}. Indeed, Lemma~\ref{lem:prop:kappa} establishes that $\bfG{h}\in \Lgb{\gamma}$ for all $\gamma>0$ (in particular for $\gamma =1$ and $\gamma =2$).
The fact that $\sum_{k=\lfloor m^{1/4}\rfloor+1}^{\infty} k^{-3/2}$ converges to $0$ as $m\to\infty$ proves the first statement.
Similarly, Item~\ref{item2:semig:estim} in Lemma~\ref{lem:semig:estim} 
proves that there exist constants $K,K'\in (0,+\infty)$ such that
$
	|\mathcal{P}_{h}^{Y}(x)| 
	\leq  |\bfG{h}(x)| + K (1+|x|) \sum_{k=1}^{\infty} k^{-3/2} 
	\leq K' (1+|x|).
$
The proof is thus completed.
\end{proof}

We are now ready to prove the finiteness of $K_h^Y$. 

\begin{lemma} \label{lem:finite:K}
	The quantity $K_h^Y$ defined by~\eqref{eq:k:th} is a finite constant. 
\end{lemma}
\begin{proof}
Let us recall that
\begin{equation}
{\scriptsize
\begin{split}
    & K_h^Y
	 = \langle \mu_Y , \bfH{h^2} + 2 \bfH{h,\mathcal{P}_{h}^Y} \rangle  
			+ \frac{ 2 \sigma_- \sigma_+}{\sigma_-+\sigma_+} \frac{8}{3 \sqrt{2\pi}} (\langle \mu_Y , \bfH{h} \rangle)^2 
			\\
	& \quad - 2 \sqrt{\frac{2}{\pi}} \frac{ 2 \sigma_- \sigma_+}{\sigma_-+\sigma_+} \langle \mu_Y , \bfH{h} \rangle
\int_{-\infty}^\infty \left( e^{-\frac{y^2}{2(\sigma(y))^2}} -\sqrt{2\pi} \frac{|y|}{\sigma(y)} \Phi\left(-\frac{|y|}{\sigma(y)}\right)\right) \mathcal{P}_{h}^Y(y) \mu_Y(\!\vd y)
	\\
	& \quad -
	2 \langle \mu_Y , \bfH{h} \rangle 
	\left(\frac{ 2 \sigma_- \sigma_+}{\sigma_-+\sigma_+}\right)^{\!2}
\int_{-\infty}^{\infty} \int_0^1 \int_{-\infty}^{\infty} \tfrac{|x| e^{-\frac{x^2}{2(\sigma(x))^2}} \Phi\left(-\frac{|y|}{\sigma(y)}\right)}{\sqrt{2\pi}\sigma(x)}  \sqrt{\frac{1}t-1} h(x \sqrt{t}, y \sqrt{1-t}) \mu_Y(\!\vd y) \vd t \mu_Y(\!\vd x),
\end{split}
}\end{equation}
where
$\Phi$ is the cumulative distribution function oh a standard Gaussian random variable,
and $\mathcal{P}_{h}^Y$ has been studied in Lemma~\ref{lem:cQ:conv}.
%
Note that $\lambda$-integrability is equivalent to $\mu_Y$-integrability.
Therefore if $f\in \Lgb{0}$ then $f$ is $\lambda$-integrable and $\langle \mu_Y , f \rangle$ is finite.
Since $\gamma\geq 1$, Lemma~\ref{lem:link:theorems} ensures that $f=\bfH{h},\bfH{h^2}\in \Lgb{1}$ which is contained in $\Lgb{0}$ by Remark~\ref{rem:spaces:relation}.
And also, in Lemma~\ref{lem:cQ:conv} it has been shown that $|\mathcal{P}_{h}^Y (y)| \leq K (1+|y|)$ for some constant $K\in (0,\infty)$,
so Lemma~\ref{lem:link:theorems} ensures that $x\mapsto \bfH{h, \mathcal{P}_{h}^Y}(x)/(1+|x|)$ belongs to $\Lgb{\gamma}$. Since $\gamma\geq 1$, we deduce that $\bfH{h, \mathcal{P}_{h}^Y}$ is $\lambda$-integrable.
Let us now consider the second last integral appearing in the expression of $K_h^Y$. 
Since $ 1- \sqrt{2\pi} \frac{|y|}{\sigma(y)} \Phi\left(-\frac{|y|}{\sigma(y)}\right) e^{\frac{y^2}{2(\sigma(y))^2}} \in [0,1]$ for all $y\in \mathbb{R}$,
we reduce to check that $e^{-y^2/2} \mathcal{P}_{h}^Y(\sigma(y)y)$ is $\lambda$-integrable. And this is the case because $e^{-y^2/2} \mathcal{P}_{h}^Y(\sigma(y)y) \leq K (1+|y|) e^{-y^2/2}$.
For the last integral we use that $h\in \Lg{\gamma}$ and exploit the inequality $h(x,y) \leq \bar{h}(x) e^{a|y-x|} \leq K e^{a|y-x|}$ for some constant $a,K\in (0,\infty)$. Then we deal with $\int_0^1 \sqrt{\frac1t -1} e^{a|\sqrt t x -\sqrt{1-t}y|} \vd t \leq e^{a|x|+a|y|} \int_0^1 \sqrt{\frac1t -1}  \vd t = e^{a|x|+a|y|}  \frac{\pi}2.$ 
By observing that $ x\mapsto e^{a|x|} \tfrac{ |x| e^{-\frac{x^2}{2(\sigma(x))^2}} }{\sqrt{2\pi}\sigma(x)} $  and $y \mapsto e^{a|y|}\Phi\left(-\frac{|y|}{\sigma(y)}\right)$ are $\lambda$-integrable, we complete the proof.
\end{proof}

To conclude the section we study some properties of the sequences $\cQ_{n,i,j}$ which intervene in the next sections.
The following facts are consequences of Lemma~\ref{lem:prop:kappa}, Lemma~\ref{lem:semig:estim} (which requires $\gamma \geq 2$), and the fact that $\sum_{k=1}^{n} \frac1k \leq 2 \log(n)$ for $n \geq 2$.
For every $\zeta\in [0,\gamma-1]$ and every $\eta \in (0,1)$, the fact that  $\langle \lambda_\sigma , \bfG{h} \rangle=0$ and Item~\ref{item2bis:semig:estim} in Lemma~\ref{lem:semig:estim} imply that for all $x\in \RR$ it holds that
\begin{equation} \label{eq:cQ:prop1}
\begin{split}
	|Q^\sigma_{n^{\eta}+1} \bfG{h} (x) | 
	& \leq  \left(1+\tfrac{|\sigma_-\sigma_+|}{\sigma_-+\sigma_+}\right)  K_\zeta n^{-\eta} \left(\tfrac1{1+|x n^{-\eta/2}/\sigma(x)|^\zeta} + \tfrac1{1+(|x|/\sigma(x))^\zeta}\right)
\end{split}
\end{equation}
and
\begin{equation} \label{eq:cQ:prop2}
	|\cQ^{(\eta)}_{n,1,0}(x)|+|\cQ^{(\eta)}_{n,1,1}(x)| \leq 2 K_\zeta \log(n) \left(\tfrac1{1+|x n^{-{\eta}/2}/\sigma(x)|^\zeta} + \tfrac1{1+|x/\sigma(x)|^\zeta}\right)
\end{equation}
for some $K_\zeta\in (0,\infty)$ depending also on $\eta$.
Hence \eqref{eq:cQ:prop2} with $\zeta=0$ and Item~\ref{item2:semig:estim} in Lemma~\ref{lem:semig:estim} imply that for all $x\in \RR$ it holds that
\begin{equation} \label{eq:cQ:prop3}
	|\cQ^{(\eta)}_{n,1,0}(x)|+|\cQ^{(\eta)}_{n,1,1}(x)| \leq 2 K \min{\{ \log(n) ,(1+|x|) \}},
\end{equation}
for some constant $K\in (0,\infty)$ depending on $\eta \in (0,1)$.
This and the fact that $\bfG{h}$ is bounded ensures that for some positive constant $K\in (0,\infty)$ it holds that
$
	|\cQ^{(\eta)}_{n,0,0}(x)|+|\cQ^{(\eta)}_{n,0,1}(x)| \leq 2 |\bfG{h}(x)| + |\cQ^{(\eta)}_{n,1,0}(x)|+|\cQ^{(\eta)}_{n,1,1}(x)| \leq 2 K (\log(n) +1)
$
and 
\begin{equation} \label{eq:cQ:prop4}
	|\cQ^{(\eta)}_{n,0,0}(x)|+|\cQ^{(\eta)}_{n,0,1}(x)| + |\cQ^{(\eta)}_{n,1,0}(x)|+|\cQ^{(\eta)}_{n,1,1}(x)| \leq 2 K ( \log(n) +1).
\end{equation}

\subsection{Dealing with the vanishing term} \label{sec:th:j2:vanish}

We now prove that $\mathcal{V}^n$ defined in \eqref{eq:mv}, satisfies that $\sup_{t\in [0,1]} n^{1/4} |\mathcal{V}^n_t| \convp[n\to\infty] 0$.
This is a consequence of~Lemma~\ref{lem:holderL} and Lemma~\ref{lem:j1:5.2}.

\begin{lemma} \label{lem:j1:5.2}
Let $h\in \Lg{3}$.
Then 
 $\sup_{s\in [0,1]} n^{1/4} |\cG^n_s-M^{n,2}_s| \convp[n\to\infty] 0$.
\end{lemma}
\begin{proof}
Using \eqref{eq:cQ} we rewrite $M^{n,2}$ in~\eqref{eq:W} as~\eqref{eq:dec:m3} and also as 
\begin{equation} \label{eq:W:2}
\begin{split}
& M^{n,2}_t  
	= \cG^n_t 
			- \frac{1}{\sqrt{n}}  \left( \cQ_{n,0,0}(\sqrt{n} Y_0) -  \cQ_{n,0,0}(\sqrt{n} Y_{\lfloor n t\rfloor/n}) \right)  
			-  \frac{1}{\sqrt{n}}  \sum_{j=0}^{\lfloor n t \rfloor -1} Q^\sigma_{\lfloor n^\frac14\rfloor +1} \bfG{h} (\sqrt n Y_{j/n}).
\end{split}
\end{equation}
%
%
Let 
\begin{equation}{m^{n}_t := \cG^n_t 
			- \frac1{\sqrt{n}} \left( \cQ^{(3/4)}_{n,0,0}(\sqrt{n} Y_0) -  \cQ^{(3/4)}_{n,0,0}(\sqrt{n} Y_{\lfloor n t\rfloor/n}) \right)  
			-  \frac1{\sqrt{n}} \sum_{j=0}^{\lfloor n t \rfloor -1} Q^\sigma_{\lfloor n^\frac34\rfloor +1} \bfG{h} (\sqrt n Y_{j/n}).
}
\end{equation}
By analogy, $m^n_t$ has an expression as~\eqref{eq:W} and so 
it is clear that $m^n$ is a martingale with respect to the filtration $(\cF_{\lfloor n t \rfloor/n})_{t\in [0,1]}$.
Therefore $n^{1/4} (M^{n,2}_t-m^n_t)$ as well.
Let us denote by $D^n_t := n^{1/4} (\cG^n_t - M^{n,2}_t)$ and $d^n_t :=  n^{1/4} (\cG^n_t - m^n_t)$.
Then $D^n-d^n = n^{1/4} (m^n-M^{n,2})$ is a martingale with respect to the filtration $(\cF_{\lfloor n t \rfloor/n})_{t\in [0,1]}$.
In this notation, the goal is to prove 
$\sup_{s\in [0,1]} |D^n_s|\convp[n\to\infty] 0$. 

{\em First step:\/} For every $t\in [0,1]$, let us show that $\lim_{n\to\infty} \EE\!\left[( D^n_t-d^n_t)^2\right]=0$ by demonstrating the stronger fact that $\lim_{n\to\infty} \EE\!\left[(D_t^n)^2+(d_t^n)^2\right]=0$.\\
Let $t\in [0,1]$ and $\eta \in \left\{ \frac14, \frac34 \right\}$ be fixed.
Inequality \eqref{eq:cQ:prop4} implies for all $\omega\in \Omega$ that
\begin{equation} \label{lem:j1:5.2.eq}
	\sup_{s\in [0,1]} n^{-\frac14} \left( \cQ^{(\eta)}_{n,0,0}(\sqrt{n} Y_0(\omega)) -  \cQ^{(\eta)}_{n,0,0}(\sqrt{n} Y_{\lfloor n s\rfloor/n}(\omega)) \right)
\xrightarrow[n\to\infty]{} 0,
\end{equation}
hence it holds also that 
$
	\lim_{n\to\infty} \EE\!\left[n^{-1/2} \left( \cQ^{(\eta)}_{n,0,0}(\sqrt{n} Y_0) -  \cQ^{(\eta)}_{n,0,0}(\sqrt{n} Y_{\lfloor n t\rfloor/n}) \right)^2 \right]=0.
$
Next observe that, since $\lambda_\sigma$ is a stationary measure, the sequences of functions $g_n^{(\eta)} := n^{1/4} Q^\sigma_{\lfloor n^{\eta}\rfloor +1} \bfG{h} $, $\eta \in \{\tfrac14,\tfrac34\}$ satisfy that $\langle \lambda_\sigma , g_n^{(\eta)} \rangle = n^{1/4} \langle \lambda_\sigma , \bfG{h} \rangle$ which is equal to 0 by Lemma~\ref{lem:prop:kappa}. This and inequality \eqref{eq:cQ:prop1} (with $\zeta=2$ since $h\in \Lg{\zeta+1}$) 
ensure that \eqref{eq:j1:4.2} holds. (A sufficient condition is $\eta \geq 1/4$.)\\
Hence, Lemma~\ref{lem:j1:4.2} shows that
$
	\lim_{n\to\infty} \EE\!\left[ n^{-1/2} \left(\sum_{j=0}^{\lfloor n t \rfloor -1} Q^\sigma_{\lfloor n^{\eta}\rfloor +1} \bfG{h} (\sqrt n Y_{{j}/{n}})\right)^2\right]=0.
$

{\em Second step:\/} It holds that $\sup_{s\in [0,1]} |D^n_s - d_s^n|\convp[n\to\infty] 0$.\\
This follows from \cite[Proposition~1.2]{aldous89} as a consequence of the previous step and the martingale property of $D^n - d^n$.

{\em Third step:\/} It holds that $\sup_{s\in [0,1]} |d_s^n| \convp[n\to\infty] 0$.
This follows from~\eqref{lem:j1:5.2.eq} and from applying Proposition~\ref{prop:th4.1a:j1} to the sequence $g_n:= n^{1/4} Q^\sigma_{\lfloor n^\frac34\rfloor +1} \bfG{h} $. (Here taking $n^{1/4} Q^\sigma_{\lfloor n^{\eta} \rfloor +1} \bfG{h}$ with $\eta>1/2$ would be sufficient.)
The assumptions are indeed satisfied: $\lim_{n\to\infty} \|g_n\|_{1} =0$ by Item~\ref{item1:semig:estim} in Lemma~\ref{lem:semig:estim}
 and in the first step we have proven that  $g_n$ satisfies \eqref{eq:j1:4.2} and in particular $\lim_{n\to\infty} \frac{g_n(\sqrt{n}x)}{\sqrt{n}}=0$.

Combining the two last steps yields the conclusion.
\end{proof}

\subsection{Dealing with the martingale term} \label{sec:th:j2:hp}
In this section we complete the proof of Proposition~\ref{th:j2:2}. The arguments are sometimes sketched because analogous to the ones in \cite[Section~6]{j1}.

Let $\gamma >3$ and $h\in \Lg{\gamma}$ be fixed.
By Definition~\ref{def:integral:gamma} of $\Lg{\gamma}$ there exist a non-negative function $\bar h \in \Lgb{\gamma}$ and a constant $a\in[0,\infty)$ such that $|h(x,y)| \leq \bar h(x) e^{a |y-x|}$. In this section $\bar h$ and $a$ are fixed.
\\
Let us also recall some notation: let $\bfH{}$ the functional in \eqref{eq:muX}, 
$Q^\sigma$ the semigroup in~\eqref{eq:def:semig}, 
$\bfG{h}$ in~\eqref{eq:kappa}, 
$\cQ$ in~\eqref{eq:cQ} and its limit $\mathcal{P}_{h}^{Y}$ in~\eqref{eq:cQ:limit} (see Lemma~\ref{lem:cQ:conv}), and $\mathbb{L}$ in~\eqref{function:G}.

For all $n\in \NN$ the $(\cF_{\lfloor n t\rfloor /n})_{t\in [0,1]}$-martingale $\cM^n$ in \eqref{eq:mv} rewrites as
$
	{ \cM^n_t =  \sum_{k=1}^{\lfloor n t \rfloor} \chi^n_k}
$
where
\begin{equation} \label{eq:chi:nk}
\begin{split}
	\chi^n_k 
	:= 
	n^{-\frac14}  
	& \left( h(\sqrt{n} Y_{(k-1)/n}, \sqrt{n} Y_{k/n} ) -   \langle \lambda_\sigma , \bfH h \rangle \, \sqrt{n}(L_{k/n}(Y) - L_{(k-1)/n}(Y)) \right)
	 \\
	+ \, & n^{-\frac14}  \left(\cQ_{n,0,0}(\sqrt{n} Y_{k/n}) - \cQ_{n,0,1}(\sqrt{n} Y_{(k-1)/n})
	\right)\!.
\end{split}
\end{equation}
Now it remains to prove Items~\ref{th:j2:item1:2}-\ref{th:j2:item4:2} in Proposition~\ref{th:j2:2}. 
%
%
%
\subsubsection{Proof of 
Item~\ref{th:j2:item1:2} in Proposition~\ref{th:j2:2} 
} 
It follows from scaling property \eqref{Yscaling} and Lemma~\ref{lem:conv:hatg}.
For all $n\in \NN$, $k\in \{1, \ldots, \lfloor n \rfloor \}$, the scaling property \eqref{Yscaling}, Lemma~\ref{lem:conv:hatg} and \eqref{eq:cQ} ensure that 
\begin{equation} \begin{split}
	& \EE\!\left[ \chi_k^n | \cF_{(k-1)/n}\right] 
	\\
	& = n^{-\frac14} \left( \bfG{h}(\sqrt{n} Y_{(k-1)/n}) + \int_{-\infty}^\infty \cQ_{n,0,0}(y) q_\sigma(1, \sqrt{n} Y_{(k-1)/n}, y) \vd y - \cQ_{n,0,1}(\sqrt{n} Y_{(k-1)/n})\right)
	\\
	& = n^{-\frac14} \left(Q^\sigma_1 \cQ_{n,0,0}(\sqrt{n} Y_{(k-1)/n}) - \cQ_{n,1,1}(\sqrt{n} Y_{(k-1)/n})\right)=0.
\end{split}
\end{equation}

\subsubsection{Proof of Item~\ref{th:j2:item2:2} in Proposition~\ref{th:j2:2}}
Let $t\in [0,1]$ be fixed.

{\it First step:}
It can be easily shown that
\begin{equation}
\begin{split}
	& \sqrt n \EE\!\left[ (\chi_k^n)^2 | \cF_{(k-1)/n}\right] 
	\\
	& =   f_{n}(\sqrt{n} Y_{(k-1)/n}) + (\langle \lambda_\sigma , \bfH h \rangle)^2 \Lm{1}{2}{\sqrt{n} Y_{(k-1)/n}}  - 2 \langle \lambda_\sigma , \bfH h \rangle  h_n(\sqrt{n} Y_{(k-1)/n})
\end{split}
\end{equation}
where $f_n$ and $h_n$ are given by 
$
f_n(x) := \bfH {h^2}(x)  
 	+ 2 \bfH{h, \cQ_{n,0,0}} (x) + g_n(x)
$ 
and
$
h_n(x):= 
	\Lm{h(x, \cdot)}{1}{x} + \Lm{\cQ_{n,0,0}}{1}{x},
$
with  
$ 
g_n(x):=Q^\sigma_1\left(  ( \cQ_{n,0,0} )^2 \right) (x)    - \left(\cQ_{n,0,1}(x)\right)^2
$.

In fact first note that, by \eqref{function:G:equal}, for all $k =1, \ldots, {\lfloor n t \rfloor}$ it holds that
\begin{equation}
	h_n(\sqrt{n} Y_{\frac{k-1}{n}})= \sqrt{n} \EE\Big[ (h(\sqrt{n} Y_{\frac{k-1}{n}}, \sqrt{n} Y_{\frac{k}{n}}) + \cQ_{n,0,0}(\sqrt{n} Y_{\frac{k}{n}}) ) (L_{\frac{k}{n}}(Y) - L_{\frac{k-1}n}(Y))  | \cF_{\frac{k-1}n}\Big] .
\end{equation}
Let us now consider $f_n$.
Clearly $f_n$ have to be the sum of all remaining terms and it has to have the desired form. 
This follows from the fact that the scaling property~\eqref{Yscaling} and equation~\eqref{eq:cQ} ensure that
$\EE\Big[ 
		(\cQ_{n,0,0}(\sqrt{n} Y_{k/n}))^2  | \cF_{(k-1)/n}\Big] 
= Q^\sigma_1 (
		(\cQ_{n,0,0}(\sqrt{n} Y_{k/n}))^2 ).
$
but also
$\EE\Big[ \cQ_{n,0,0}(\sqrt{n} Y_{k/n})  | \cF_{(k-1)/n}\Big] = \cQ_{n,1,1}(\sqrt n Y_{(k-1)/n})$
which, together with the definition of $\bfG{h}$ in \eqref{eq:kappa}, \eqref{function:G:equal}, and the fact that
$ \cQ_{n,1,1}(\sqrt{n} Y_{(k-1)/n})   =  \cQ_{n,0,1}(\sqrt{n} Y_{(k-1)/n}) -  \bfG{h}(\sqrt{n} Y_{(k-1)/n}) $
show that $\cQ_{n,0,1}(\sqrt{n} Y_{(k-1)/n})$ can be rewritten as 
\begin{equation}
\begin{split}
	\cQ_{n,0,1}(\sqrt{n} Y_{\frac{k-1}{n}})  
	=  
	 \EE\Big[ h(\sqrt{n} Y_{\frac{k-1}n}, \sqrt{n} Y_{\frac{k}{n}}) + \cQ_{n,0,0}(\sqrt{n} Y_{\frac{k}n}) -\langle \lambda_\sigma , \bfH h \rangle \sqrt{n} ( L_{\frac{k}n}(Y) - L_{\frac{k-1}n}(Y) ) | \cF_{\frac{k-1}n}\Big] .
\end{split}	
\end{equation}
The proof of the first step is thus completed.

{\it Second step:} 
It follows from applying Proposition~\ref{prop:th4.1b:j1} to the constant sequence of functions $\Lm{1}{2}{\cdot}$ that
$
	\frac1{\sqrt n}  \sum_{k=1}^{\lfloor n t \rfloor} \Lm{1}{2}{\sqrt{n} Y_{(k-1)/n}} \convp[n\to\infty]
	\frac{2 \sigma_- \sigma_+}{\sigma_-+\sigma_+} \frac{8}{3 \sqrt{2 \pi}}
 L_t(Y).
$
%
The fact that the assumptions of the proposition (i.e.~$\Lm{1}{2}{\cdot} \in \Lgb{2}$) are satisfied and 
$\langle \lambda_\sigma , \Lm{1}{2}{\cdot} \rangle = \frac{2 \sigma_- \sigma_+}{\sigma_-+\sigma_+} \frac{8}{3 \sqrt{2 \pi}}$
follows from the fact that $\EE\!\left[ \left( L_1(W)\right)^2\right]= 1$ 
and Lemma~\ref{lem:prop:G}.

In the two final steps, we want to apply Proposition~\ref{prop:th4.1b:j1} to the sequences $f_n$ and $h_n$.

{\em Third step:\/} We show that
$
	\frac1{\sqrt n}  \sum_{k=1}^{\lfloor n t \rfloor} f_n (\sqrt{n} Y_{(k-1)/n}) \convp[n\to\infty]
	\langle \lambda_\sigma , \bfH {h^2} + 2 \bfH{h, \mathcal{P}_{h}^{Y} } \rangle L_t(Y).
$
applying Proposition~\ref{prop:th4.1b:j1}.
To do so we check that the sequence $f_n$ satisfies \eqref{eq:j1:4.2} and that 
$
	\lim_{n\to\infty}\langle \lambda_\sigma , f_n \rangle=\langle \lambda_\sigma , \bfH {h^2} + 2 \bfH{h, \mathcal{P}_{h}^{Y} } \rangle
$.

The fact that $\lim_{n\to\infty} \langle \lambda_\sigma , g_n \rangle$=0 follows from the fact that $\lambda_\sigma$ is a stationary measure, inequality~\eqref{eq:cQ:prop1} with $\zeta=\gamma-1 >2$, and inequality~\eqref{eq:cQ:prop2}.
%
Note that
$
	\bfH{h, \cQ_{n,0,0}}= \bfH{h, \bfG{h}}+\bfH{h, \cQ_{n,1,0}}.
$
Lemma~\ref{lem:prop:kappa} (in particular the fact that $\bfG{h}$ is bounded) ensures that there exists a constant $K\in (0,\infty)$ 
such that
$
	|\bfH{h, \bfG{h}}(x)| 
	\leq K \bar h(x) .
$
By \eqref{eq:cQ:prop3} 
there exists constants $K_1,K_2\in (0,\infty)$ (all depending on $\sigma_\pm$ and $K_2$ depending also on the constant $a\in [0,\infty)$) such that
\begin{equation}
\begin{split}
	|\bfH{h, \cQ_{n,1,0}}(x)| 
	& 
	\leq K_1 \bar h(x) \left(  \int_{-\infty}^\infty  (1+|x| \ind{[-x,x]}(y)+ |y| \ind{\RR\setminus[-x,x]}(y)) e^{a |y-x|} q_\sigma (1,x,y) \vd y \right)
	\\
	& 
	\leq K_2 \bar h(x) (1+|x|). 	 
\end{split}
\end{equation} 
The fact that $h\in \Lg{\gamma}\subseteq \Lg{2}$ ensures that $\bar h \in \Lgb{1}$ (by Lemma~\ref{lem:link:theorems}) and so $\bfH{h, \cQ_{n,0,0}} 
\in \Lone$. 
Hence, dominated convergence and Lemma~\ref{lem:cQ:conv} show that 
$
	\lim_{n\to\infty}\langle \lambda_\sigma , f_n \rangle=\langle \lambda_\sigma , \bfH {h^2} + 2 \bfH{h, \mathcal{P}_{h}^{Y} } \rangle
$.

Let us now show that $f_n$ satisfies equation \eqref{eq:j1:4.2}.
By~Lemma~\ref{lem:link:theorems},
$\bfH {h^2} \in \Lgb{2}$. 
%
\\
Let us explore the contribution to \eqref{eq:j1:4.2} of the other parts of $f_n$. 
Let us first consider 
$
	\bfH{h, \cQ_{n,0,0}}= \bfH{h, \bfG{h}}+\bfH{h, \cQ_{n,1,0}}.
$
Above we saw that $|\bfH{h, \bfG{h}}|\leq K\bar h$ with $K$ non negative constant. 
Cauchy-Schwarz inequality implies that
\begin{equation}
	( \bfH{h, \cQ_{n,1,0}}(x) )^2
	\leq  \int_{-\infty}^\infty h(x,y)^2 q_\sigma (1,x,y) \vd y \int_{-\infty}^\infty  (\cQ_{n,1,0}(y))^2 q_\sigma (1,x,y) \vd y.
\end{equation}
It can be easily shown that the first factor is uniformly bounded by a finite constant.
And by \eqref{eq:cQ:prop2} (taking $\zeta = \gamma-1>2$) there exist constants $K_1,K_2\in (0,\infty)$ 
such that
\begin{equation} \label{eq:recicle} 
	\begin{split}
	Q^\sigma_1 \left(  (\cQ_{n,1,0})^2 \right) & = \int_{-\infty}^{\infty} (\cQ_{n,1,0}(y))^2 q_\sigma (1,x,y) \vd y
	\\
	& \leq  K_1
	 \int_{-\infty}^{\infty}   \frac{(\log n)^2  q_\sigma (1,x,y)}{(1+|y n^{-\frac18}/\sigma(y)|^{\gamma-1})^2}  \vd y
	\leq  \frac{ K_2 (\log (n))^2}{1+|x n^{-\frac18}/\sigma(x)|^{2(\gamma-1)}}.
\end{split}
\end{equation}
The last inequality are consequences of the upper bound for the transition density of OBM
$
	q_\sigma(1,x,y) \leq \frac{\left(1+\left|\frac{\sigma_--\sigma_+}{\sigma_-+\sigma_+}\right|\right)}{\sqrt{2\pi}\sigma(y)} e^{-\frac12(\frac{y}{\sigma(y)}-\frac{x}{\sigma(x)})^2}
$
and of \cite[Lemma~3.2]{j1} (or some computations).
Therefore 
$
	\bfH{h, \cQ_{n,0,0}}(x) \leq K \left(\bar h(x) + \frac{ \log (n)}{1+|x n^{-\frac18}/\sigma(x)|^{\gamma-1}} \right)\!.$
Finally we consider the auxiliary function $g_n$: note that
$
	|g_n|(x) \leq Q^\sigma_1 (  2 (\cQ_{n,0,0})^2 + 2 (\bfG{h})^2)(x) + 2 (\cQ_{n,1,1}(x))^2 +  2(\bfG{h}(x))^2.
$
Inequality~\eqref{eq:cQ:prop2} and inequality~\eqref{eq:recicle} imply that
$Q^\sigma_1 \left(  (\cQ_{n,0,0})^2 \right) (x)+ (\cQ_{n,1,1}(x))^2 \leq \frac{ K (\log (n))^2}{1+|x n^{-\frac18}/\sigma(x)|^{2(\gamma-1)}}$
for some non negative constant $K$.
Lemma~\ref{lem:prop:kappa}, the fact that $\gamma \geq 2$ and Item~\ref{item2bis:semig:estim} of Lemma~\ref{lem:semig:estim} yields $(\bfG{h})^2\leq K  \bfG{h}$ and
$
	Q^\sigma_1 \left((\bfG{h})^2\right) (x) \leq K (e^{-x^2/2} + \frac{1}{1+({|x|}/{\sigma(x)})^\gamma})
$
for some non-negative constant $K$.
%
%
%
Combining all terms yields that $f_n$ satisfies~\eqref{eq:j1:4.2}. (Here we used the fact that $\gamma \geq 3$ and $\eta=1/4<1/3$ in the definition of $\cQ$.)

{\em Fourth step:\/} 
We show that
$
	\frac1{\sqrt n}  \sum_{k=1}^{\lfloor n t \rfloor} h_n(\sqrt{n} Y_{(k-1)/n}) \convp[n\to\infty]
	c_h L_t(Y)
$
applying Proposition~\ref{prop:th4.1b:j1},
where 
$ 2 \langle \lambda_\sigma , \bfH h \rangle c_h:= - K_h + 
	\langle \lambda_\sigma , \bfH {h^2} + 2 \bfH{h,\mathcal{P}_{h}^{Y} } \rangle 
	+ \frac{ 2 \sigma_- \sigma_+}{\sigma_-+\sigma_+} \frac{8}{3 \sqrt{2\pi}} (\langle \lambda_\sigma , \bfH h \rangle)^2 $.
To do so we check that the sequence $h_n$ satisfies \eqref{eq:j1:4.2} and that 
$
	\lim_{n\to\infty} \langle \lambda_\sigma , h_n \rangle=c_h
$.

Inequality~\eqref{eq:cQ:prop3}, 
the fact that $\bfG{h}$ is bounded (see Lemma~\ref{lem:prop:kappa}), 
the fact that
\begin{equation}
	\EE\!\left[ L_1(W)  (1 + |W_1|+ e^{a \sigma(W_1) |W_1|}  \ind{\{\sigma(W_1)\in \RR\}}) \right]<\infty
\end{equation}
(see the joint density of BM and its local time~\eqref{eq:BMl}),
boundedness of $\bar h$, and the change of variable $s= \frac{x^2}{ (\sigma(x))^2 t}$ show that there exists $K\in (0,\infty)$ such that 
\begin{equation}
\begin{split}
	& |\Lm{h(x,\cdot)}{1}{x}| + |\Lm{\cQ_{n,0,0}}{1}{x}| 
	\leq 
	2 K |x| e^{-\frac{x^2}{(\sigma(x))^2}} \in \Lone.
\end{split}
\end{equation}
Hence dominated convergence, and Lemma~\ref{lem:cQ:conv} demonstrates that 
$\lim_{n\to\infty} \langle \lambda_\sigma , \Lm{\cQ_{n,0,0}}{1}{\cdot} \rangle 
= \langle \lambda_\sigma , \Lm{ \mathcal{P}_{h}^{Y} }{1}{\cdot} \rangle$.
Moreover 
the latter inequalities ensure also that 
$h_n$ 
satisfy \eqref{eq:j1:4.2}.

Lemma~\ref{lem:prop:G} allows us to rewrite
$c_h:= \langle \lambda_\sigma , \Lm{ h(\cdot,) }{1}{\cdot} \rangle + \langle \lambda_\sigma , \Lm{ \mathcal{P}_{h}^{Y} }{1}{\cdot} \rangle$ as
\begin{equation}\begin{split}	& c_h = \sqrt{\frac{2}{\pi}} \frac{(\sigma_- + \sigma_+)}{2\sigma_-\sigma_+}
			 \int_0^1  \sqrt{\frac{1}t-1}  \int_{-\infty}^\infty \int_0^\infty  \rho_1^\sigma(y,\ell) \ell \mathcal{P}_{h}^{Y} (y \sqrt{1-t}) \vd \ell \vd y \vd t
	\\
	& \qquad +
\int_{-\infty}^\infty \int_0^1 \int_{-\infty}^\infty \int_0^\infty 
\frac{|x| e^{-\frac{x^2}{2 (\sigma(x))^2}}}{\sqrt{2\pi}\sigma(x)}  \sqrt{\frac{1}t-1}  \rho_1^\sigma(y,\ell) \ell h(x \sqrt{t},y \sqrt{1-t})  
\vd \ell \vd y \vd t \lambda_\sigma(\!\vd x)
\end{split}\end{equation}
where 
$\rho_1^\sigma$ is the joint density of a standard OBM and its local time at time $1$ (given in \eqref{eq:OBMl}).
This expression can be easily checked to be 
the desired expression for $c_h$.

\subsubsection{Proof of Item~\ref{th:j2:item3:2} in Proposition~\ref{th:j2:2}}
Let $t\in [0,1]$ be fixed.

{\it First step:}
Let us show that 
\begin{equation}
\begin{split}
	& 
\EE\!\left[\chi_k^n (Y_{k/n} - Y_{(k-1)/n})  | \cF_{(k-1)/n}\right] 
	\\
	& = \frac1{\sqrt n}
\left(  n^{-\frac14} f_1(\sqrt{n} Y_{(k-1)/n}) + g_{n}(\sqrt{n} Y_{(k-1)/n}) + \langle \lambda_\sigma , \bfH h \rangle   n^{-\frac14} f_2(\sqrt{n} Y_{(k-1)/n})\right)
\end{split}
\end{equation}
where $f_1,f_2,g_n$ are given by 
\begin{equation}
\begin{split}
	f_1(x)& := \int_{-\infty}^\infty (h(x,y)+\bfG{h}(y)) (y-x) q_\sigma(1,x,y) \vd y , 
	\ f_2(x):= x \Lm{1}{1}{x}, \ \text{and}
	\\
	g_n(x) & := n^{-\frac14} \int_{-\infty}^\infty \cQ_{n,1,0}(y) (y-x) q_\sigma(1,x,y) \vd y.
\end{split}	
\end{equation}

Throughout the proof of this step let $\mathbb{Id} \colon x\mapsto x$ denote the identity function.
It follows from \eqref{eq:chi:nk}, the fact that $Y$ is a martingale, \eqref{function:G:equal},
and the scaling property~\eqref{Yscaling}
that
\begin{equation}
\begin{split}
	& n^{\frac34} \, \EE\!\left[\chi_k^n (Y_{k/n} - Y_{(k-1)/n})  | \cF_{(k-1)/n}\right] 
	\\
	& =  f_1(\sqrt{n} Y_{(k-1)/n}) +  \langle \lambda_\sigma , \bfH h \rangle f_2(\sqrt{n} Y_{(k-1)/n})
	- \Lm{ \mathbb{Id} }{1}{\sqrt{n} Y_{(k-1)/n}} )
	+ n^{1/4}g_n(\sqrt{n}Y_{(k-1)/n}).
\end{split}
\end{equation}
The proof of this step is completed if
$ \Lm{ \mathbb{Id} }{1}{\cdot} =0$.
This equality follows from the fact that for a standard BM, say $W$, it holds that $\EE\!\left[ L_1(W)  W_1  \right]=0$
and from Lemma~\ref{lem:prop:G}. 

In the next steps we want to check that Proposition~\ref{prop:th4.1a:j1} can be applied to the sequences $n^{-\frac14} f_1$, $n^{-\frac14} f_2$ and $g_n$. 

{\it Second step:} We show that $f_1$ is bounded and integrable: $f_1\in \Lgb{0}$.
  
The proof follows from Lemma~\ref{lem:prop:kappa}, and Item~\ref{item2bis:semig:estim} in Lemma~\ref{lem:semig:estim}. In this proof the fact that $\gamma>3$ is strongly used.

The fact that $h\in \Lg{\gamma}\subseteq \Lg{3}$ and Lemma~\ref{lem:prop:kappa} ensure that 
$\bar h$ and $\bfG{h}$ are in $\Lgb{2}$.
This, Cauchy-Schwarz, the fact that 
$
	\left(\int_{-\infty}^{\infty} (y-x)^2 q_\sigma(1,x,y) \vd y\right) 
	\leq \sigma_-^2 \ind{\{\sigma_-\in\RR\}} +\sigma_+^2 \ind{\{\sigma_+\in\RR\}}
$,
and the fact that that $\int_{-\infty}^{\infty}  (\bfG{h}(y))^2 q_\sigma(1,x,y) \vd y = Q^\sigma_1 \bfG{h}^2 (x)$
yield
\begin{equation}
\begin{split} 
	|f_1(x)|  
	&
	\leq \bar h(x) \int_{-\infty}^{\infty}  e^{a |y-x|} |y-x| q_\sigma(1,x,y) \vd y
	+ \sqrt{\sigma_-^2 \ind{\{\sigma_-\in\RR\}} +\sigma_+^2 \ind{\{\sigma_+\in\RR\}}} \left( |Q^\sigma_1 \bfG{h}^2 (x)| \right)^{\!1/2}.
\end{split}
\end{equation}
By Lemma~\ref{lem:prop:kappa} $\bfG{h}$ is bounded, i.e.~there exists $K\in (0,\infty)$ such that 
$\sup_{x\in \RR}|\bfG{h}(x)|\leq K$ and $\langle \lambda_\sigma , \bfG{h} \rangle=0 $, and so
$
	| \langle \lambda_\sigma , \bfG{h}^2 \rangle | \leq K  | \langle \lambda_\sigma , \bfG{h} \rangle | =0.
$
And so Item~\ref{item2bis:semig:estim} in Lemma~\ref{lem:semig:estim} implies, up to increase the constant $K\in [1,\infty)$, 
that
\begin{equation}
	|Q^\sigma_1 \bfG{h}^2 (x)| \leq K \frac{1}{1+|x/\sigma(x)|^{\gamma-1}} + K \frac1{\sqrt{2 \pi}} e^{-\frac{x^2}{2 (\sigma(x))^2 }} 
\leq 2 K^2 \frac{1}{(1+|x/\sigma(x)|^{\frac{\gamma-1}{2}})^2} .
\end{equation}
Since $\gamma >3 $ it holds that $\left( |Q^\sigma_1 \bfG{h}^2 (x)| \right)^{\!\frac12}\in \Lgb{0}$ and so $|f_1|\in \Lgb{0}$.

{\it Third step:} The fact that $f_2$ is bounded and integrable follows from Lemma~\ref{lem:prop:G}.

This and the change of variable $s \mapsto r= \frac{x^2}{(\sigma(x))^2 s}$ show that there exists a positive constant $K$ 
such that
\begin{equation}
\begin{split} 
	f_2(x)  
	& = \left(\frac{2 \sigma_-\sigma_+}{\sigma_-+\sigma_+}\right)^{\! 2} \frac{x |x|}{\sigma(x)} \EE\!\left[ \frac{L_1(W)}{\sigma(W_1)} \right] 
	\int_0^1  \frac{\sqrt{1-s}}{\sqrt{2\pi} s^\frac32} e^{-\frac{x^2}{2 (\sigma(x))^2 s}} 
	\! \vd s
	\\
	& = \left(\frac{2 \sigma_-\sigma_+}{\sigma_-+\sigma_+}\right)^{\! 2} \frac{x |x|}{\sigma(x)} \EE\!\left[ \frac{L_1(W)}{\sigma(W_1)} \right] 
	\int_{\frac{x^2}{(\sigma(x))^2}}^\infty  \frac1{\sqrt{2\pi}}
 \sqrt{\frac{(\sigma(x))^2}{x^2}-\frac1r} e^{-\frac{r}{2}} 
	\frac1{\sqrt{r}}\! \vd r \ \ind{\{\sigma(x)\in \RR\}}
	\\
	& \leq \left(\frac{2 \sigma_-\sigma_+}{\sigma_-+\sigma_+}\right)^{\! 2} \frac{|x| \EE\!\left[L_1(W)\right] }{\min\{\sigma_-,\sigma_+\}}
	\int_{\frac{x^2}{(\sigma(x))^2}}^\infty  \frac1{\sqrt{2\pi}}
  e^{-\frac{r}{2}} 
	\frac1{\sqrt{r}}\! \vd r \ \ind{\{\sigma(x)\in \RR\}}
	\leq K e^{-\frac{x^2}{4(\sigma(x))^2}} \in \Lgb{0}.
\end{split}
\end{equation}
In the last inequality we used that $\EE\!\left[ L_1(W) \right] \in [0,\infty)$.
This step is thus completed.

{\it Fourth step:} We prove that $\int_{-\infty}^{\infty} |g_n(x)| \vd x$ 
and 
$\frac1{\sqrt{n}} g_n(\sqrt{n} x)$ converge to 0.

Cauchy-Schwarz inequality 
ensures
\begin{equation}
\begin{split}
	|g_n(x)|^2 
	& \leq n^{-1/2} \left( \int_{-\infty}^{\infty} (\cQ_{n,1,0}(y))^2 q_\sigma (1,x,y) \vd y\right) \left(\int_{-\infty}^{\infty} (y-x)^2 q_\sigma(1,x,y) \vd y\right).
\end{split}
\end{equation}
Note that
$
	\left(\int_{-\infty}^{\infty} (y-x)^2 q_\sigma(1,x,y) \vd y\right) 
\leq \sigma_-^2 \ind{\{\sigma_-\in \RR\}} +\sigma_+^2 \ind{\{\sigma_+\in \RR\}}
$ 
and inequality~\eqref{eq:recicle}
yield that there exists a constant $K\in (0,\infty)$ depending on $\gamma$ and $\sigma_\pm$ such that
$
 |g_n(x)|
	\leq  \frac{n^{-\frac14} K_2 \log (n)}{1+|x n^{-\frac18}/\sigma(x)|^{\gamma-1}}.
$

In the last steps we proved that Proposition~\ref{prop:th4.1a:j1} can be applied and this completes the proof.

\subsubsection{Proof of Item~\ref{th:j2:item4:2} in Proposition~\ref{th:j2:2}}
Let $\varepsilon \in (0,\infty)$ be fixed. 
For every $k$, H\"older's inequality and Markov's inequality show
\begin{equation}
\begin{split}
	& \EE\!\left[ |\chi_k^n|^2 \ind{\{|\chi_k^n|\geq  \varepsilon \}}| \cF_{(k-1)/n}\right] 
	 \leq 
\EE\!\left[ |\chi_k^n|^{5} | \cF_{(k-1)/n}\right]  \varepsilon^{-3} .
\end{split}
\end{equation}
The fact that $h\in \Lg{\gamma}\subseteq \Lg{0}$ ensures that $\bar h$ is bounded and integrable. This combined with Jensen's inequality and \eqref{eq:cQ:prop4} 
ensures that $\sup_{x\in \RR} \bar h(x) +  |\langle \lambda_\sigma , \bfH h \rangle| $ is bounded by a constant $K\in (0,\infty)$ and for all $n\in \NN\setminus\{0,1,2,3\}$ it holds that
\begin{equation}
\begin{split}
	 & \EE\!\left[ |\chi_k^n|^{5} | \cF_{(k-1)/n}\right] 
	  \leq  4^{4} n^{-\frac{5}{4}} K^{5} \Big(  \EE\!\left[ e^{5 a \sqrt n | Y_{k/n} - Y_{(k-1)/n}|}  | \cF_{(k-1)/n}\right]  
	\\
	& \quad \qquad 	+ \EE\!\left[  \left(\sqrt{n}|L_{k/n}(Y) - L_{(k-1)/n}(Y)| \right)^{5} | \cF_{(k-1)/n}\right] 
	+ (\log n)^{5} \Big).
\end{split}
\end{equation}
The fact that density of the OBM has a Gaussian behavior together with the scaling property~\eqref{Yscaling} ensure that
$ 
	\EE\!\left[ e^{5 a \sqrt n | Y_{k/n} - Y_{(k-1)/n}|}  | \cF_{(k-1)/n}\right]  = \int_{-\infty}^{\infty} e^{5 a |y-x|} q_\sigma(1,x,y) \vd y
$
is bounded.
By \eqref{function:G} and Lemma~\ref{lem:prop:G} we can show, similarly to the third step of the proof of Item~\ref{th:j2:item3}, that there exist constants $K_1,K_2 \in (0,\infty)$ such that 
\begin{equation}
\begin{split}
	& \EE\!\left[  \left(\sqrt{n}|L_{k/n}(Y) - L_{(k-1)/n}(Y)| \right)^{5} | \cF_{(k-1)/n}\right] 
	= \Lm {1}{5}{\sqrt{n}Y_{(k-1)/n}} 
	\\
	& \leq K_1 \int_0^1 \frac{(1-t)^{\frac52}}{\sqrt{2 \pi} t^{\frac32}}  \frac{\sqrt{n}|Y_{(k-1)/n}|}{\sigma(Y_{(k-1)/n})} e^{-\frac{n (Y_{(k-1)/n})^2}{2 (\sigma(Y_{(k-1)/n})) t } } \vd t
\leq K_2.
\end{split}
\end{equation}
We conclude that there exists a constant $K\in (0,\infty)$ such that
\begin{equation}
\begin{split}
	& \sum_{k=1}^n \EE\!\left[ |\chi_k^n|^2 \ind{\{|\chi_k^n|\geq  \varepsilon \}}| \cF_{(k-1)/n}\right] 
	\leq n K n^{-\frac54} \log(n)  \varepsilon^{-3} \xrightarrow[n\to\infty]{} 0.
\end{split}
\end{equation}

\noindent \textbf{Acknowledgments:} 
The author thanks Sylvie R{\oe}lly and Paolo Pigato for fruitful discussions on the topic of this paper and 
Andrey Pilipenko for the reference~\cite{port1}. 
The author is grateful to the anonymous referees of different reviewing processes for their valuables suggestions which pushed her to improve the document, for instance the presentation of the technical proofs.

\addcontentsline{toc}{section}{Bibliography}
\bibliographystyle{abbrvnat}
\bibliography{../Bernoulli/2024_submitted_version/Bibliography}


\end{document}